\documentclass{memo-l}

\makeatletter

\numberwithin{section}{chapter} 
\theoremstyle{plain}
\newtheorem{thmc}{Theorem}[chapter]
\numberwithin{equation}{section} 
\numberwithin{figure}{section} 
\theoremstyle{plain}
\newtheorem*{thmc*}{Theorem}
\theoremstyle{plain}
\theoremstyle{plain}
\newtheorem{lemc}[thmc]{Lemma} 
\theoremstyle{plain}
\newtheorem{propc}[thmc]{Proposition} 
\theoremstyle{definition}
\newtheorem{defnc}[thmc]{Definition}
\theoremstyle{remark}

\theoremstyle{remark}

\theoremstyle{remark}

\theoremstyle{remark}

\theoremstyle{definition}

\theoremstyle{remark}

\theoremstyle{plain}
\newtheorem{thm}{Theorem}[section]
\numberwithin{equation}{section} 
\numberwithin{figure}{section} 
\theoremstyle{plain}
\newtheorem*{thm*}{Theorem}
\theoremstyle{plain}
\newtheorem{cor}[thm]{Corollary} 
\theoremstyle{plain}
\newtheorem{lem}[thm]{Lemma} 
\theoremstyle{plain}
\newtheorem{prop}[thm]{Proposition} 
\theoremstyle{definition}
\newtheorem{defn}[thm]{Definition}
\theoremstyle{remark}
\newtheorem{rem}[thm]{Remark}
\theoremstyle{remark}

\theoremstyle{remark}

\theoremstyle{remark}

\theoremstyle{definition}

\theoremstyle{remark}
  \newtheorem*{acknowledgement*}{Acknowledgement}

\theoremstyle{plain}

\theoremstyle{plain}
\theoremstyle{plain}
\theoremstyle{plain}
\theoremstyle{definition}

\theoremstyle{remark}

\theoremstyle{remark}

\theoremstyle{remark}

\theoremstyle{plain}

\newcommand{\id}{\operatorname{id}}

\makeatletter
\def\freeprodsize@{\@setfontsize\freeprodsize\@xxpt\@xxpt}
\def\freeprod@{\mathop{\hbox{\freeprodsize@ $*$}}}
\def\freeprod{\freeprod@\displaylimits}
\makeatother
\newcommand{\bigast}{\freeprod}

\newcommand{\TA}{\operatorname{{\mathrm T}(A)}}
\newcommand{\TAim}{\operatorname{{\mathrm AT}(A)}}
\newcommand{\TAuim}{\operatorname{{\mathrm UAT}(A)}}
\newcommand{\TAQD}{\operatorname{{\mathrm AT}(A)_{QD}}}
\newcommand{\TAuQD}{\operatorname{{\mathrm UAT}(A)_{QD}}}
\newcommand{\TAlfd}{\operatorname{{\mathrm AT}(A)_{LFD}}}
\newcommand{\TAulfd}{\operatorname{{\mathrm UAT}(A)_{LFD}}}

\newcommand{\TB}{\operatorname{{\mathrm T}(B)}}
\newcommand{\TBim}{\operatorname{{\mathrm AT}(B)}}
\newcommand{\TBuim}{\operatorname{{\mathrm UAT}(B)}}
\newcommand{\TBQD}{\operatorname{{\mathrm AT}(B)_{QD}}}
\newcommand{\TBuQD}{\operatorname{{\mathrm UAT}(B)_{QD}}}
\newcommand{\TBlfd}{\operatorname{{\mathrm AT}(B)_{LFD}}}
\newcommand{\TBulfd}{\operatorname{{\mathrm UAT}(B)_{LFD}}}

\usepackage{amssymb}

\makeatother
 
\begin{document}

\title[Finite Representation Theory]{Invariant Means and Finite
Representation Theory of C$^*$-algebras}

\author{Nathanial P. Brown}

\address{Department of Mathematics, Penn State University, State College, PA 16802}

\email{nbrown@math.psu.edu}

\renewcommand{\thepage}{\roman{page}}

\thanks{This research was supported by MSRI and NSF Postdoctoral Fellowships.}
\date{May 13, 2004.}


\begin{abstract}
Various subsets of the tracial state space of a unital C$^*$-algebra
are studied.  The largest of these subsets has a natural
interpretation as the space of invariant means.  II$_1$-factor
representations of a class of C$^*$-algebras considered by Sorin Popa
are also studied. These algebras are shown to have an unexpected
variety of II$_1$-factor representations.  In addition to developing
some general theory we also show that these ideas are related to
numerous other problems in operator algebras.
\end{abstract}
\maketitle
\cleardoublepage
  \thispagestyle{empty}
  \vspace*{13.5pc}
  \begin{center}
  To my big, beautiful family -- on both sides of the Pacific.
  \end{center}
  \cleardoublepage

\setcounter{page}{7}
 \tableofcontents

\chapter{Introduction}
\renewcommand{\thepage}{\arabic{page}}
\setcounter{page}{1}

One of von Neumann's motivations for initiating the study of operator
algebras was to provide an abstract framework for unitary
representation theory of locally compact groups.  Hence it is no
surprise that representation theory of C$^*$-algebras attracted the
attention of many experts over the years. However, Glimm's deep work,
essentially closing the book on global (irreducible) representation
theory, combined with the emergence of exciting new fields such as
Connes' noncommutative geometry, Jones' theory of subfactors,
Elliott's classification program and, most recently, Voiculescu's
theory of free probability, have kept representation theory out of the
limelight for the last couple decades.

In this paper we revisit representation theory with the goal of
convincing the reader that (a) {\em finite} representation theory is
an important subject where much work is still needed and (b)
future advances may well provide the key to unlocking
some important open problems.  Indeed, in addition to developing some
basic general theory, these notes contain the solutions to several
problems, give streamlined proofs of (generalizations of) some known
results and give  new insight into other problems.  In
contrast to irreducible representations or factor representations of
type III (e.g.\ the celebrated work of Powers) it seems that finite
representations have not received the attention they deserve.  We
hope to lay the groundwork for future study in this direction.

By finite representation theory we mean GNS representations arising
from tracial states. It turns out that some traces are better than
others and the `good' ones are precisely those which can be
interpreted as invariant means.

\begin{defnc}
\label{thm:invariantmean} Let $A \subset B(H)$ be a concretely
represented, unital C$^*$-algebra.  A state $\tau$ on $A$ is
called an {\em amenable trace} if there exists a state $\phi$ on
$B(H)$ such that (1) $\phi|_A = \tau$ and (2) $\phi(uTu^*) =
\phi(T)$ for every unitary $u \in A$ and $T \in B(H)$.
\end{defnc}

Amenable traces really are
traces (i.e.\ $\tau(ab) = \tau(ba)$ for all $a,b \in A$). In more
classical language, $\tau$ is the restriction of a state which
contains $A$ in its centralizer and it is well known that this
procedure yields a tracial state.  In case one hasn't seen this
argument, first note that for every $a \in A$ and {\em unitary} $u \in
A$ we have $$\tau(au) = \tau(u(au)u^*) = \tau(ua).$$ Since every
element $b \in A$ is a linear combination of unitaries, it follows that 
$\tau(ab) = \tau(ba)$ for all $a,b \in A$. 

Next we observe that these traces are analogous to invariant means on
groups.  Recall that amenable groups are, by definition, those which
have an invariant mean -- i.e.\ there exists a state on
$L^{\infty}(G)$ which is invariant under the left translation action
of $G$ on $L^{\infty}(G)$.  If $A \subset B(H)$ is a C$^*$-algebra
then there is a natural action of the unitary group of $A$ on $B(H)$
given by $T \mapsto uTu^*$, where $T \in B(H)$ is arbitrary and $u \in
A$ is unitary.  Hence condition (2) in the definition above is
precisely the statement that {\em $\phi$ is a state on $B(H)$ which is
invariant under the action of the unitary group of $A$}.  In other
words, $\tau$ is an amenable trace if and only if it is the
restriction of an invariant mean on $B(H)$.\footnote{At this point one
may wonder if this notion depends on the particular choice of
representation $A \subset B(H)$. It doesn't, as we will observe later.}

This notion is well known (under different names) and has already
shown its importance (cf.\ \cite{connes:classification},
\cite{connes:compactmetricspaces}, \cite{kirchberg:invent},
\cite{kirchberg:propertyTgroups}, \cite{bekka}, \cite{popa:simpleQD},
\cite{bedos:hypertraces}).  Indeed, in Connes' remarkable paper
\cite{connes:classification} he showed (using the terminology
`hypertrace') that a II$_1$-factor is isomorphic to the hyperfinite
II$_1$-factor if and only if its unique trace is amenable.  Continuing the
philosophy that this notion should correspond to some kind of
amenability, Bekka defined a unitary representation of a locally
compact group to be amenable if the C$^*$-algebra generated by the
image of the representation has an amenable trace (cf.\ \cite{bekka}).
The most important result concerning the structure of amenable traces,
however, is due to Kirchberg who (building on work of Connes and using
the terminology `liftable') showed that amenable traces are precisely
those which enjoy a natural finite dimensional approximation
property (cf.\ \cite{kirchberg:propertyTgroups}). This is reminiscent
of the approximation properties which appear in the theories of
nuclear, exact or quasidiagonal C$^*$-algebras and, just as for
operator algebras, it proves very useful.

In addition to the study of amenable traces and their
corresponding GNS representations, we will also investigate
II$_1$-factor representations of Popa algebras.  These algebras
are defined via an internal finite dimensional approximation
property (see Definition \ref{thm:popaalgebradefn}) which is the
C$^*$-analogue of an approximation property which characterizes
the hyperfinite II$_1$-factor $R$.  In fact, the definition of a
Popa algebra appears to be so close to the characterization of $R$
that it was believed for some time that there could only be one
II$_1$-factor which arose from a representation of a Popa algebra
(namely $R$).  We will see that this is not the case as the finite
representation theory of these algebras is quite rich (containing
all McDuff factors, for example).

While we believe that both amenable traces and factor
representations of Popa algebras are topics of independent
interest, perhaps the most surprising part of this work is that
combining these two aspects of representation theory leads to a
variety of new results which provide a common thread between
several important problems in operator algebras.  Our results are
most strongly connected to the classification program and
questions around free probability but there are also relations
with geometric group theory (see Proposition
\ref{thm:residuallyfinite} and questions (8) and (9) in Section
\ref{thm:questions}), theoretical numerical analysis and operator
theory (see Section \ref{thm:numericalanalysis}), K-homology
(see Section \ref{thm:k-homology}) and the C$^*$-algebraic
structure of the hyperfinite II$_1$-factor (see Section
\ref{thm:stablyfinitevsQD}).

Since we deal with a broad range of topics, it may be worthwhile 
to give a detailed overview before proving any results. 

\vspace{2mm}
\noindent{1.1. \bf Approximating Traces on C$^*$-algebras}

As mentioned above, amenable traces are precisely those which
arise from a natural finite dimensional approximation property.
Strengthening this property in various ways leads to natural
definitions of subspaces of the space of amenable traces (see
Sections \ref{thm:uatraces}, \ref{thm:qdtraces} and
\ref{thm:lfdtraces}). Studying these stronger approximation
properties is roughly the analogue of passing from the class of
nuclear C$^*$-algebras to a sub-class such as nuclear,
quasidiagonal C$^*$-algebras or homogeneous C$^*$-algebras.  One
subspace of the amenable traces precisely characterizes those GNS
representations which have hyperfinite von Neumann algebras after
closing in a weak topology (see Theorem \ref{thm:mainthmUTAwafd}).
The relation between the other subspaces and representation theory
is not yet clear. However, we will see that these other subspaces
play a critical role in Elliott's classification program (see
Section \ref{thm:classificationprogram}).

As mentioned above, basic facts about amenable traces have a
number of consequences. Sections \ref{thm:numericalanalysis}
and \ref{thm:k-homology} work out some applications to the finite
section method from numerical analysis and some basic
K-homological questions, respectively. We also observe that these
traces provide natural obstructions to the existence of (unital)
$*$-homomorphisms between certain classes of operator algebras
(cf.\ Corollary \ref{thm:nohomomorphisms}).

\vspace{2mm}
\noindent{1.2. \bf II$_1$-Factor Representations of Popa Algebras}

Another goal of these notes is to study II$_1$-factor representations
of Popa algebras.

\begin{defnc}
\label{thm:popaalgebradefn}
A simple, separable, unital C$^*$-algebra, $A$, is called a {\em Popa
algebra} if for every finite subset $\mathfrak{F} \subset A$ and
$\varepsilon > 0$ there exists a nonzero finite dimensional
C$^*$-subalgebra $B \subset A$ with unit $e$ such that $\| ex - xe \|
\leq \varepsilon$ for all $x \in \mathfrak{F}$ and $e\mathfrak{F}e
\subset^{\varepsilon} B$ (i.e.\ for each $x \in \mathfrak{F}$ there
exists $b \in B$ such that $\| exe - b \| \leq \varepsilon$).
\end{defnc}

Popa algebras are always quasidiagonal. Using his local quantization
technique in the C$^*$-algebra setting Popa nearly provides a converse
in \cite{popa:simpleQD}: Every simple, unital, quasidiagonal
C$^*$-algebra with `sufficiently many projections' (e.g.\ real rank
zero) is a Popa algebra.  Thus the class of Popa algebras is much
larger than one might first guess.

For some time there was speculation that quasidiagonality may be
closely related to nuclearity.  For example, in \cite[pg.\
157]{popa:simpleQD} Popa asked whether every Popa algebra with {\em
unique} trace is necessarily nuclear. (Counterexamples were first
constructed by Dadarlat in \cite{dadarlat:nonnuclearsubalgebras}.)
More generally, he asked in \cite[Remark 3.4.2]{popa:simpleQD} whether
the hyperfinite II$_1$-factor $R$ was the {\em only} II$_1$-factor
which could arise from a GNS representation of a Popa algebra.
Support for a positive answer was provided by the following
observation of Popa (cf.\ \cite{popa:injectiveimplieshyperfiniteI},
\cite{popa:injectiveimplieshyperfiniteII}):

\begin{thm*}
Let $M$ be a separable II$_1$-factor.  Then $M \cong R$ if and only if
for every finite subset $\mathfrak{F} \subset M$ and $\varepsilon > 0$
there exists a nonzero finite dimensional C$^*$-subalgebra $B \subset
M$ with unit $e$ such that $\| ex - xe \|_2 \leq \varepsilon \| e
\|_2$ for all $x \in \mathfrak{F}$ and for each $x \in \mathfrak{F}$
there exists $b \in B$ such that $\| exe - b \|_2 \leq \varepsilon
\|e\|_2$, where $\| \cdot \|_2$ is the 2-norm on $M$ coming from the
unique trace.
\end{thm*}

Note that the definition of a Popa algebra is the C$^*$-analogue
of the approximation property above which characterizes $R$.
Moreover, the 2-norm version appears, at first glance, to allow
one to check the approximation property on a weakly dense
subalgebra.  In other words, it was quite natural to expect that
if a II$_1$-factor $M$ contained a weakly dense Popa algebra then
it would have to be isomorphic to $R$. Indeed, suppose $A \subset
M$ is a weakly dense Popa algebra, $\mathfrak{F} \subset A$, a
finite set in the unit ball of $A$, and $\varepsilon > 0$ are
given.  By definition we can find a nonzero finite dimensional
C$^*$-subalgebra $B \subset A$ with unit $e$ such that $\| ex - xe
\| \leq \varepsilon$ for all $x \in \mathfrak{F}$ and
$e\mathfrak{F}e \subset^{\varepsilon} B$.  From the general
inequality $\|ab\|_2 \leq \|a\| \|b\|_2$ it follows that $$\| ex -
xe \|_2 \leq \|ex - exe\|_2 + \|exe - xe\|_2 = \| e(ex - xe)\|_2 +
\|(ex - xe)e\|_2 \leq 2 \varepsilon \| e \|_2$$ for all $x \in
\mathfrak{F}$ and for each $x \in \mathfrak{F}$ there exists $b
\in B$ such that $$\| exe - b \|_2 = \| e(exe - b)\|_2 \leq
\varepsilon \|e\|_2.$$ Hence $M$ has the desired approximation
on a weakly dense subalgebra and one ``should'' be able
to pass from $A$ to all of $M$ via some additional 
argument (which would imply $M\cong R$).  Since Popa algebras are
simple all their GNS representations are faithful and hence it was
expected that $R$ would be the only II$_1$-factor which could arise 
from a GNS representation.

It turns out, however, that this is not the case (i.e.\ it is
impossible to use an approximation argument to pass from $A$ to
all of $M$).  For example, we will construct a Popa algebra $A$
with the property that for every (separable) II$_1$-factor $M$
there exists a tracial state $\tau$ on $A$ such that $\pi_{\tau}
(A)^{\prime\prime} \cong M \bar{\otimes} R$ (see Theorem
\ref{thm:arbitraryMcDuff}).  Popa has also asked (private
communication) if a Popa algebra with unique trace must
necessarily yield the hyperfinite II$_1$-factor.  We will see that
this question also has a negative answer: There exists a Popa
algebra with Dixmier property (hence unique trace) such that the
GNS construction yields a
non-hyperfinite II$_1$-factor (see Theorem
\ref{thm:uniquetracenotTAF}). On the other hand we will show that
if $A$ is a locally reflexive (e.g.\ exact) Popa algebra with
unique trace $\tau$ then $\pi_{\tau} (A)^{\prime\prime} \cong R$
thus giving a positive answer to Popa's question in this case (see
Theorem \ref{thm:locallyreflexiveuniquetrace}).

\vspace{2mm}
\noindent{1.3. \bf Applications}

We already mentioned that amenable traces have applications
to certain questions in theoretical numerical analysis, K-homology
and are related to some natural questions concerning the
C$^*$-algebraic structure of the hyperfinite II$_1$-factor. We now
describe, in more detail, connections with some important open
problems in other areas of operator algebra theory.

\vspace{2mm} {\noindent\bf Elliott's Classification Program}

There are several new results which those in the classification
program may find of interest.  For example, we will observe that
certain questions around tracial approximation properties provide
both necessary conditions and sufficient conditions for various
cases of Elliott's conjecture to hold (cf.\ Propositions
\ref{thm:elliott} and \ref{thm:BFD}). As motivation for the other
results, we briefly describe the state of affairs (as we see it)
in the real rank zero case of the classification program.

That all approximately finite dimensional (AFD) II$_1$-factors are
isomorphic was known to Murray and von Neumann.  In
\cite{connes:classification} Connes proved that all injective
II$_1$-factors are AFD and hence fall under the Murray-von Neumann
classification theorem.  Elliott's classification program, in our
opinion, is in a state similar to the II$_1$-factor case prior to
Connes' work -- the classification theorems (for very large classes of
real rank zero algebras) exist and we should try to decide whether or
not the various hypotheses are always satisfied.  Of course, more
general classification theorems would always be welcome, but we
already have very good results which we should also be trying to
exploit.

For example, Huaxin Lin has succeeded in classifying his so-called
tracially AF algebras.  Roughly speaking, the only difference between
a Popa algebra and a (simple) tracially AF algebra is that the finite
dimensional algebra $B$ from Definition 1.2 is required (for tracially
AF algebras) to be `large in trace' (see \cite{lin:TAF} for the
precise definition and \cite{lin:TAFclassification} for the
classification theorem).  Lin's classification result is so exciting
because Popa proved in \cite{popa:simpleQD} that every simple, unital,
quasidiagonal C$^*$-algebra of real rank zero is a Popa algebra and hence,
`almost' tracially AF.  (Note that Popa requires no nuclearity or
K-theoretic hypotheses!)  This is the great advantage of tracially AF
algebras as opposed to the AH algebras considered in
\cite{elliott-gong} as there are no general hypotheses (yet) which
allow one to deduce an AH-type structure.

Thus Lin's theorem, combined with Popa's work, may (someday) nearly
complete the simple, quasidiagonal, real rank zero case of Elliott's
conjecture -- there are rather mild K-theoretic restrictions.  The
K-theory of a tracially AF algebra must be weakly unperforated and
satisfy the Riesz interpolation property.  However, every such
invariant arises from a tracially AF algebra and hence this case of
Elliott's conjecture is {\em equivalent} to the following question
(modulo a UCT assumption):

\vspace{2mm} {\em Is every simple, unital, nuclear, quasidiagonal
C$^*$-algebra with real rank zero, weakly unperforated K-theory and
having the Riesz decomposition property necessarily tracially AF?}
\vspace{2mm}

Thus the real question is what other general hypotheses (in addition
to quasidiagonality and real rank zero) are needed to ensure that the
finite dimensional algebras, which exist by Popa's work, can be taken
large in trace?  Are assumptions on K-theory enough?  Some experts
have felt that this question would have little to do with nuclearity,
primarily because Popa's work requires no such hypotheses.  (See, for
example, \cite[page 694]{lin:TAF} where it was ``tempting to
conjecture that every quasidiagonal, simple C$^*$-algebra of real rank
zero, stable rank one and with weakly unperforated $K_0$ is tracially
AF.'')  We will see, however, that K-theoretic assumptions are not
enough and, moreover, nuclearity must play a role (if the question
above is to be resolved affirmatively) as there exists an {\em exact}
Popa algebra with very nice K-theory which is not tracially AF (cf.\
Corollary \ref{thm:notTAF}). The one drawback of this example is that
it does not have a unique tracial state.  Popa has asked (private
communication) if a unique trace would be enough to ensure that Popa
algebras are tracially AF.  We will see that this is also not enough
(see Theorem \ref{thm:uniquetracenotTAF}). It is interesting to note,
however, that the case of an exact Popa algebra with unique trace
remains open.

On the other hand, inspired by recent work of Huaxin Lin
\cite{lin:ACtraces}, we will show that approximation properties of
traces provide an abstract hypothesis which does ensure that Popa
algebras are tracially AF.  That is, under certain technical
assumptions, we will show that (not necessarily nuclear) tracially
AF algebras are easily {\em characterized} by tracial
approximation properties (see Proposition
\ref{thm:characterization}).  Moreover, in the presence of
nuclearity and a unique trace it may turn out that this
approximation property is always satisfied (see the discussion in
Section \ref{thm:classificationprogram}).  It is the case that
type I C$^*$-algebras {\em always} have this nice tracial
approximation property and hence we are able to classify many
C$^*$-algebras which are built up out of type I algebras (for
example, thanks to the deep work of Q.\ Lin and Phillips, all
crossed products of compact manifolds by diffeomorphisms which
have real rank zero and unique trace).  Our classification results
are similar to  recent work  of Huaxin Lin, however the proofs are
significantly shorter and the present approach strikes us as
technically and conceptually simpler.

\vspace{2mm} {\noindent\bf Free Probability}

There are two new results arising from this work which are relevant to
free probability.  One is related to Connes' embedding problem and the
other is related to the semicontinuity and invariance questions for
Voiculescu's free entropy dimension.

Regarding Connes' embedding problem (i.e.\ the question of whether or
not microstates always exist) we obtain the following result (see
Theorem \ref{thm:embeddableMcDuff}): A II$_1$-factor $M$ embeds into
the ultraproduct of the hyperfinite II$_1$-factor if and only if there
exists a weakly dense C$^*$-subalgebra $A \subset M$ and a u.c.p.\ map
$\Phi : B(L^2(M)) \to M$ such that $\Phi(a) = a$ for all $a \in A$
(i.e.\ $M$ has Lance's WEP relative to a weakly dense subalgebra). We
also show that this is equivalent to the following: for every finite
subset ${\mathfrak F} \subset M$ and $\varepsilon > 0$ there exists an
operator system $X \subset M$ such that ${\mathfrak F}$ is
$\varepsilon$-contained in $X$ (in 2-norm) and $X$ is completely order
isomorphic to the hyperfinite II$_1$-factor (i.e.\ $M$ is quite
literally built out of the hyperfinite II$_1$-factor).  Hence Connes'
embedding problem predicts that the hyperfinite II$_1$-factor is {\em
the} basic building block for all II$_1$-factors.  Our own feeling is
that the world of II$_1$-factors is too exotic to expect that
everything is built up out of the nicest possible II$_1$-factor, but
we have not yet been able to construct a counterexample.  On the other
hand, this result also shows that many well known II$_1$-factors
(e.g.\ free group factors or property T factors coming from residually
finite groups) have a {\em dense} internal structure reminiscent of
the hyperfinite II$_1$-factor (but, of course, are not themselves
hyperfinite).\footnote{Though it is standard, it is a bit misleading
to use the terminology ``hyperfinite'' in this paragraph as
``injective'' would be more appropriate.  Of course, Connes'
uniqueness theorem asserts that these notions coincide for
II$_1$-factors, but the isomorphism identifying $R$ with a subspace $X
\subset M$ is not normal and hence the finite dimensional subalgebras
which are weakly dense in $R$ will map over to some small portion of
$X$.  However, injectivity is preserved so that $X$ is still an
injective operator system.}

As previously mentioned, we will show that every McDuff factor
contains a weakly dense Popa algebra (Theorem
\ref{thm:arbitraryMcDuff}).  Based on the proof of this result, in
joint work with Ken Dykema, we have been able to show that all the
(interpolated) free group factors $L({\mathbb F}_n)$ ($1 < n <
\infty$) contain finitely generated, weakly dense Popa algebras (see
\cite{BD}).  The point is that Popa algebras are fundamentally
different, at least from a C$^*$-algebraic point of view, than the
canonical generators of free group factors.  For example, Popa
algebras {\em never} contain two Haar unitaries which are free with
respect to some tracial state.\footnote{This is because such unitaries
would generated a copy of $C^*_r({\mathbb F}_2)$, but Popa algebras
are quasidiagonal and hence can't contain a copy of the
non-quasidiagonal C$^*$-algebra $C^*_r({\mathbb F}_2)$.}  In other
words, the techniques of this paper lead to new constructions of
generators of some well known II$_1$-factors.  This is relevant to the
invariance question for free entropy dimension as all new generators
give candidates for computation.  We make no attempt to compute free
entropy dimension in this paper -- our only point is that the results
presented here give new tools to construct ``exotic'' generators of
well known II$_1$-factors and we hope our work provides the foundation for
more examples in the future.

\begin{acknowledgement*}
This work began during the year-long program in operator algebras
at MSRI, 2000-2001.  I spoke to nearly everyone I encountered
during that year about various aspects of this work, as well as
most everyone I have encountered since, and it would be impossible
to recall all of the people who contributed remarks and ideas.
Instead I express my sincerest thanks to MSRI, the organizers and
participants and everyone else over the last three years who
patiently endured conversations on this topic. However, I must
specifically thank Marius Dadarlat and Dimitri Shlyakhtenko for a
(seemingly infinite) number of helpful discussions.  Finally I
thank the University of Tokyo, Yasu Kawahigashi and Taka Ozawa for
their hospitality as a significant part of the (re)writing of
these notes took place during a one year visiting position at the
University of Tokyo.
\end{acknowledgement*}

\chapter{Notation, definitions and useful facts}

The purpose of this chapter is to set our notation and list a
number of facts which will be used throughout.  Our notation should be
standard and the facts we collect here are all well known, very simple
or minor variations of the statements most commonly seen in the
literature.  In particular, there should be no harm in jumping to the
next chapter and referring back as necessary.

{\em Unless otherwise noted or obviously false,
all C$^*$-algebras are assumed to be unital and
separable.}  Similarly, all von Neumann algebras will be assumed to
have separable preduals (with the exception of $R^{\omega}$, which is
well known to be non-separable).

For a Hilbert space $H$, we will let $B(H)$ and $\mathcal{K}(H)$
denote the bounded and, respectively, compact operators on $H$.  The
canonical (unbounded, densely defined) trace on $B(H)$ will be denoted
by $\mathrm{Tr}$ while $\| \cdot \|$ will be the operator norm on
$B(H)$, $\| \cdot \|_{HS}$ and $\| \cdot \|_1$ will be the
Hilbert-Schmidt and $L^1$-norms, respectively\footnote{i.e.\
$\|T\|_{HS}^2 = \mathrm{Tr}(T^*T)$ and $\|T\|_1 = \mathrm{Tr}(|T|)$.},
and $\langle \cdot ,\cdot \rangle_{HS}$ will be the inner product on
the Hilbert space of Hilbert-Schmidt operators.\footnote{i.e.\
$\langle S, T \rangle_{HS} = \mathrm{Tr}(T^*S)$.}

When $A$ is a C$^*$-algebra with state $\eta$ we will denote the
associated GNS Hilbert space, representation and von Neumann
algebra by $L^2 (A,\eta)$, $\pi_{\eta} : A \to B(L^2 (A,\eta))$
and $\pi_{\eta}(A)^{\prime\prime}$, respectively.  Given $a \in
A$, $\hat{a} \in L^2 (A,\eta)$ will be the canonical image of $a$
in the GNS Hilbert space.

The symbols $\odot$, $\otimes$ and $\bar{\otimes}$ will denote the
algebraic, minimal and W$^*$-tensor products, respectively.

If $A$ is a C$^*$-algebra we will let $A^{op}$ denote the opposite
algebra (i.e.\ $A^{op} = A$ as involutive normed linear spaces, but
multiplication in $A^{op}$ is defined by $a \circ b = ba$; the latter
multiplication being the given multiplication in $A$).  $A^{**}$ will
denote the enveloping von Neumann algebra of $A$ (i.e.\ the Banach
space double dual of $A$).  Contrary to our standing assumption that
von Neumann algebras should have separable preduals, $A^{**}$ usually
has a non-separable predual.

If $\tau \in \TA$ is a tracial state then there is a canonical
antilinear isometry $J: L^2(A,\tau) \to L^2(A,\tau)$ defined by
$J(\hat{a}) = \hat{a^*}$.  One defines a `right regular
representation' (i.e.\ a $*$-homomorphism $\pi_{\tau}^{op} :
A^{op} \to B(L^2(A,\tau))$) by $\pi_{\tau}^{op} (a) = J\pi_{\tau}
(a^*)J$. Since $J\pi_{\tau}(A)J \subset \pi_{\tau}(A)^{\prime}$
one then gets an algebraic homomorphism $\pi_{\tau} \odot
\pi_{\tau}^{op} : A \odot A^{op} \to B(L^2(A,\tau))$ defined on
elementary tensors by $\pi_{\tau} \odot \pi_{\tau}^{op} (a\otimes
b) = \pi_{\tau} (a) \pi_{\tau}^{op} (b)$.  It is an important
fact, essentially due to Murray and von Neumann, that
$J\pi_{\tau}(A)^{\prime\prime}J = \pi_{\tau}(A)^{\prime}$ and
(hence) $\pi_{\tau}(A)^{\prime\prime} =
J\pi_{\tau}^{op}(A)^{\prime}J$.

Completely positive maps (cf.\ \cite{paulsen:cbmaps}) will play an
important role in these notes.  We will use the abbreviations c.p.\ and
u.c.p.\ for `completely positive' and `unital completely positive',
respectively.  We will often use two fundamental results concerning
such maps: Stinespring's dilation theorem and Arveson's extension
theorem (cf.\ \cite{paulsen:cbmaps}).

Multiplicative domains of u.c.p.\ maps will appear several times.  For
a u.c.p.\ map $\phi: A \to B$ we will let $A_{\phi}$ denote the
multiplicative domain (cf.\ \cite{paulsen:cbmaps}).  By definition,
$$A_{\phi} = \{a \in A: \phi(a^*a) = \phi(a^*)\phi(a) \ \mathrm{and}
\phi(aa^*) = \phi(a)\phi(a^*)\}$$ and it is not too hard to show that
this set can also be described as the set of elements which commute
with the Stinespring projection (in the Stinespring representation).
The key fact for us is that {\em u.c.p.\ maps are
bimodule maps over their multiplicative domains} -- i.e.\ if $a,c \in
A_{\phi}$ and $b \in A$ then $\phi(abc) = \phi(a)\phi(b)\phi(c)$.

The hyperfinite II$_1$-factor will appear many times and will
always be denoted by $R$.  We will use $\mathrm{tr}_n$ to denote
the unique tracial state on the $n\times n$ matrices.  For a von
Neumann algebra, $M$, with faithful, normal, tracial state $\tau$,
we will let $\| \cdot \|_{2,\tau}$ be the associated 2-norm (i.e.\
$\| x \|_{2,\tau} = \tau(x^* x)^{1/2}$).  If $M$ has a unique
trace (i.e.\ is a factor) then we will drop the dependence on
$\tau$ and simply write $\| \cdot \|_2$.

The ultraproduct of the hyperfinite II$_1$ factor $R^{\omega}$
will also appear several times.  That is, given a free ultrafilter
$\omega \in \beta {\mathbb N} \backslash {\mathbb N}$ one defines
an ideal $I_{\omega} \subset l^{\infty}(R) = \{ (x_n) \in \Pi_{n
\in {\mathbb N}} R: \sup_{n \in {\mathbb N}} \| x_n \| < \infty
\}$ by $I_{\omega} = \{ (x_n) \in l^{\infty}(R) : \lim_{n \to
\omega} \|x_n \|_2 = 0 \}$. Then the ultraproduct of $R$ with
respect to $\omega$ is defined to be the (C$^*$-algebraic)
quotient: $R^{\omega} = l^{\infty}(R)/I_{\omega}$.  $R^{\omega}$
is a II$_1$ factor with trace $\tau_{\omega} ((x_n) + I_{\omega})
= \lim_{n \to \omega} \tau_{R} (x_n)$.

The following simple fact will be useful.

\begin{lemc}
\label{thm:commutatorestimate}
Let $A \subset B(H)$ be a C$^*$-algebra and $P \in B(H)$ be a finite
rank projection. Then $$\| P a - a P\| = \max \{ \| Paa^*P - PaPa^*P
\|^{1/2}, \| Pa^*aP - Pa^*PaP \|^{1/2} \},$$ and $$\frac{\| P a - a P
\|_{HS}}{\| P \|_{HS}} = \bigg(\frac{\mathrm{Tr}
\big(Paa^*P - PaPa^*P\big) + \mathrm{Tr} \big(Pa^*aP -
Pa^*PaP\big)}{\mathrm{Tr}(P)} \bigg)^{1/2}.$$
\end{lemc}

\begin{proof}
Use the identity $$Pa - aP = Pa(1 -P) - (1 - P)aP,$$ the fact that
$Pa(1 -P)$ and $(1 - P)aP$ have orthogonal domains and ranges (hence
are orthogonal Hilbert-Schmidt vectors) and compute away.
\end{proof}

We will need the following version of Voiculescu's Theorem.

\begin{thmc}
\label{thm:Voiculescu'sthm}
Let $A \subset B(H)$ be in general position (i.e.\ $A \cap
\mathcal{K}(H) = \{ 0 \}$).  If $\phi : A \to M_n ({\mathbb C})$ is a
u.c.p.\ map then there exist isometries $V_k : {\mathbb C}^n \to H$
such that $\| \phi(a) - V_k^* a V_k \| \to 0$, for all $a \in A$, as
$k \to \infty$.  Moreover, letting $P_k = V_k V_k^*$, we have 
$$\lim_{k\to \infty}   \| P_k a - a P_k\| = \max \{ \| \phi(aa^*) -
\phi(a)\phi(a^*) \|^{\frac{1}{2}}, \| \phi(a^*a) - \phi(a^*)\phi(a) \|^{\frac{1}{2}}
\}$$ and $$\lim_{k\to \infty} \frac{\| P_k a - a P_k\|_{HS}}{\| P_k \|_{HS}} =
\bigg( \mathrm{tr}_n \big( \phi(aa^*) - \phi(a)\phi(a^*)\big) +
\mathrm{tr}_n \big(\phi(a^*a) - \phi(a^*)\phi(a)\big)
\bigg)^{\frac{1}{2}}.$$
\end{thmc}

\begin{proof}  That there exist isometries $V_k : {\mathbb C}^n \to H$
such that $\| \phi(a) - V_k^* a V_k \| \to 0$ is the first step in
proving the usual version of Voiculescu's Theorem (cf.\
\cite{arveson:extensions} or \cite{davidson}).  The commutator
estimates, which are the important part for us, follow from the
previous lemma.
\end{proof}

We will also need the following technical version of
Voiculescu's Theorem.  A proof can be found in \cite{brown:QDsurvey}
or \cite{brown:herrero}.

\begin{propc}
\label{thm:technicalVoiculescuThm}
Let $A \subset B(H)$ be in general position and $\Phi : A \to B(K)$ be
a u.c.p.\ map which is a faithful $*$-homomorphism modulo the compacts
(i.e.\ composing with the quotient map to the Calkin algebra yields a
faithful $*$-monomorphism $A \hookrightarrow Q(K)$).  Then there
exists a sequence of unitaries $U_n : K \to H$ such that for every $a
\in A$ we have
$$\limsup \| a - U_n \Phi(a) U_n^* \| \leq 2\max \{ \| \Phi(aa^*) -
\Phi(a)\Phi(a^*) \|^{\frac{1}{2}} , \| \Phi(a^* a) - \Phi(a^*)\Phi(a) \|^{\frac{1}{2}}
\}.$$
\end{propc}

We now list four simple facts which will also be useful. The first two
are well known and the second two are just a bit of trickery.  We
begin with a simple adaptation of the fact that if $M$ is a von
Neumann algebra with faithful, normal tracial state $\tau$ and $1_M
\in N \subset M$ is a sub-von Neumann algebra then there always exists a
$\tau$-preserving (hence faithful) conditional expectation $M \to N$
(cf.\ \cite[Exercise 8.7.28]{kadison-ringrose}).

\begin{lemc}
\label{thm:tracepreservingconditionalexpectation}
Let $A$ be a C$^*$-algebra, $\tau \in \TA$ be a tracial state and $1_A
\in B \subset A$ be a finite dimensional subalgebra. Then there exists
a conditional expectation $\Phi_B : A \to B$ such that $\tau \circ
\Phi_B = \tau$.
\end{lemc}

\begin{proof} Assume first that $\tau|_B$ is faithful.
Write $B \cong M_{n(1)} ({\mathbb C}) \oplus \cdots \oplus M_{n(k)}
({\mathbb C})$, and let $\{ e^{(1)}_{i,j} \}_{1 \leq i,j \leq n(1)}
\cup \ldots \cup \{ e^{(k)}_{i,j} \}_{1 \leq i,j \leq n(k)}$ be a
system of matrix units for $B$. Then the desired conditional
expectation is given by $$\Phi_B (x) = \sum\limits_{s = 1}^{k}
\sum\limits_{i,j = 1}^{n(s)}
\frac{\tau(xe_{i,j}^{(s)})}{\tau(e_{i,i}^{(s)})} e_{j,i}^{(s)}.$$

When $\tau|_B$ is not faithful, the formula above no longer makes
sense.  However, one can decompose $B$ as the direct sum of two
finite dimensional algebras, $B_0 \oplus B_f$, where $\tau|_{B_0}
= 0$ and $\tau|_{B_f}$ is faithful.  Letting $e_0$ (resp.\ $e_f$)
be the unit of $B_0$ (resp.\ $B_f$) we get a $\tau$-preserving
conditional expectation by mapping each $a \in A$ to $E_{B_0} (e_0
a e_0) + E_{B_f} (e_f a e_f)$, where $E_{B_0} : e_0 A e_0 \to B_0$
is any conditional expectation (which exists by finite
dimensionality) and $E_{B_f} : e_f A e_f \to B_f$ is a $\tau|_{e_f
A e_f}$ preserving conditional expectation as in the first part of
the proof.
\end{proof}

\begin{lemc}
\label{thm:tracepreservingembedding} If $B$ is a finite
dimensional C$^*$-algebra with tracial state $\tau$ then for every
$\varepsilon > 0$ there exists $n \in {\mathbb N}$ and a unital
$*$-monomorphism $\rho : B \hookrightarrow M_n({\mathbb C})$ such
that $| \tau (x) - \mathrm{tr}_n \circ \rho (x) | < \varepsilon \|
x \|$ for every $x \in B$.
\end{lemc}

\begin{proof}
If $\tau$ is a {\em rational} convex combination of extreme traces
then one can find an honestly trace preserving embedding by inflating
the summands of $B$ according to the rational numbers appearing in the
convex combination. The general case then follows by approximation.
\end{proof}

We remind the reader of the following theorem of Voiculescu (cf.\
\cite[Theorem 1]{dvv:QDhomotopy}): A separable, unital C$^*$-algebra
$A$ is quasidiagonal if and only if there exists a sequence of
u.c.p.\ maps $\varphi_n : A \to M_{k(n)} ({\mathbb C})$ which are
asymptotically multiplicative (i.e.\ $\| \varphi_n (ab) -
\varphi_n(a)\varphi_n(b) \| \to 0$ for all $a,b \in A$) and
asymptotically isometric (i.e.\ $\|a \| = \lim \| \varphi_n (a) \|$,
for all $a \in A$).

\begin{lemc}
\label{thm:trickery} Assume that $A$ is a quasidiagonal
C$^*$-algebra and $\psi_n : A \to M_{l(n)} ({\mathbb C})$ is an
asymptotically multiplicative (but not necessarily asymptotically
isometric) sequence of u.c.p.\ maps. Then there exists a sequence
of u.c.p.\ maps $\Phi_n : A \to M_{t(n)} ({\mathbb C})$ which are
asymptotically multiplicative, asymptotically isometric and such
that $| \mathrm{tr}_{t(n)}(\Phi_n (a)) - \mathrm{tr}_{l(n)}
(\psi_n (a)) | \to 0$ for all $a \in A$.
\end{lemc}

\begin{proof}
Let $\varphi_n : A \to M_{k(n)} ({\mathbb C})$ be an asymptotically
multiplicative, asymptotically isometric sequence of u.c.p.\ maps.
Choose integers $s(n)$ such that $\frac{k(n)}{s(n)} \to 0$; thus 
$\frac{s(n)l(n)}{s(n)l(n) + k(n)} \to
1$.  Then one defines
$\Phi_n : A \to M_{s(n)l(n) + k(n)} ({\mathbb C})$ to be the block
diagonal map with one summand equal to $\varphi_n$ and $s(n)$ summands
equal to $\psi_n$.
\end{proof}

Note that the previous lemma can be formulated in terms of
$*$-homomorphisms when $A$ is a residually finite dimensional
$C^*$-algebra.

Though we will try to keep everything unital, non-unital maps are
sometimes unavoidable.  The next lemma keeps everything running
smoothly.

\begin{lemc}
\label{thm:nonunital} Let $\phi_n :A \to M_{k(n)} ({\mathbb C})$
be contractive c.p.\ (c.c.p.) maps (though we still assume $A$ is
unital).
\begin{enumerate}
\item If $\| \phi_n(ab) - \phi_n(a)\phi_n(b) \|_2 \to 0$ for all
$a,b \in A$ and $\mathrm{tr}_{k(n)} (\phi_n (1_A)) \to 1$ then
there exist u.c.p.\ maps $\psi_n : A \to M_{k(n)} ({\mathbb C})$
which are also asymptotically multiplicative with respect to
2-norms and such that $| \mathrm{tr}_{k(n)}(\psi_n (a)) -
\mathrm{tr}_{k(n)}(\phi_n(a)) | \to 0$ as $n \to \infty$, with
convergence being uniform on the unit ball of $A$.

\item If $\| \phi_n(ab) - \phi_n(a)\phi_n(b) \| \to 0$ for all
$a,b \in A$ and $\mathrm{tr}_{k(n)} (\phi_n (1_A)) \to 1$ then
there exist integers $l(n) \leq k(n)$ and asymptotically multiplicative 
u.c.p.\ maps $\psi_n : A
\to M_{l(n)} ({\mathbb C})$ such that $| \mathrm{tr}_{l(n)}(\psi_n
(a)) - \mathrm{tr}_{k(n)}(\phi_n(a)) | \to 0$ as $n \to \infty$,
with convergence being uniform on the unit ball of $A$.  
\end{enumerate}
\end{lemc}

\begin{proof} For the proof of (1) we will need an observation of
Choi-Effros (cf.\ \cite[Lemma
2.2]{choi-effros:injectivityandoperatorspaces}) which says that we
can find u.c.p.\ maps $\psi_n : A \to M_{k(n)} ({\mathbb C})$ such
that $\phi_n (a) = c_n \psi_n (a) c_n$ for all $n$ and $a \in A$,
where $c_n = \phi_n (1_A)^{1/2}$.  Since $\mathrm{tr}_{k(n)}
(\phi_n (1_A)) \to 1$ it follows that $\| c_n - 1_{M_{k(n)}} \|_2
\to 0$.   Also we note the following general inequality for all $x
\in M_{k(n)} ({\mathbb C})$: $$\| x - c_n x c_n \|_2 \leq 2 \| x
\| \| c_n - 1_{M_{k(n)}} \|_2.$$  Hence we find that
\begin{eqnarray*}
\|\psi_n(ab) - \psi_n(a)\psi_n(b) \|_2
&\leq&
\|\psi_n(ab) - c_n \psi_n (ab) c_n\|_2 + \|\phi_n(ab) -
\phi_n(a)\phi_n(b) \|_2\\
& \ &
+ \|\phi_n(a)\phi_n(b) - \psi_n(a)\psi_n(b) \|_2\\
&\leq&
2\|ab\|\| c_n - 1_{M_{k(n)}} \|_2 + \|\phi_n(ab) -
\phi_n(a)\phi_n(b) \|_2\\
& \ &
+ 4\|a\|\|b\|\| c_n - 1_{M_{k(n)}} \|_2\\
&\leq&
6\|a\|\|b\|\| c_n - 1_{M_{k(n)}} \|_2 + \|\phi_n(ab) -
\phi_n(a)\phi_n(b) \|_2.\\
\end{eqnarray*}
This shows asymptotic multiplicativity of the $\psi_n$'s. 

The tracial approximation part is a simple application of the
Cauchy-Schwartz inequality: since $| \mathrm{tr}_{k(n)}(\psi_n (a)) -
\mathrm{tr}_{k(n)}(c_n\psi_n(a) c_n)
|=|\mathrm{tr}_{k(n)}(\psi_n(a)(1_{M_{k(n)}} - c_n^2))|$ we see that
$$| \mathrm{tr}_{k(n)}(\psi_n (a)) -
\mathrm{tr}_{k(n)}(c_n\psi_n(a) c_n)
| \leq
\|a\|\| c_n^2 - 1_{M_{k(n)}} \|_2.$$

The proof of part (2) uses similar estimates as those above but
starts off a little different.  Hence we only describe how to get
the $\psi_n$'s and leave the remaining details to the reader. In
the case that $\|\phi_n(ab) - \phi_n(a)\phi_n(b) \| \to 0$ we have
that the $\phi_n(1_A)$'s are getting closer and closer to
projections (by functional calculus), say $P_n$.  One checks that 
the c.c.p.\ maps $A \to
M_{l(n)}(\mathbb{C})$ defined by $a \mapsto P_n \phi_n(a) P_n$ are 
also asymptotically multiplicative -- $\| P_n \phi_n(a) - \phi_n(a) P_n\| 
\to 0$, for all $a \in A$, since $\phi_n(1_A)$ is asymptotically central. 
Unfortunately they still aren't unital, but $\phi_n(1)$
is invertible in $P_n M_{k(n)}(\mathbb{C}) P_n$, since it is very
close in norm to the identity, and hence we get u.c.p.\ maps by
defining $$\psi_n(a) = (P_n \phi_n(1))^{-\frac{1}{2}}(P_n
\phi_n(a) P_n)(P_n \phi_n(1))^{-\frac{1}{2}}.$$  That these maps
$\psi_n$ are still asymptotically multiplicative in norm is
similar to the 2-norm case above as is the tracial approximation
property.
\end{proof}

Finally, we will need the following consequence of the
Hahn-Banach theorem.  This result is a good exercise (for those not
already familiar with it).

\begin{lemc}
\label{thm:HB}
Let $S(A)$ denote the state space of a C$^*$-algebra $A$.  Assume that
${\mathcal S} \subset S(A)$ is a set of states with the property that
for each {\em self-adjoint} $a \in A$ we have $$\|a\| = \sup_{\phi \in
{\mathcal S}} \{ |\phi(a)| \}.$$ Then the (weak-$*$) closed, convex hull of
${\mathcal S}$ is equal to $S(A)$.
\end{lemc}

\chapter{Amenable traces and stronger approximation properties}

We will now develop some of the basic theory of amenable traces.  Our
first goal is to show that this notion has numerous equivalent
formulations, just as amenable groups admit numerous
characterizations.  These fundamental results, which are essentially
due to Connes in the unique trace case and Kirchberg in general, are
non-trivial and have already appeared in the literature (cf.\
\cite{connes:classification}, \cite{kirchberg:propertyTgroups}).
However, due to its importance, we feel that a detailed presentation
of Ozawa's (substantially simplified but still quite technical) proof
is worthwhile (cf.\ \cite{ozawa:QWEP}).

As previously mentioned, one of the characterizations of amenable
traces is a natural finite dimensional approximation property.  In
later sections we will define several subsets of the amenable
traces by simply strengthening this approximation property in various
ways.  At this point our goal is just to define these subsets and make a
few general observations.  We will discuss examples and applications in the
later.

We remind the reader that, unless obviously false,  all
C$^*$-algebras (resp.\ von Neumann algebras) are assumed to be
separable and unital (resp.\ have separable preduals).  We also
remind the reader that all notation not defined below should have
been explained in the previous chapter.

\section{Characterizations of amenable traces}

We first recall the definition of an amenable trace.

\begin{defn}
Let $A \subset B(H)$ be a concretely represented, unital
C$^*$-algebra.  A state $\tau$ on $A$ is called an {\em amenable
trace} if there exists a state $\phi$ on $B(H)$ such that (1)
$\phi|_A = \tau$ and (2) $\phi(uTu^*) = \phi(T)$ for every unitary
$u \in A$ and $T \in B(H)$.  We will denote the set of amenable
traces on $A$ by $\TAim$.
\end{defn}

Though it appears that this set of traces may depend on the choice of
representation $A \subset B(H)$ it turns out that this is not the case.

\begin{prop}
Let $A \subset B(H)$ be a C$^*$-algebra, $\tau \in \TAim$ and $\pi:A\to
B(K)$ be any other faithful representation of $A$. Then $\tau\circ
\pi^{-1}$ is an amenable trace on $\pi(A)$.  In other words, the set
$\TAim$ does not depend on the choice of faithful representation $A
\subset B(H)$ and hence being an amenable trace is a natural {\em
abstract} property of a tracial state.
\end{prop}

\begin{proof}
By Arveson's extension theorem we can find a u.c.p.\ map $\Psi:B(K)
\to B(H)$ such that $\Psi(\pi(a)) = a$ for all $a \in A$. Note that
$\pi(A)$ belongs to the multiplicative domain of $\Psi$ since
$\Psi|_{\pi(A)}$ is a $*$-homomorphism (namely, $\pi^{-1}$).  Define a
state $\psi$ on $B(K)$ by $\psi = \phi\circ\Psi$.  Now for arbitrary
$T \in B(K)$ and unitary $u \in A$ we have $$\psi(\pi(u)T\pi(u^*)) =
\phi\big(\Psi(\pi(u)) \Psi(T)\Psi(\pi(u^*)) \big) = \phi(u\Psi(T)u^*)
= \phi(\Psi(T)) = \psi(T).$$  Hence $\tau\circ
\pi^{-1}$ is an amenable trace on $\pi(A)$.
\end{proof}

We now begin the long and technical, though essentially elementary,
journey to the various characterizations of amenable traces.  This
will require numerous calculations\footnote{The diligent reader should
go find a very large chalkboard before reading further!} as well as
the {\em Powers-St{\o}rmer} inequality.

\begin{prop}[Powers-St{\o}rmer]
If $a,b \in B(H)$ are positive trace class operators then $\| a - b
\|_{HS}^2 \leq \| a^2 - b^2 \|_{1}$.  In
particular, if $u \in B(H)$ is a unitary and $h \geq 0$ has finite
rank then $\| uh^{1/2} - h^{1/2}u \|_{HS} = \| uh^{1/2}u^*
- h^{1/2} \|_{HS} \leq \| uhu^* - h
\|_{1}^{1/2}$.
\end{prop}

\begin{proof} We could refer to the original paper, of course, but the
proof is short and elementary so we include it.

Let $\{v_k\}$ be an orthonormal basis of $H$ consisting of
eigenvectors of $a - b$ and let $\lambda_k$ be the corresponding
(real) eigenvalues.  Note that since $a+b \geq a-b$ and $a+b \geq
-(a-b)$ it follows that $$\langle (a+b)v,v \rangle \geq |\langle
(a-b)v,v \rangle|,$$ for all $v \in H$.  Note also that for any
self-adjoint $x$ we have the inequality $|\langle x v,v\rangle | \leq
\langle |x| v,v \rangle.$

Now we compute
\begin{eqnarray*}
\| a - b\|_{HS}^2
&=&
\mathrm{Tr}(|a - b|^2)\\
&=&
\sum_{k} \langle |a - b|^2 v_k,v_k \rangle\\
&=&
\sum_k |\lambda_k|^2\\
&\leq&
\sum_k |\lambda_k|\langle (a+b)v_k,v_k \rangle\\
&=&
\sum_k |\langle (a+b)(\frac{1}{2}\lambda_k v_k),v_k \rangle +
\langle (a+b)v_k, \frac{1}{2}\lambda_k v_k\rangle|\\
&=&
\sum_k |\langle \frac{1}{2}\big((a+b)(a-b) +
(a-b)(a+b)\big)v_k,v_k \rangle|\\
&=&
\sum_k |\langle (a^2 - b^2)v_k,v_k \rangle|\\
&\leq&
\sum_k \langle |a^2 - b^2|v_k,v_k \rangle\\
&=&
\| a^2 - b^2 \|_{1}.
\end{eqnarray*}
\end{proof}

Our next lemma is the key technical ingredient. The proof given
below is due to Ozawa \cite[Theorem 6.1]{ozawa:QWEP} -- with a few
more details thrown in.

\begin{lem}
\label{thm:ozawalemma}
Let $h \in B(H)$ be a positive, finite rank operator with rational
eigenvalues and $\mathrm{Tr}(h) = 1$.  Then there exists a u.c.p.\ map
$\phi : B(H) \to M_q({\mathbb C})$ such that $\mathrm{tr_q}(\phi(T)) =
{\rm Tr}(hT)$ for all $T \in B(H)$ and $|\mathrm{tr_q}\big(\phi(uu^*) -
\phi(u)\phi(u^*)\big)| \leq 2\| uhu^* - h \|_{1}^{1/2}$ for
every unitary operator $u \in B(H)$.
\end{lem}

\begin{proof}
Let $v_1, \ldots, v_k \in H$ be the eigenvectors of $h$ and
$\frac{p_1}{q}, \ldots, \frac{p_k}{q}$ the corresponding eigenvalues.
Note that $\sum p_j = q$ since Tr$(h) = 1$. Let $\{ w_m \}$ be any
orthonormal basis of $H$ and consider the following orthonormal subset
of $H\otimes H$: $$\{ v_1 \otimes w_1, \ldots, v_1 \otimes w_{p_1} \}
\cup \{ v_2 \otimes w_1, \ldots, v_2 \otimes w_{p_2} \} \cup \ldots
\cup \{ v_k \otimes w_1, \ldots, v_k \otimes w_{p_k} \}.$$ Let
$P \in B(H\otimes H)$ be the orthogonal projection onto the span of
these vectors.

Our first task is to write down the matrix of $P (T\otimes 1) P$ (in
the basis above), for an arbitrary $T \in B(H)$.  Though we have no
intention of doing this completely we will make a few remarks.

The matrix of $P (T\otimes 1) P$ decomposes into $k\times k$ blocks
(which are not square!) as follows
$$P (T\otimes 1) P =
\begin{pmatrix}
A_{1,1} & A_{1,2} & \cdots & A_{1,k}\\
A_{2,1} & A_{2,2} & \cdots & A_{2,k}\\
\vdots & \vdots & \ddots & \vdots\\
A_{k,1} & A_{k,2} & \cdots & A_{k,k},\\
\end{pmatrix}$$
where the matrix $A_{i,j}$ has $p_i$ rows and $p_j$ columns.
The matrices $A_{i,i}$ look like
$$
\begin{pmatrix}
\langle Tv_i,v_i \rangle & 0 & \cdots & 0\\
0 & \langle Tv_i,v_i \rangle & \cdots & 0\\
\vdots & \vdots & \ddots & \vdots\\
0 & 0 & \cdots & \langle Tv_i,v_i\rangle\\
\end{pmatrix}.$$
For $i \neq j$ the matrices $A_{i,j}$ look like
$$
\begin{pmatrix}
\langle Tv_j,v_i \rangle & 0 & \cdots & 0\\
0 & \langle Tv_j,v_i \rangle & \cdots & 0\\
\vdots & \vdots & \ddots & \vdots\\
\end{pmatrix},$$
but are not necessarily square (unless $p_i = p_j$).  In particular
note that the number of $\langle Tv_j,v_i \rangle$'s appearing is
equal to $\min\{p_i,p_j\}$.

Now one computes the matrix of $P(T\otimes 1)P(T^* \otimes 1)P$.
Having done so the following facts become obvious.

\begin{enumerate}
\item $\frac{1}{q} \mathrm{Tr}(P(T\otimes 1)P) = \frac{1}{q}\big(\sum_{i =
1}^{k} p_i \langle Tv_i, v_i \rangle \big) = \sum_{i = 1}^{k} \langle
Thv_i, v_i \rangle = \mathrm{Tr}(Th)$.

\item $\frac{1}{q} \mathrm{Tr}(P(T\otimes 1)P(T^* \otimes 1)P) =
\frac{1}{q}\sum_{i,j = 1}^{k} |\langle Tv_j, v_i \rangle|^2
\min\{p_i, p_j\}$.
\end{enumerate}

As if that wasn't bad enough, one should now write down the
matrices of $h^{1/2}T$, $h^{1/2}T^*$ and $h^{1/2}Th^{1/2}T^*$ (in any
orthonormal basis which begins with $\{ v_1, \ldots, v_k \}$).  Having
done so one immediately sees that, letting $T_{i,j} = \langle Tv_j,
v_i \rangle$,
$$\mathrm{Tr}(h^{1/2}Th^{1/2}T^*) = \sum_{i,j = 1}^{k} \frac{1}{q}(p_i
p_j)^{1/2} |T_{i,j}|^2.$$ Hence, if we define a u.c.p.\ map $\phi :
B(H) \to M_q({\mathbb C})$ by $\phi(T) = P(T\otimes 1)P$ then
$\mathrm{tr_q}(\phi(T)) = \mathrm{Tr}(hT)$ for all $T \in B(H)$.  For TeXnical 
reasons we define $(\ast) = |\mathrm{Tr}(h^{1/2}Th^{1/2}T^*) - 
\mathrm{tr_q}(\phi(T)\phi(T^*))|$ and observe the following estimates:
\begin{eqnarray*}
(\ast)
& = & \sum_{i,j = 1}^{k}
\frac{1}{q}|T_{i,j}|^2 \bigg((p_i p_j)^{1/2} - \min\{p_i, p_j\} \bigg)\\
& \leq &  \sum_{i,j = 1}^{k} \frac{1}{q}|T_{i,j}|^2 p_i^{1/2}|p_i^{1/2}
          -  p_j^{1/2}|\\
& \leq & \bigg(\sum_{i,j = 1}^{k} \frac{1}{q}|T_{i,j}|^2 p_i   \bigg)^{1/2}
         \bigg(\sum_{i,j = 1}^{k} \frac{1}{q}|T_{i,j}|^2 (p_i^{1/2}
          -  p_j^{1/2})^2 \bigg)^{1/2}\\
& = & \| Th^{1/2} \|_{HS} \| h^{1/2}T - Th^{1/2}
\|_{HS}.
\end{eqnarray*}
Now if $T$ happens to be a unitary operator then $\| Th^{1/2}
\|_{HS} = \| h^{1/2} \|_{HS} = 1$ and $\|
h^{1/2}T - Th^{1/2} \|_{HS} = \| Th^{1/2}T^* - h^{1/2}
\|_{HS}$ and hence we can apply the Powers-St{\o}rmer
inequality after the inequalities above to get:
$$|\mathrm{Tr}(h^{1/2}Th^{1/2}T^*) - \mathrm{tr}(\phi(T)\phi(T^*))|
\leq \| ThT^* - h \|_{1}^{1/2}.$$ Finally, the
Cauchy-Schwartz inequality applied to the Hilbert-Schmidt operators
implies that for every unitary operator $T \in B(H)$,
\begin{eqnarray*}
\mathrm{tr_q}(\phi(TT^*) - \phi(T)\phi(T^*)) & \leq & |1 -
\mathrm{Tr}(h^{1/2}Th^{1/2}T^*)|
         + \| ThT^* - h \|_{1}^{1/2}\\
& = & |\mathrm{Tr}(ThT^*) - \mathrm{Tr}(h^{1/2}Th^{1/2}T^*)| +
         \| ThT^* - h \|_{1}^{1/2}\\
& = & |\mathrm{Tr}((Th^{1/2} - h^{1/2}T)h^{1/2}T^*)| +
         \| ThT^* - h \|_{1}^{1/2}\\
& \leq & \| h^{1/2}T^* \|_{HS} \| Th^{1/2} - h^{1/2}T
\|_{HS}
         + \| ThT^* - h \|_{1}^{1/2}\\
& \leq & 2\| ThT^* - h \|_{1}^{1/2}.
\end{eqnarray*}
\end{proof}

\begin{lem}
\label{thm:estimate}
If $\phi:A \to M_n({\mathbb C})$ is a u.c.p.\ map then for all $a,b
\in A$ we have
$$\| \phi(ab) - \phi(a)\phi(b) \|_{2} \leq \|b\|\bigg(
\mathrm{tr}_n \big( \phi(aa^*) - \phi(a)\phi(a^*)\big) + \mathrm{tr}_n
\big(\phi(a^*a) - \phi(a^*)\phi(a)\big) \bigg)^{1/2}.$$
\end{lem}

\begin{proof} Thanks to Stinespring's dilation theorem we may assume
that $A \subset B(H)$ and $\phi(a) = PaP$, for all $a \in A$, where $P
\in B(H)$ is a finite rank projection. (If Stinespring's
representation is not faithful just dilate it further.)  Note that tr$_n$
on $M_n({\mathbb C}) = PB(H)P$ is then identified with the linear
functional $\frac{\mathrm{Tr}(\cdot)}{\mathrm{Tr}(P)}$ and hence the
2-norm $\| \cdot \|_{2}$ gets identified with $$\frac{\|
\cdot \|_{HS}}{\|P \|_{HS}}.$$

Hence, if $a,b \in A$ are given we can apply Lemma \ref{thm:commutatorestimate} 
to $\| \phi(ab) - \phi(a)\phi(b) \|_{2} = (\ast)$ (again, for TeXnical reasons) and 
deduce that 
\begin{eqnarray*}
(\ast)
&=&
\frac{\|PabP - PaPbP \|_{HS}}{\|P \|_{HS}}\\
&=&
\frac{\|P(Pa - aP)bP \|_{HS}}{\|P \|_{HS}}\\
&\leq&
\|b\| \frac{\|Pa - aP \|_{HS}}{\|P \|_{HS}}\\
&=&
\|b\|\bigg(\frac{\mathrm{Tr}
\big(Paa^*P - PaPa^*P\big) + \mathrm{Tr} \big(Pa^*aP -
Pa^*PaP\big)}{\mathrm{Tr}(P)} \bigg)^{1/2}\\
&=&
\|b\|\bigg(\mathrm{tr}_n \big( \phi(aa^*) - \phi(a)\phi(a^*)\big) +
\mathrm{tr}_n \big(\phi(a^*a) - \phi(a^*)\phi(a)\big) \bigg)^{1/2}.
\end{eqnarray*}
\end{proof}

Recall that if $\tau$ is a tracial state on $A$ then there is a
``right regular'' representation $\pi_{\tau}^{op}:A^{op} \to
B(L^2(A,\tau))$ with the property that $\pi_{\tau}(A)^{\prime} =
\pi_{\tau}^{op}(A^{op})^{\prime\prime}$ and
$\pi_{\tau}(A)^{\prime\prime} = \pi_{\tau}^{op}(A^{op})^{\prime}$.  In
particular, there is a natural $*$-homomorphism $$A \odot A^{op} \to
B(L^2(A,\tau)), a\otimes b \to \pi_{\tau}(a)\pi_{\tau}^{op}(b).$$
Composing this representation with the vector state $x \mapsto \langle
x \hat{1},\hat{1} \rangle$, where $\hat{1} \in L^2(A,\tau)$ denotes
the natural image of the unit of $A$, we get a positive linear
functional $\mu_{\tau}$ on $A \odot A^{op}$ which will play a
distinguished role in what follows.  The reader not familiar with this
construction is advised to work out the case $A = M_n({\mathbb C})$,
$\tau = \mathrm{tr}_n$ as it is not only instructive but will also be
used below.

We are now in a position to prove the main characterizations of
amenable traces.

\begin{thm}[\cite{kirchberg:propertyTgroups}]
\label{thm:amenabletraces} Let $\tau$ be a tracial state on $A$.
Then the following are equivalent:
\begin{enumerate}
\item $\tau$ is amenable.

\item There exists a sequence of u.c.p.\ maps $\phi_n : A \to M_{k(n)}
(\mathbb C)$ such that $\| \phi_n (ab) - \phi_n (a) \phi_n (b)
\|_{2} \to 0$ and $\tau (a) = \lim_{n \to \infty}
\mathrm{tr}_{k(n)} \circ \phi_n (a)$, for all $a,b \in A$.

\item The positive linear functional $\mu_{\tau}: A \odot A^{op} \to
{\mathbb C}$ is continuous with respect to the minimal tensor product norm
(i.e.\ extends to a state on $A\otimes A^{op}$).

\item The natural $*$-homomorphism $A \odot A^{op} \to B(L^2(A,\tau))$
is continuous with respect to the minimal tensor product norm (i.e.\ extends
to a representation of $A\otimes A^{op}$).

\item For any faithful representation $A \subset B(H)$ there
exists a u.c.p.\ map $\Phi : B(H) \to
\pi_{\tau}(A)^{\prime\prime}$ such that $\Phi(a) =
\pi_{\tau}(a)$.\footnote{The proof will show that this is also
equivalent to knowing that {\em there exists} a faithful
representation $A \subset B(H)$ and there exists a u.c.p.\ map
$\Phi : B(H) \to \pi_{\tau}(A)^{\prime\prime}$ such that $\Phi(a)
= \pi_{\tau}(a)$.}
\end{enumerate}
\end{thm}

\begin{proof}
(1) $\Longrightarrow$ (2).\footnote{The main idea in this
implication comes directly from Connes' celebrated uniqueness
theorem for the injective II$_1$-factor
\cite{connes:classification}.} Let $A \subset B(H)$ be a faithful
representation. Since $\tau$ is an amenable trace we can find a
state $\psi$ on $B(H)$ which extends $\tau$ and such that $\psi(u
T u^*) = \psi(T)$ for all unitaries $u \in A$ and operators $T \in
B(H)$. Since the normal states on $B(H)$ are dense in the set of
all states on $B(H)$ we can find a net of positive operators
$h_{\lambda} \in {\mathcal T}$ such that $\mathrm{Tr}(h_{\lambda}
T) \to \psi(T)$ for all $T \in B(H)$. Since $\psi(u^* T u) =
\psi(T)$ it follows that $\mathrm{Tr}(h_{\lambda} T) -
\mathrm{Tr}((uh_{\lambda}u^*) T) \to 0$ for every $T \in B(H)$ and
unitary $u \in A$.  In other words, for every unitary $u \in A$
the net of trace class operators $h_{\lambda} - uh_{\lambda}u^*$
tends to zero in the weak topology (coming from $B(H)$).  Hence,
by the Hahn-Banach theorem, there are convex combinations which
tend to zero in L$^1$-norm.  In fact, taking finite direct sums
(i.e.\ considering n-tuples $(u_1 h_{\lambda} u_1^* - h_{\lambda},
\ldots, u_n h_{\lambda} u_n^* - h_{\lambda})$) one applies a
similar argument to show that if ${\mathfrak F} \subset A$ is a
finite set of unitaries then for every $\epsilon > 0$ we can find
a positive trace class operator $h$ such that $\mathrm{Tr}(h) =
1$, $|\mathrm{Tr}(uh) - \tau(u)| < \epsilon$ and $\| h - uhu^*
\|_1 < \epsilon$ for all $u \in {\mathfrak F}$.  Since finite rank
operators are norm dense in the trace class operators we may
further assume that $h$ is finite rank with rational eigenvalues.

Applying Lemma \ref{thm:ozawalemma} to bigger and bigger finite
sets of unitaries and smaller and smaller $\epsilon$'s we can
construct a sequence of u.c.p.\ maps $\phi_n : B(H) \to
M_{k(n)}({\mathbb C})$ such that tr$_{k(n)}(\phi_n (u)) \to
\tau(u)$ and $|\mathrm{tr}_{k(n)}(\phi_n(uu^*) -
\phi_n(u)\phi_n(u^*))| \to 0$ for every unitary $u$ in a countable
set with dense linear span in $A$. Of course, we may further
assume that this set of unitaries is closed under the adjoint
operation.  From Lemma \ref{thm:estimate} it follows that for
every unitary in this set and every $a \in A$ we have
$$\|\phi_n (ua) - \phi_n (u) \phi_n (a) \|_{2} \to 0.$$ Since every
element in $A$ is a linear combination of four unitaries it follows
that $$\|\phi_n (ab) - \phi_n (a) \phi_n (b) \|_{2} \to 0$$ and
$$\mathrm{tr}_{k(n)}(\phi_n(a)) \to \tau(a)$$ for all $a,b \in A$.

(2) $\Longrightarrow$ (3). We first note that it suffices to show
that $\mu_{\tau}$ is the weak-$*$ limit of {\em min}-continuous
functionals.  That is, if we can find a sequence of functionals
$\sigma_n:A\otimes A^{op} \to {\mathbb C}$ such that
$\sigma_n(a\otimes b) \to \mu_{\tau}(a\otimes b)$ for all $a \in
A, b \in A^{op}$ then it will follow that any weak-$*$ cluster
point of the $\sigma_n$'s will be an extension of $\mu_{\tau}$ to
$A\otimes A^{op}$ (i.e.\ $\mu_{\tau}$ is continuous with respect
to the minimal tensor product norm).

So let $\phi_n : A \to M_{k(n)} (\mathbb C)$ be a sequence of
u.c.p.\ maps with the property that $\| \phi_n (ab) - \phi_n (a)
\phi_n (b) \|_{2} \to 0$ and $\tau (a) = \lim_{n \to \infty}
\mathrm{tr}_{k(n)} \circ \phi_n (a)$, for all $a,b \in A$.  Note
we can also regard these maps as going from $A^{op}$ to $M_{k(n)}
(\mathbb C)^{op}$ and they are still u.c.p.\ maps.  To distinguish
them, however, we let $\phi_n^{op} : A^{op} \to M_{k(n)} (\mathbb
C)^{op}$ be the ``opposite'' maps (though they are literally the
same as $\phi_n$ as maps of linear spaces).  Since u.c.p.\ maps
behave well with respect to minimal tensor products we may
consider the u.c.p.\ maps
$$\phi_n\otimes \phi_n^{op}: A\otimes A^{op} \to M_{k(n)} (\mathbb
C)\otimes M_{k(n)} (\mathbb C)^{op}.$$ All we need to show is that
there is a state $\mu_n$ on $M_{k(n)} (\mathbb C)\otimes M_{k(n)}
(\mathbb C)^{op}$ such that $\mu_n\circ \phi_n\otimes \phi_n^{op} \to
\mu_{\tau}$ in the weak-$*$ topology.  In this picture it may not be
so clear what the right state is but if we identify $M_{k(n)} (\mathbb
C)\otimes M_{k(n)} (\mathbb C)^{op}$ with $B(L^2(M_{k(n)} (\mathbb
C),\mathrm{tr}_{k(n)}))$, $M_{k(n)} (\mathbb C)\otimes 1$ with the
image of $M_{k(n)} (\mathbb C)$ in the GNS representation with respect to the
unique tracial state and $1\otimes M_{k(n)} (\mathbb C)^{op}$ with the
commutant of $M_{k(n)} (\mathbb C) \subset B(L^2(M_{k(n)} (\mathbb
C),\mathrm{tr}_{k(n)}))$ then it should be clear which state to pick.
Let $\mu_n$ be the vector state on $B(L^2(M_{k(n)} (\mathbb
C),\mathrm{tr}_{k(n)}))$ given by $x \mapsto \langle x \hat{1},\hat{1}
\rangle$. Observe that $$\mu_n(\phi_n\otimes \phi_n^{op}(a\otimes b))
= \langle \phi_n(a)J\phi_n(b^*)J \hat{1},\hat{1} \rangle =
\mathrm{tr}_{k(n)}(\phi_n(a)\phi_n(b)).$$ The Cauchy-Schwartz
inequality shows that for all $x \in M_{k(n)} (\mathbb C)$,
$|\mathrm{tr}_{k(n)}(x)| \leq \|x\|_{2}$ and hence
$$|\mathrm{tr}_{k(n)}(\phi_n(a)\phi_n(b)) -
\mathrm{tr}_{k(n)}(\phi_n(ab))| \to 0.$$ Thus we see that
$$\mu_n(\phi_n\otimes \phi_n^{op}(a\otimes b)) \to \tau(ab) = \langle
\pi_{\tau}(a)\pi_{\tau}(b) \hat{1},\hat{1} \rangle =
\mu_{\tau}(a\otimes b),$$ for all $a,b \in A$.

(3) $\Longrightarrow$ (4) follows from (the proof of) uniqueness of GNS
representations.  Indeed, on the norm dense $*$-subalgebra $A\odot
A^{op} \subset A\otimes A^{op}$ we can construct a unitary operator
which conjugates the representation $A\odot A^{op} \to B(l^2(A,\tau))$
to the GNS representation of $A\otimes A^{op}$ with respect to $\mu_{\tau}$.

(4) $\Longrightarrow$ (5) uses an argument of Lance which we feel is
the single most useful trick in the theory of tensor products. Since
$A\otimes A^{op} \subset B(H) \otimes A^{op}$ we can extend the
$*$-homomorphism $\pi_{\tau} \otimes \pi_{\tau}^{op} : A\otimes A^{op}
\to B(L^2(A,\tau))$ to a completely positive map $\Phi : B(H) \otimes
A^{op} \to B(L^2(A,\tau))$.  Since $\Phi|_{A\otimes A^{op}}$ is a
homomorphism it follows that $A \otimes A^{op}$ (and, in particular,
$1\otimes A^{op}$) is in the multiplicative domain of $\Phi$.  Hence,
for every $T \in B(H)$, it follows that $\Phi(T \otimes 1) \in
\Phi(1\otimes A^{op})^{\prime} = \pi_{\tau}^{op} (A^{op})^{\prime} =
\pi_{\tau}(A)^{\prime\prime}$.

(5) $\Longrightarrow$ (1). Since $A$ is contained in the
multiplicative domain of $\Phi$, it is easy to verify that $\varphi(T)
= \langle\Phi(T)\hat{1},\hat{1} \rangle$ defines a state on $B(H)$
which both extends $\tau$ and is invariant under the action of the
unitary group of $A$ on $B(H)$ and hence $\tau$ is an amenable trace.
\end{proof}

Though the previous theorem will get the most use we should also point
out a few other characterizations of amenable traces.

\begin{thm}
\label{thm:amenabletracesII}
Let $A \subset B(H)$ be in general position (i.e.\ $A \cap
\mathcal{K}(H) = \{0 \}$) and $\tau \in \TA$ be a tracial state.  Then
the following are equivalent:
\begin{enumerate}
\item $\tau$ is an amenable trace.

\item  There exists a $*$-monomorphism $\sigma : \pi_{\tau}
(A)^{\prime\prime} \hookrightarrow R^{\omega}$ such that
$\tau^{\prime\prime} = \tau_{\omega}\circ\sigma$ and $\sigma \circ
\pi_{\tau} : A \to R^{\omega}$ can be lifted to a u.c.p.\ map $\Phi:A \to
l^{\infty}(R)$, where $\tau^{\prime\prime}$ is the vector trace on
$\pi_{\tau} (A)^{\prime\prime}$ induced by $\tau$. (Kirchberg called
these ``liftable'' traces
in \cite{kirchberg:propertyTgroups}.)

\item There exist finite rank projections $P_n \in B(H)$ (not
necessarily increasing) such that $$\frac{\| aP_n - P_n a \|_{HS}}{\|
P_n \|_{HS}} \to 0 \ {\rm and \ } \tau(a) = \lim_{n \to \infty}
\frac{\langle aP_n, P_n\rangle_{HS}}{\langle P_n, P_n\rangle_{HS}},$$
for all $a \in A.$ (This is Connes' ``F{\o}lner'' condition which he
used to characterize the hyperfinite II$_1$-factor
\cite{connes:classification}.)
\end{enumerate}
\end{thm}

\begin{proof} We first show that the second statement above is
equivalent to the finite dimensional approximation property which is
statement (2) in Theorem \ref{thm:amenabletraces} and then we observe
that statement (3) above is also equivalent to the finite dimensional
approximation property.

(1) $\Longrightarrow$ (2).\footnote{Though this is implication is
well known to the experts and fairly straightforward, the reader
unfamiliar with this type of argument is advised to nail down
every detail as we will see several variations on this argument
later on.} Assume there exists a sequence of u.c.p.\ maps $\phi_n
: A \to M_{k(n)} (\mathbb C)$ such that $\| \phi_n (ab) - \phi_n
(a) \phi_n (b) \|_{2} \to 0$ and $\tau (a) = \lim_{n \to \infty}
\mathrm{tr}_{k(n)} \circ \phi_n (a)$, for all $a,b \in A$. Regard
each matrix algebra $M_{k(n)} (\mathbb C)$ as a subfactor of $R$
and consider the direct sum u.c.p.\ map $\Phi:A \to l^{\infty}(R)$
given by $$\Phi(a) = (\phi_n(a)).$$ It is not hard to see that
composing with the quotient map $l^{\infty}(R) \to R^{\omega}$
yields a $\tau$ preserving $*$-homomorphism $A \to R^{\omega}$
(with u.c.p.\ lifting $\Phi$). Finally one checks (essentially due
to uniqueness of GNS representations) that the weak closure (in
$R^{\omega}$) of this $*$-homomorphism is isomorphic to
$\pi_{\tau}(A)^{\prime\prime}$.

(2) $\Longrightarrow$ (1). For this direction we have a couple of
options.  We prefer the approximation ideas (and these ideas will
keep resurfacing throughout the paper) so our proof will use that
technology.\footnote{A more elegant approach, left to the reader,
would be to use statement (5) of Theorem \ref{thm:amenabletraces}.
The idea is that if $A \subset B(H)$ and $\Phi:A\to l^{\infty}(R)$
is the assumed u.c.p.\ lifting then we may extend $\Phi$ to all
$B(H)$, by the injectivity of $l^{\infty}(R)$, then pass to the
quotient $R^{\omega}$ and, finally, compose with a conditional
expectation $R^{\omega} \to
\sigma(\pi_{\tau}(A)^{\prime\prime}$).} So assume the existence of
a u.c.p.\ map $\Phi:A \to l^{\infty}(R)$ such that composition
with the quotient map $l^{\infty}(R) \to R^{\omega}$ yields a
$\tau$-preserving $*$-homomorphism $A \to R^{\omega}$.  It
suffices to show that if a finite set $\mathfrak{F} \subset A$ and
$\epsilon > 0$ are given then we can find a u.c.p.\ map $\phi:A
\to M_n(\mathbb{C})$ which is $\epsilon$-multiplicative (in
2-norm) and almost recovers $\tau$ on the finite set
$\mathfrak{F}$.

Let $\Phi_n:A\to R$ be the map $\Phi:A \to l^{\infty}(R) =
\Pi_{\mathbb{N}}$ composed with the projection map from
$l^{\infty}(R)$ onto the $n^{th}$ summand.  Since $\Phi$ is
multiplicative modulo the ideal $$I_{\omega} = \{ (x_n): \lim_{n
\to \omega} \|x_n\|_2 = 0\}$$ it follows that $$\|\Phi_n(ab) -
\Phi_n(a) \Phi_n(b) \|_2 \to 0 \ \mathrm{as} \ n \to \omega.$$
Hence for large enough $n$ the maps $\Phi_n$ are almost 2-norm
multiplicative on $\mathfrak{F}$ and recover $\tau$.  The proof is
then completed by replacing $R$ with a finite dimensional matrix
algebra which is possible since $R$ is hyperfinite.

That statement (3) above implies $\tau$ is amenable is quite
simple since cutting by the projections $P_n$ yields u.c.p.\ maps
to matrix algebras which are asymptotically multiplicative and
asymptotically recover $\tau$.  The converse is a consequence of
Voiculescu's Theorem (version \ref{thm:Voiculescu'sthm}) since we
assumed $A$ is in general position.
\end{proof}

\section{Uniform amenable traces}
\label{thm:uatraces}

Having seen that amenable traces are characterized by a finite
dimensional approximation property we can now define smaller subsets
of traces by demanding more of the approximants.

\begin{defn}
A trace $\tau$ will be called {\em uniform amenable} if there
exists a sequence of u.c.p.\ maps $\phi_n : A \to M_{k(n)}
(\mathbb C)$ such that $\| \phi_n (ab) - \phi_n (a) \phi_n (b)
\|_{2} \to 0$, for all $a,b \in A$, and $$\| \tau -
\mathrm{tr}_{k(n)}\circ \phi_n \|_{A^*} \to 0$$ where $\| \cdot
\|_{A^*}$ is the natural norm on the dual Banach space of $A$. The
set of all such traces will be denoted $\TAuim$.
\end{defn}

Similar to the case of amenable traces, the space $\TAuim$ admits a
number of nice characterizations.  We thank Yasuyuki Kawahigashi,
Sergei Neshveyev and Narutaka Ozawa for discussions in Oberwolfach
which  shortened our original proof of (1) $\Longrightarrow$
(2).  Ozawa also added (6) to the list below and showed us the elegant
proof.

\begin{thm}
\label{thm:mainthmUTAwafd}

For a trace $\tau \in \TA$, the following are equivalent:

\begin{enumerate}
\item $\tau \in \TAuim$.

\item There exist u.c.p.\ maps $\psi_n : A^{**} \to M_{k(n)} ({\mathbb
C})$ such that for each free ultrafilter $\omega \in \beta{\mathbb N}
\backslash {\mathbb N}$ we have $$\lim\limits_{n \to \omega} \| \psi_n
(xy) - \psi_n(x) \psi_n (y) \|_{2} = 0 \ {\rm and} \ \lim\limits_{n
\to \omega} \mathrm{tr}_{k(n)} \circ \psi_n (x) = \tau^{**} (x),$$
for all $x,y \in A^{**}$, where $\tau^{**}$ is the normal trace on
$A^{**}$ induced by $\tau$.

\item There exist u.c.p.\ maps $\psi_n : \pi_{\tau} (A)^{\prime\prime}
\to M_{k(n)} ({\mathbb C})$ such that for each free ultrafilter
$\omega \in \beta{\mathbb N} \backslash {\mathbb N}$ we have
$$\lim\limits_{n \to \omega} \| \psi_n (xy) - \psi_n(x) \psi_n (y)
\|_{2} = 0 \ {\rm and} \ \lim\limits_{n \to \omega} \mathrm{tr}_{k(n)}
\circ \psi_n (x) = \tau^{\prime\prime} (x),$$ for all $x,y \in
\pi_{\tau} (A)^{\prime\prime}$, where $\tau^{\prime\prime}$ is the
normal trace on $\pi_{\tau} (A)^{\prime\prime} $ induced by $\tau$.

\item There exists a u.c.p.\ liftable, normal $*$-monomorphism $\sigma
: \pi_{\tau} (A)^{\prime\prime} \hookrightarrow R^{\omega}$.

\item $\pi_{\tau} (A)^{\prime\prime}$ is hyperfinite.

\item $\pi_{\tau} : A \to \pi_{\tau} (A)^{\prime\prime}$ is weakly
nuclear (i.e.\ there exists u.c.p.\ maps $\phi_n : A \to M_{k(n)}
({\mathbb C})$, $\psi_n : M_{k(n)} ({\mathbb C}) \to \pi_{\tau}
(A)^{\prime\prime}$ such that $\psi_n \circ \phi_n (a) \to \pi_{\tau}
(a)$ in the $\sigma$-weak topology for every $a \in A$).
\end{enumerate}
\end{thm}

\begin{proof}

(1) $\Longrightarrow$ (2). Let $\phi_n : A \to M_{k(n)} ({\mathbb
C})$ be 2-norm asymptotically multiplicative u.c.p.\ maps such
that $\mathrm{tr}_{k(n)} \circ \phi_n \to \tau$ in the norm on
$A^*$.  Let $\phi_n^{**} : A^{**} \to M_{k(n)} ({\mathbb C})$ be
the canonical normal extensions to the double dual.  As is well
known, and easily checked, $\phi_n^{**}$ are also u.c.p.\ maps.
Hence $\mathrm{tr}_{k(n)} \circ \phi_n^{**}$ is a normal state on
$A^{**}$ and is in fact equal to the normal state on $A^{**}$
induced by the functional $\mathrm{tr}_{k(n)} \circ \phi_n$ on $A$
(they are both normal and agree on $A$).  Since $\|
\mathrm{tr}_{k(n)} \circ \phi_n - \tau \|_{A^*} \to 0$, it follows
that $\mathrm{tr}_{k(n)} \circ \phi_n^{**} (x) \to \tau^{**}(x)$
for every $x \in A^{**}$.  Hence, as in the proof of (1)
$\Longrightarrow$ (2) from Theorem \ref{thm:amenabletracesII},  we
can construct a u.c.p.\ map $\Phi : A^{**} \to R^{\omega}$ such
that $i)$ $\tau_{\omega} \circ \Phi = \tau^{**}$, $ii)$
$\Phi|_{A}$ is a $*$-homomorphism and $iii)$ $\Phi$ has a u.c.p.\
lifting $A^{**} \to l^{\infty} (R)$.\footnote{As we have seen
before, one first identifies each $M_{k(n)} ({\mathbb C})$ with a
subfactor of $R$ and then considers the direct sum map $A^{**} \to
l^{\infty}(R)$, $x \mapsto (\phi_n^{**} (x))$, and finally
composes with the quotient map to $R^{\omega}$.} If we knew that
$\Phi$ was a homomorphism then it would follow that
$\lim\limits_{n \to \omega} \| \phi_n^{**} (xy) - \phi_n^{**} (x)
\phi_n^{**} (y) \|_{2} = 0$ for all $x,y \in A^{**}$ and hence
this is what we will show.

First note that $\Phi$ is normal: if $\{ x_{\lambda} \} \subset
A^{**}_{sa}$ is a norm bounded, increasing net of self adjoint
elements with strong operator topology limit $x$ then $\{
\Phi(x_{\lambda})\}$ is increasing up to $\Phi(x)$ (in the strong
operator topology -- i.e.\ 2-norm) since $\Phi(x_{\lambda}) \leq
\Phi(x)$ and $\tau_{\omega}\circ (\Phi(x_{\lambda})) = \tau^{**}
(x_{\lambda}) \to \tau^{**} (x) = \tau_{\omega} \circ \Phi(x)$. It
follows that $\Phi$ is continuous from the $\sigma$-weak topology
on $A^{**}$ to the $\sigma$-weak topology on $R^{\omega}$ (i.e.\
with respect to the weak-$*$ topologies coming from the preduals).
Letting $\Psi : A^{**} \to R^{\omega}$ be the (weak-$*$
continuous) $*$-homomorphism which extends $\Phi|_{A}$ (and which
exists by universality of $A^{**}$) it follows that $\Phi = \Psi$
since they are continuous and agree on $A$.  Hence $\Phi$ is also
multiplicative.

(2) $\Longrightarrow$ (3).  Assuming (2), we can use the maps $\psi_n$
to construct a u.c.p.\ liftable, $*$-homomorphism $\sigma : A^{**} \to
R^{\omega}$ such that $\tau_{\omega} \circ \sigma = \tau^{**}$.  It
follows that $\sigma(A^{**}) \cong \pi_{\tau} (A)^{\prime\prime}$ and,
hence, $A^{**} \cong ker(\sigma) \oplus \pi_{\tau}
(A)^{\prime\prime}$.  Restricting the maps $\psi_n$ to this non-unital
copy of $\pi_{\tau} (A)^{\prime\prime}$ gives c.p.\ maps with the
desired properties.  Then applying Lemma \ref{thm:nonunital} we can
replace these nonunital maps with unital ones and we get (3).

(3) $\Longrightarrow$ (4) is similar to the proof of (1)
$\Longrightarrow$ (2) from Theorem \ref{thm:amenabletracesII}.

(4) $\Longrightarrow$ (5).  Identify $\pi_{\tau} (A)^{\prime\prime}$
    with $\sigma( \pi_{\tau} (A)^{\prime\prime}) \subset R^{\omega}$
    and let $\Phi:\pi_{\tau} (A)^{\prime\prime} \to l^{\infty}(R)$ be
    a u.c.p.\ splitting and $E:R^{\omega} \to \pi_{\tau}
    (A)^{\prime\prime}$ be any conditional expectation.  Now assume
    that $X \subset Y$ is an inclusion of operator systems and $\phi:X
    \to \pi_{\tau} (A)^{\prime\prime}$ is a u.c.p.\ map.  Since
    $l^{\infty}(R)$ is injective we may extend the map
    $\Phi\circ\phi:X \to l^{\infty}(R)$ to all of $Y$.  Composing this
    extension with the quotient map $l^{\infty}(R) \to R^{\omega}$
    followed by the conditional expectation $E:R^{\omega} \to
    \pi_{\tau} (A)^{\prime\prime}$ yields a u.c.p.\ map $Y \to
    \pi_{\tau} (A)^{\prime\prime}$.  This map extends $\phi$ since
    $\Phi$ was a lifting.

(5) $\Longrightarrow$ (1) is a simple consequence of Lemmas
\ref{thm:tracepreservingconditionalexpectation} and
\ref{thm:tracepreservingembedding} since hyperfinite von Neumann
algebras contain weakly dense, finite dimensional subalgebras.

At this point we have shown that (1) - (5) are equivalent.

(5) $\Longrightarrow$ (6) is trivial\footnote{Use the fact that finite
dimensional subalgebras are weakly dense in hyperfinite von Neumann
algebras.  Perhaps we should also remark that there is no difference
between only considering matrix algebras, as opposed to general finite
dimensional algebras, since every finite dimensional algebra unitally
embeds into a matrix algebra and there will always be a conditional
expectation from the larger matrix algebra back onto the original
finite dimensional algebra.} and hence we are left to prove (6)
$\Longrightarrow$ (5).  So assume that $\pi_{\tau} : A \to \pi_{\tau}
(A)^{\prime\prime}$ is weakly nuclear and $\phi_n : A \to M_{k(n)}
({\mathbb C})$, $\psi_n : M_{k(n)} ({\mathbb C}) \to \pi_{\tau}
(A)^{\prime\prime}$ are u.c.p.\ maps whose composition converges to
$\pi_{\tau}$ in the point-$\sigma$-weak topology. Using these maps it
is not hard to see that the canonical homomorphism $A \odot \pi_{\tau}
(A)^{\prime} \to B(L^2(A,\tau))$, $a\otimes x \mapsto \pi_{\tau}
(a)x$, is continuous with respect to the minimal tensor product
norm. (Use the fact that the natural map on the maximal tensor product
approximately factorizes through $M_{k(n)} \otimes_{max} \pi_{\tau}
(A)^{\prime} = M_{k(n)} \otimes \pi_{\tau} (A)^{\prime}$ and hence
factors through the minimal tensor product.)  As in the proof of (4)
$\Longrightarrow$ (5) from Theorem \ref{thm:amenabletraces}, it
follows that there exists a conditional expectation $B(L^2(A,\tau))
\to \pi_{\tau} (A)^{\prime}$ and hence $\pi_{\tau} (A)^{\prime}$ is
injective.  This implies that $\pi_{\tau} (A)^{\prime\prime}$ is also
injective and the proof is complete.
\end{proof}

Note that part (3) in the previous theorem could be used as an
abstract (i.e.\ representation free) definition of
quasidiagonality, analogous to Voiculescu's abstract
characterization, in the setting of tracial von Neumann algebras.
The equivalence of (3) and (5) would then say that
quasidiagonality is equivalent to hyperfiniteness. Note that this
is in stark contrast to the C$^*$-case where Dadarlat has
constructed non-nuclear Popa -- hence quasidiagonal -- algebras
\cite{dadarlat:nonnuclearsubalgebras}.

\begin{rem}
Though we prefer to state all approximation properties in this
paper in terms of matrix algebras there is one possible advantage
to using general finite dimensional approximants.  Namely, if
$\tau$ is a uniform amenable trace on $A$ then there exist finite
dimensional C$^*$-algebras $B_n$, traces $\tau_n \in
\mathrm{T}(B_n)$ and u.c.p.\ maps $\phi_n:A\to B_n$ such that
$\|\phi_n(ab) - \phi_n(a)\phi_n(b) \|_{2,\tau_n} \to 0$ for all
$a,b \in A$ and (here is the advantage)
$$\tau = \tau_n\circ\phi_n$$ for all $n \in {\mathbb N}$.  The proof
of this uses the fact that $\pi_{\tau} (A)^{\prime\prime}$ is
hyperfinite and hence there are $\tau$-preserving conditional
expectations onto the finite dimensional subalgebras of $\pi_{\tau}
(A)^{\prime\prime}$.
\end{rem}

\section{Quasidiagonal traces}
\label{thm:qdtraces}

In the previous section we strengthened the approximation property
of amenable traces by requiring norm convergence, as opposed to
weak-$*$ convergence, to the trace $\tau$.  We now define a subset of
traces by requiring a stronger notion of asymptotic multiplicativity.
Since we have two choices for convergence to $\tau$ this actually
leads to two subsets of $\TAim$.

\begin{defn}
We say a trace $\tau \in \TA$ is {\em quasidiagonal} if there exist
u.c.p.\ maps $\phi_n:A\to M_{k(n)}({\mathbb C})$ such that
$\mathrm{tr}_{k(n)}\circ \phi_n \to \tau$ in the weak-$*$ topology and
$$\|\phi_n(ab) - \phi_n(a)\phi_n(b)\| \to 0$$ for all $a,b \in
A$.\footnote{This terminology is inspired by Voiculescu's abstract
characterization of quasidiagonal C$^*$-algebras: $A$ is
quasidiagonal if and only there exist u.c.p.\ maps $\phi_n:A\to
M_{k(n)}({\mathbb C})$ such that $\|\phi_n(ab) -
\phi_n(a)\phi_n(b)\| \to 0$ and $\|a\| = \lim \|\phi_n(a)\|$ for
all $a,b \in A$.}  A trace $\tau \in \TA$ will be called {\em
uniform quasidiagonal} if one can further arrange that $$\| \tau -
\mathrm{tr}_{k(n)}\circ \phi_n\|_{A^*} \to 0.$$ These sets of
traces will be denoted $\TAQD$ and $\TAuQD$, respectively.
\end{defn}

These sets are far more mysterious than the amenable
traces.\footnote{But also more important to the classification
program.}  For example, we don't know whether every amenable trace
is quasidiagonal. Though we doubt that this is the case, at the
moment we have no good way of distinguishing the two sets.
(Compare with Blackadar and Kirchberg's question
\cite{blackadar-kirchberg}: Is every nuclear, stably finite
C$^*$-algebra necessarily quasidiagonal?) If one starts with a
quasidiagonal C$^*$-algebra, however, then it turns out that
$\TAQD$ is {\em precisely} the set of traces which can be encoded
in the definition of quasidiagonality. Recall that, by definition,
a C$^*$-algebra $A$ is quasidiagonal if there exists a faithful
representation $\pi : A \to B(H)$ such that one can find an
increasing sequence of finite rank projections $P_1 \leq P_2 \leq
\ldots$ with the property that $\| \pi(a)P_n - P_n \pi(a) \| \to
0$ for all $a \in A$ and $P_n \to 1_H$ in the strong operator
topology. (This is not the right definition for non-separable
algebras.)

\begin{prop}
\label{thm:QDcase} Let $A \subset B(H)$ be in general position
(i.e.\ $A \cap {\mathcal K}(H) = 0$).  If $A$ is quasidiagonal
then $\TAQD \neq \emptyset$ and, moreover, there exists an
increasing sequence of finite rank projections $P_1 \leq P_2 \leq
\ldots$, converging strongly to the identity, which asymptotically
commutes (in norm) with every element in $A$ and such that for
each $\tau \in \TAQD$ there exists a subsequence $\{ n_k \}$ such
that $$\frac{\langle aP_{n_k}, P_{n_k}\rangle_{HS}}{\langle
P_{n_k}, P_{n_k}\rangle_{HS}} \to \tau(a), \ as \ k \to \infty,$$
for all $a \in A$.
\end{prop}

\begin{proof} That $\TAQD \neq \emptyset$ is well known (cf.\
\cite[2.4]{dvv:QDsurvey}, \cite[Proposition 6.1]{brown:QDsurvey}),
but we remind the reader of the proof.  If $\pi : A \to B(H)$ and
$P_1 \leq P_2 \leq \ldots$ are as in the definition of
quasidiagonality then one defines u.c.p.\ maps by $$\phi_n(a) =
P_n \pi(a) P_n.$$ The asymptotic commutativity of $P_n$ ensures
that $\phi_n$ are asymptotically multiplicative in norm.  Note
also that $P_n B(H) P_n \cong M_{rank(P_n)} (\mathbb C)$. Finally,
a straightforward calculation shows that any weak-$*$ cluster
point of the sequence of states $\{ \mathrm{tr}_{rank(P_n)} \circ
\phi_n \}$ is necessarily a tracial state and hence $\TAQD \neq
\emptyset$.  (Note that we used the assumption that $A$ is unital
here; in the non-unital case it can happen that $\TA = \emptyset$
-- e.g.\ the suspension of a Cuntz algebra).

To prove the rest of the proposition we claim that it suffices to show
that for every finite set $\mathfrak{F} \subset A$, finite dimensional
subspace $X \subset H$, $\varepsilon > 0$ and trace $\tau \in \TAQD$
there exists a finite rank projection $P \in B(H)$ such that
\begin{enumerate}
\item $\| [a,P] \| < \varepsilon$ for all $a \in \mathfrak{F}$.

\item $P(x) = x$ for all $x \in X$.

\item $| \frac{\langle aP, P\rangle_{HS}}{\langle P, P\rangle_{HS}} -
\tau(a) | < \varepsilon$ for all $a \in \mathfrak{F}$.
\end{enumerate}

Assume for the moment that we were able to prove this local version.
Then, if $\{ a_n \} \subset A$ is a sequence which is dense in the
unit ball of $A$ and $\{ \tau_j \}$ is any sequence of traces in
$\TAQD$ we could apply the above local approximation property to
construct a sequence $P_1 \leq P_2 \leq \ldots$ which was converging
strongly to the identity, asymptotically commuting in norm with $A$
and such that $$| \frac{\langle a_i P_n, P_n\rangle_{HS}}{\langle P_n,
P_n\rangle_{HS}} - \tau_n(a_i) | < 1/n$$ for all $n \in {\mathbb N}$
and $1 \leq i \leq n$.  Since $A$ is separable, the weak-$*$ topology
on $\TA$ is metrizable and hence we can always find a sequence of
traces $\{ \tau_j \} \subset \TAQD$ such that there exists a weak-$*$
dense subset $Y \subset \TAQD$ with the property that every element of
$Y$ appears infinitely many times in the sequence $\{ \tau_j \}$.  The
sequence of projections associated with such a sequence of traces will
have all the properties asserted in the statement of the
proposition. Hence it suffices to prove the local statement in the
first paragraph of the proof.

The required local statement is now a consequence of Voiculescu's
Theorem (version \ref{thm:technicalVoiculescuThm}) and a little
trickery.  Let $\tau \in \TAQD$ be arbitrary. Since $A$ is
quasidiagonal, by Lemma \ref{thm:trickery} we can find a sequence
of u.c.p.\ maps $\phi_n : A \to M_{k(n)}({\mathbb C})$ which are
asymptotically multiplicative, asymptotically isometric and such
that $\mathrm{tr}_{k(n)} \circ \phi_n \to \tau$ in the weak-$*$
topology.  If a finite set $\mathfrak{F} \subset A$ and $\epsilon
> 0$ are given then, by passing to a subsequence if necessary, we
may assume that $\| \phi_n (ab) - \phi_n(a)\phi_n(b) \| <
\varepsilon$ and $|\mathrm{tr}_{k(n)} \circ \phi_n (a) - \tau(a) |
< \varepsilon$ for all $n$ and for all $a,b \in \mathfrak{F}$.
Letting $K = \oplus_n {\mathbb C}^{k(n)}$ and $\Phi = \oplus_n
\phi_n : A \to B(K)$ we have that $\Phi$ is a faithful
$*$-homomorphism modulo the compacts.  Hence we can find a unitary
operator $U : K \to H$ such that $U \Phi(a) U^*$ is nearly equal
(in norm) to $a$, for all $a \in \mathfrak{F}$.  Hence, if $Q_s
\in B(K)$ is the orthogonal projection onto $\oplus_1^s {\mathbb
C}^{k(n)}$ we have that $\| [UQ_s U^*, a]\|$ is small for all $s
\in {\mathbb N}$ and for all $a \in A$.  Moreover, compressing $a
\in \mathfrak{F}$ to the range of $UQ_s U^*$ will almost recover
the trace $\tau$ (for all $s$).  Finally, if $X \subset H$ is any
finite dimensional subspace then $X$ will almost be contained in
the range of $UQ_s U^*$ for sufficiently large $s$ and hence a
tiny (norm) perturbation of a sufficiently large $UQ_s U^*$ will
actually be the identity on $X$ (and still almost commute with
$\mathfrak{F}$ and almost recover the trace $\tau$).
\end{proof}

It is natural to wonder whether or not every trace on a
quasidiagonal C$^*$-algebra arises as above; i.e.\ whether or not
$\TA = \TAQD$ when $A$ is quasidiagonal.  This is the case for $A
= C^*({\mathbb F}_{\infty})$ (the universal C$^*$-algebra of a
countably generated free group) if and only if Connes' embedding
problem has an affirmative answer (see Proposition
\ref{thm:BFDII}). Unfortunately this is not the case in general;
there exist (exact) quasidiagonal C$^*$-algebras such that $\TA
\neq \TAQD$ (see Corollary \ref{thm:exactQD}).

\section{Locally finite dimensional traces}
\label{thm:lfdtraces}

There is at least one more natural set of traces which one may define
-- those arising as limits of finite dimensional traces (i.e.\ limits
of traces whose GNS representation is finite dimensional).  Though
very natural, this set of traces would be empty for every infinite
dimensional, simple C$^*$-algebra (or any other C$^*$-algebra which
has no finite dimensional representations).  However it turns out that
there is another possibility for strengthening the definition of
quasidiagonal trace which contains a number of non-trivial examples,
though it is more technical to formulate.  Recall that if $\phi:A\to
B$ is a u.c.p.\ map then we denote the multiplicative domain of $\phi$
by $A_{\phi}$.

\begin{defn}
A trace $\tau$ will be called {\em locally finite dimensional} if
there exist u.c.p.\ maps $\phi_n:A\to M_{k(n)}({\mathbb C})$ such that
$\mathrm{tr}_{k(n)}\circ \phi_n \to \tau$ in the weak-$*$ topology and
$$d(a,A_{\phi_n}) \to 0$$ for all $a \in A$.\footnote{$d(a,A_{\phi_n})
\to 0$ means there exist $a_n \in A_{\phi_n}$ such that $\|a -
a_n\| \to 0$.} A trace $\tau \in \TA$ will be called {\em uniform
locally finite dimensional} if one can further arrange that $$\|
\tau - \mathrm{tr}_{k(n)}\circ \phi_n\|_{A^*} \to 0.$$ These sets
of traces will be denoted $\TAlfd$ and $\TAulfd$, respectively.
\end{defn}

Note that being locally finite dimensional is stronger than
quasidiagonality of a trace.  Indeed, if $\tau$ is locally finite
dimensional, with $\phi_n:A\to M_{k(n)}({\mathbb C})$ as in the
definition, and $a,b \in A$ are given then we have
$$\phi_n(ab) \approx \phi_n(a_n b_n) = \phi_n(a_n)\phi_n(b_n)
\approx \phi_n(a)\phi_n(b),$$ where $a_n, b_n \in A_{\phi_n}$ are
chosen so that $\|a - a_n \|, \|b - b_n\| \to 0$.

At this point it is hard to motivate this definition but the
interested reader may find it fun to prove that every trace on an
AF algebra is uniform locally finite dimensional. (Even though
there are plenty of simple, infinite dimensional AF algebras where
limits of finite dimensional traces would be impossible.) In fact,
we will see many more examples as well as observe that these
traces are intimately related to future progress in the
classification program. However, we having nothing more to say
about this definition at the moment.

\section{Miscellaneous remarks and permanence properties}

A few remarks about all of these sets are probably in order.
First note that we have the following inclusions:

$$\begin{array}{ccccccc}
\TA & \supset & \TAim & \supset & \TAQD & \supset & \TAlfd\\
& & \cup & & \cup & & \cup\\
& & \TAuim & \supset & \TAuQD & \supset & \TAulfd
\end{array}$$

A natural question to ask is whether or not any of these
inclusions are proper.  We will see that all of the vertical
inclusions are proper in general, but for large classes of
C$^*$-algebras (e.g.\ exact C$^*$-algebras) one has the equalities
$\TAim = \TAuim$ and $\TAQD = \TAuQD$ (see Theorem
\ref{thm:locallyreflexive}).\footnote{We don't know if every exact
C$^*$-algebra satisfies the equation $\TAlfd = \TAulfd$ and this
turns out to be a very important open question (see Section
\ref{thm:classificationprogram}).} While it is not hard to give
examples showing $\TA \neq \TAim$ (see Proposition
\ref{thm:reducedgroupalgebras}) we have not yet found an example
which shows that the set of amenable traces can be strictly larger
than the set of quasidiagonal traces.\footnote{Even for nuclear
C$^*$-algebras $A$ it appears to be very difficult to decide
whether or not one always has $\TAim = \TAQD$. If every nuclear
C$^*$-algebra satisfies the equation $\TAim = \TAQD$ then every
nuclear C$^*$-algebra with a faithful trace would be quasidiagonal
-- in particular, Rosenberg's conjecture that $C^*_r(\Gamma)$ is
quasidiagonal, for every discrete amenable group $\Gamma$, would
follow.  On the other hand, nuclear examples showing $\TAim \neq
\TAQD$ would likely have significant consequences for Elliott's
classification program.} We believe such examples should exist but
it is not clear at this point where to look. Similarly, the
relation between the quasidiagonal and the locally finite
dimensional traces is not yet well understood.

Another natural question is whether or not any of these sets are convex
and/or weak-$*$ closed subsets of $\TA$.

\begin{prop}
\label{thm:convexity} The sets $\TAim$, $\TAQD$, $\TAuim$ and
$\TAuQD$ are all convex. Both $\TAim$ and $\TAQD$ are closed in
the weak-$*$ topology while the sets $\TAuim$ and $\TAuQD$ are
closed in norm (i.e.\ the norm on $A^*$) and thus, by the
Hahn-Banach theorem, closed in the weak topology coming from
$A^{**}$.
\end{prop}

\begin{proof} Since $\TAQD$ and
$\TAim$ (resp.\ $\TAuim$ and $\TAuQD$) are defined via weak-$*$
convergence (resp.\ norm convergence), it is easy to see that
these sets are closed in this topology.

The proof of the convexity assertion is a reasonable exercise.
However, we will need little perturbations of this argument
several times so perhaps it is worthwhile pointing out the main
things to keep in mind for future arguments.\footnote{One of the
irritating aspects of our various definitions of approximation is
that we have required the use of matrix algebras as opposed to
general finite dimensional C$^*$-algebras.  This is mostly because
it makes the paper easier to write!  There is no loss of
generality in requiring matrix algebras, but it is sometimes
convenient in proofs to use general finite dimensional algebras.
The argument we are about to sketch allows one to pass freely from
one setting to another.}  We only give the proof of convexity for
$\TAuim$ as all other cases use similar arguments.  So assume
$\tau_1, \tau_2 \in \TAuim$, $0 < \lambda < 1$ and $\tau = \lambda
\tau_1 + (1 - \lambda)\tau_2$.  If we fix a finite set
$\mathfrak{F} \subset A$ and $\epsilon > 0$ then it is possible to
find u.c.p.\ maps $\phi_i:A \to M_{n(i)}(\mathbb{C})$, $i = 1,2$,
such that $$\|\phi_i(ab) - \phi_i(a)\phi_i(b) \|_{2,
\mathrm{tr}_{n(i)}} < \epsilon$$ for all $a,b \in A$ and $$\|
\tau_i - \mathrm{tr}_{n(i)}\circ \phi_i \|_{A^*} < \epsilon.$$  We
now consider the finite dimensional algebra $$B =
M_{n(1)}(\mathbb{C}) \oplus M_{n(2)}(\mathbb{C})$$ and the tracial
state $$\gamma = \lambda\mathrm{tr}_{n(1)} \oplus (1 -
\lambda)\mathrm{tr}_{n(2)}.$$  Let $\phi = \phi_1 \oplus \phi_2:A
\to B$  and some straightforward calculations show that
$\|\phi(ab) - \phi(a)\phi(b) \|_{2,\gamma}^2$ is equal to 
$$\lambda\|\phi_1(ab) -
\phi_1(a)\phi_1(b)\|_{2,\mathrm{tr}_{n(1)}}^2 + (1-\lambda)
\|\phi_2(ab) - \phi_2(a)\phi_2(b)\|_{2,\mathrm{tr}_{n(2)}}^2$$ and
$$\| \tau - \gamma\circ\phi \|_{A^*} \leq \lambda \|\tau_1 -
\mathrm{tr}_{n(1)}\circ\phi_1 \|_{A*} + (1-\lambda)\|\tau_2 -
\mathrm{tr}_{n(2)}\circ\phi_2 \|_{A*}.$$  It follows that $\tau$
has the right approximation properties if we allow general finite
dimensional C$^*$-algebras.  Hence we can complete the proof by
appealing to Lemma \ref{thm:tracepreservingembedding} to replace
$(B,\gamma)$ with a full matrix algebra.
\end{proof}

Evidently $\TAlfd$ (resp.\ $\TAulfd$) is also weak-$*$ (resp.\
norm) closed but we don't know whether or not it is convex.

We will soon see that if $\Gamma$ is a non-amenable, residually
finite, discrete group then the canonical trace on $C^* (\Gamma)$
(which gives the left regular representation in the GNS
construction) is always a weak-$*$ limit of (uniform locally)
finite dimensional traces, but is not itself a uniform amenable
trace.  Thus, the sets $\TAuim$, $\TAuQD$ and $\TAulfd$ need not
be weak-$*$ closed in general.

Finally, one may wonder if any of these sets define a face in
$\TA$. We are not sure about the other four sets, but this is the
case for the (uniformly) amenable traces.

The proofs are a simple consequence of Theorems
\ref{thm:amenabletraces} and \ref{thm:mainthmUTAwafd} and the
following fact.

\begin{lem} Assume $\tau, \gamma \in \TA$ and there exists $0 <
s < 1$ such that $s\tau \leq \gamma$.\footnote{This means
$s\tau(a^*a) \leq \gamma(a^*a)$ for all $a \in A$.}  Then there
exists a projection $P \in \pi_{\gamma}(A)^{\prime}$ such that $a
\mapsto P\pi_{\gamma}(a)$ is unitarily equivalent to $\pi_{\tau}$.
\end{lem}

\begin{proof}  By \cite[Proposition 3.3.5]{pedersen:book} we can
find a positive element $y \in \pi_{\gamma}(A)^{\prime}$ such that
$$\tau(a) = \langle \pi_{\gamma}(a) y \hat{1}, \hat{1} \rangle,$$
for all $a \in A$.  Uniqueness of GNS representations implies that
the projection onto the closure of the subspace $\{
\pi_{\gamma}(a) y\hat{1}:a \in A\}$ will do the trick.
\end{proof}

That $\TAim$ is a face was first observed by Kirchberg.

\begin{prop}{\em (cf.\ \cite[Lemma
3.4]{kirchberg:propertyTgroups})} $\TAim$ is a face in $\TA$.
\end{prop}

\begin{proof}  If $\tau_1, \tau_2 \in \TA$, $0 < s < 1$ are given
and $\gamma = s\tau_1 + (1-s)\tau_2 \in \TAim$ then we can find
projections $P_i \in \pi_{\gamma}(A)^{\prime\prime}$ such that $a
\mapsto P_i \pi_{\gamma}(a)$ is unitarily equivalent to
$\pi_{\tau_i}$ and (by amenability) we can find a u.c.p.\ map
$\Phi:B(H) \to \pi_{\gamma}(A)^{\prime\prime}$ such that $\Phi(a)
= \pi_{\tau}(a)$ (for any fixed inclusion $A \subset B(H)$). Hence
we can reapply condition (5) in Theorem \ref{thm:amenabletraces}
to deduce amenability of each $\tau_i$ by first taking the u.c.p.\
map $\Phi$ and composing with compression by $P_i$ (identifying
the latter map with the GNS representation coming from $\tau_i$).
\end{proof}

\begin{prop}
$\TAuim$ is a face in $\TA$.
\end{prop}

\begin{proof}  The proof is very similar to the previous
proposition.   We use Theorem \ref{thm:mainthmUTAwafd} and the
fact that all von Neumann subalgebras of a finite, hyperfinite von
Neumann algebra are again hyperfinite.
\end{proof}

We now wish to discuss a few permanence properties of the traces
we have defined. We begin with two easy situations.  In the first
we consider restrictions of traces to subalgebras.  The proof of
the following proposition is a very simple exercise -- just take
the finite dimensional approximating maps which are given on the
larger algebra and restrict them to the smaller algebra.  Note
that this obvious procedure does not work when considering locally
finite dimensional traces as it is not clear that the
multiplicative domain of a u.c.p.\ map defined on the big algebra
would have to intersect the smaller algebra.\footnote{We are not
claiming that the restriction of a locally finite dimensional
trace to a subalgebra need not be locally finite dimensional!
Indeed, we don't know whether or not this is true.  We are only
claiming that the obvious proof doesn't seem to work.}

\begin{prop}[Restriction to subalgebras]
Assume $1_A \in B \subset A$ is a C$^*$-subalgebra and $\tau \in
\TA$.  If $\tau$ is amenable, uniform amenable, quasidiagonal or
uniform quasidiagonal then $\tau|_B$ enjoys the same approximation
property.
\end{prop}

Another trivial fact, again left to the reader, is that nice
traces defined on quotients always lift to nice traces.  In this
situation, even the locally finite dimensional traces behave well.

\begin{prop}[Lifting from a quotient] If $I \triangleleft A$ is
an ideal and $\tau \in \mathrm{T}(A/I)$ is a trace enjoying any
one of the six approximation properties defined in these notes
then the trace on $A$ gotten by composing $\tau$ with the quotient
map $A \to A/I$ inherits the same approximation property.
\end{prop}

A marginally less trivial situation is that of tensor products.

\begin{prop}
\label{thm:tensorproduct}  Assume $\tau \in \TA$ and $\gamma \in
\TB$.
\begin{enumerate}
\item If $\tau \in \TAim$ (resp.\ $\tau \in \TAuim$) and $\gamma
\in \TBim$ (resp.\ $\gamma \in \TBuim$) then $\tau \otimes \gamma
\in \mathrm{AT}(A\otimes B)$ (resp.\ $\tau \otimes \gamma \in
\mathrm{UAT}(A\otimes B)$).

\item If $\tau \in \TAQD$ (resp.\ $\tau \in \TAuQD$) and $\gamma
\in \TBQD$ (resp.\ $\gamma \in \TBuQD$) then $\tau \otimes \gamma
\in \mathrm{AT}(A\otimes B)_{QD}$ (resp.\ $\tau \otimes \gamma \in
\mathrm{UAT}(A\otimes B)_{QD}$).

\item If $\tau \in \TAlfd$ (resp.\ $\tau \in \TAulfd$) and $\gamma
\in \TBlfd$ (resp.\ $\gamma \in \TBulfd$) then $\tau \otimes
\gamma \in \mathrm{AT}(A\otimes B)_{LFD}$ (resp.\ $\tau \otimes
\gamma \in \mathrm{UAT}(A\otimes B)_{LFD}$).
\end{enumerate}
\end{prop}

\begin{proof}  We think the details of the proof are best left as
an exercise. One should recall, however, that one of the nice
things about minimal tensor products and completely bounded maps
is that if $\phi:A \to C$ and $\psi:B \to D$ are completely
bounded maps (e.g.\ linear functionals -- which, together with the
case of u.c.p.\ maps, is what one must handle in the proof of the
present proposition) then there is a well defined completely
bounded map $\phi\otimes \psi:A\otimes B \to C\otimes D$ such that
$$\| \phi \otimes \psi \|_{CB} = \|\phi\|_{CB} \| \psi \|_{CB}.$$
One can refer to any book on operator space theory for this
fundamental fact.

Given this fact all six cases are handled in the same fashion.
Take u.c.p.\ maps on $A$ and $B$ with appropriate approximation
properties and consider the tensor product of these maps.
\end{proof}

Determining when a trace with nice approximation properties which
happens to vanish on an ideal descends to a nice trace on the
quotient is sometimes easy, sometimes hard and sometimes not at
all clear.\footnote{For example, it is not clear what happens with
the locally finite dimensional traces.}

\begin{prop}[Passage to a quotient]
\label{thm:passagetoquotient} Let $\tau \in \TA$ be given, $I =
ker(\pi_{\tau}(A))$, $B = A/I \cong \pi_{\tau}(A)$ and $\dot{\tau}
\in \mathrm{T}(B)$ be the (faithful, vector) trace induced by
$\tau$.
\begin{enumerate}
\item If $\tau \in \TAim$ and there exists a contractive c.p.\
(c.c.p.) splitting $\Phi:B \to A$ of the quotient map $A \to B$
then $\dot{\tau} \in \mathrm{AT}(B)$.

\item If $\tau \in \TAuim$ then $\dot{\tau} \in \mathrm{UAT}(B)$.

\item If $\tau \in \TAQD$, there exists a c.c.p.\ splitting
$\Phi:B \to A$ of the quotient map $A \to B$ and the extension $$0
\to I \to A \to B \to 0$$ is quasidiagonal then $\dot{\tau} \in
\TBQD$.\footnote{An extension $0 \to I \to A \to B \to 0$ is
called quasidiagonal if there exists an approximate unit $\{p_n\}
\subset I$ which (a) is quasicentral in $A$ and (b) each $p_n$ is
a projection.}

\item If $\tau \in \TAuQD$, there exists a c.c.p.\ splitting
$\Phi:B \to A$ of the quotient map $A \to B$ and the extension $$0
\to I \to A \to B \to 0$$ is quasidiagonal then $\dot{\tau} \in
\TBuQD$.
\end{enumerate}
\end{prop}

\begin{proof}  Though none of these assertions are particularly
difficult note that only the second one is obvious since the fifth
statement in Theorem \ref{thm:mainthmUTAwafd} asserts that the
{\em image} of the GNS representation completely determines when a
trace is uniform amenable (and the images of the GNS
representations of $(A,\tau)$ and $(B, \dot{\tau})$ are
canonically unitarily equivalent).

For the proof of (1) we will apply the fifth condition from
Theorem \ref{thm:amenabletraces}.  Indeed, by assumption, if $A
\subset B(H)$ and $B \subset B(K)$ are faithful representations
then we can find a u.c.p.\ map $$\Psi: B(H) \to
\pi_{\tau}(A)^{\prime\prime} =
\pi_{\dot{\tau}}(B)^{\prime\prime}$$ such that $\Psi(a) =
\pi_{\tau}(a)$ for all $a \in A$.  By Arveson's extension theorem
we may assume that the u.c.p.\ splitting $\Phi:B \to A$ is
actually defined on all of $B(K)$ and takes values in $B(H)$.  One
readily verifies that the u.c.p.\ map $\Psi\circ\Phi:B(K) \to
\pi_{\dot{\tau}}(B)^{\prime\prime}$ also satisfies condition (5)
in Theorem \ref{thm:amenabletraces} and hence $\dot{\tau} \in
\mathrm{AT}(B)$.

The proofs of (3) and (4) are quite similar so we only give the
proof of (4). Since $0 \to I \to A \to B \to 0$ is a quasidiagonal
extension we can find a quasicentral approximate unit of
projections, say $\{p_n \} \subset I$.  If we consider the
sequence of c.p.\ splittings $$\Phi_n:B \to A, \ \Phi_n(b) = (1 -
p_n)\Phi(b)(1 - p_n)$$ then it is well known that these maps will
be asymptotically multiplicative in norm.\footnote{Since $\|
\Phi_n(xy) - \Phi_n(x)\Phi_n(y) \| \leq \| (1 - p_n)[\Phi(xy) -
\Phi(x)\Phi(y)](1 - p_n)\| +   \| [p_n,\Phi(x)]\|\|y\|$ and the
latter two norms tend to zero since $\{p_n\}$ is an approximate
unit and quasicentral, respectively.}  Hence if a finite set
$\mathfrak{F} \subset B$ and $\epsilon > 0$ are given then we can
find $n$ large enough that $\Phi_n$ is as close to multiplicative
on $\mathfrak{F}$ as one desires.  Then one can apply the
definition of uniform quasidiagonal trace to the finite set
$\Phi_n(\mathfrak{F}) \subset A$ to find a u.c.p.\ map $\phi:A \to
M_n(\mathbb{C})$ which is  almost multiplicative on
$\Phi_n(\mathfrak{F})$ and almost recaptures $\tau$ (in the
uniform norm) after composing with the trace on matrices.  It is a
simple exercise to check that $\phi\circ\Phi_n:B \to
M_n(\mathbb{C})$ is almost multiplicative on $\mathfrak{F}$ and
that $$\| \mathrm{tr}_n \circ\phi\circ\Phi_n - \dot{\tau} \|_{B^*}
\leq \| \mathrm{tr}_n \circ\phi - \tau \|_{A^*}.$$ Of course we
really should have unital maps instead of just contractive c.p.\
maps but Lemma \ref{thm:nonunital} allows us to circumvent this
technicality and thus we deduce that $\dot{\tau}$ is uniform
quasidiagonal on $B$.
\end{proof}

\begin{rem}  In some instances we don't actually need the full
power of a c.c.p.\ splitting $\Phi:B \to A$ -- it may suffice to
just have {\em local liftability} which asserts that for each
finite dimensional operator system $X \subset B$ there should
exist a u.c.p.\ map $\Phi_X:X \to A$ which lifts $X$.  For example
with this hypothesis the conclusions in (1) and (3) above still
hold with a little more work.  For (1) the idea of the proof is
the same except one has to take a limit of maps using the fact
that the set of u.c.p.\ maps between two von Neumann algebras is
compact in the point ultraweak topology.  For (3) the proof is
again quite similar but one must apply Arveson's extension theorem
at one point to extend maps which are only defined on finite
dimensional operator subsystems of $B$ to all of $B$.  Since we
won't need these results we leave the details to the interested
reader.  Finally we recall that the Effros-Haagerup lifting
theorem asserts that  local liftability is equivalent to knowing
that for every C$^*$-algebra $D$ the sequence $$0 \to I \otimes D
\to A \otimes D \to B \otimes D \to 0$$ is exact (cf.\
\cite{effros-haagerup}).
\end{rem}

Though we will try our hardest to avoid non-unital C$^*$-algebras
they are occasionally forced upon us.  Hence a few words about
unitizations and amenable traces on ideals are probably in order.
However, we don't want to spend time systematically studying this
nonunital setting so we only present the results that we will
need.

Assume that $\tau$ is a tracial state on a {\em nonunital}
C$^*$-algebra $A$ and let $\tilde{A}$ be the unitization of $A$
and $\tilde{\tau}$ denote the unique extension of $\tau$ to a
tracial state on $\tilde{A}$.  When trying to extend the six
definitions of traces given in the previous sections one is
faced with two natural possibilities.  On the one hand it would be
very natural to simply replace u.c.p.\ maps with c.c.p.\ maps,
mimic the previous approximation properties and take these as the
nonunital definitions. On the other hand it would also be natural
to just define your way out of this issue by declaring that $\tau$
is quasidiagonal, for example, if $\tilde{\tau}$ is a
quasidiagonal trace on $\tilde{A}$.  Not surprisingly, in all six
cases these two approaches yield the same sets of traces on $A$.

\begin{prop}[Traces on nonunital C$^*$-algebras] If $A$ is
nonunital, $\tau \in \TA$, $\tilde{A}$ is the unitization of $A$
and $\tilde{\tau} \in \mathrm{T}(\tilde{A})$ is the unique
extension then the following statements hold.
\begin{enumerate}
\item $\tilde{\tau}$ is amenable (resp.\ uniform amenable) if and
only if there exist c.c.p.\ maps $\phi_n:A \to
M_{k(n)}(\mathbb{C})$ such that $\|\phi_n(ab) - \phi_n(a)\phi_n(b)
\|_2 \to 0$ and $\mathrm{tr}_{k(n)}\circ \phi_n \to \tau$ weak$-*$
(resp.\ $\|\mathrm{tr}_{k(n)}\circ \phi_n - \tau\|_{A^*} \to 0$).

\item $\tilde{\tau}$ is quasidiagonal (resp.\ uniform
quasidiagonal) if and only if there exist c.c.p.\ maps $\phi_n:A
\to M_{k(n)}(\mathbb{C})$ such that $\|\phi_n(ab) -
\phi_n(a)\phi_n(b) \| \to 0$ and $\mathrm{tr}_{k(n)}\circ \phi_n
\to \tau$ weak$-*$ (resp.\ $\|\mathrm{tr}_{k(n)}\circ \phi_n -
\tau\|_{A^*} \to 0$).

\item $\tilde{\tau}$ is locally finite dimensional (resp.\ uniform
locally finite dimensional) if and only if there exist c.c.p.\
maps $\phi_n:A \to M_{k(n)}(\mathbb{C})$ such that $d(a,A_{\phi})
\to 0$ and $\mathrm{tr}_{k(n)}\circ \phi_n \to \tau$ weak$-*$
(resp.\ $\|\mathrm{tr}_{k(n)}\circ \phi_n - \tau\|_{A^*} \to 0$).
\end{enumerate}
\end{prop}

\begin{proof}  The key technical fact we will need is again due to
Choi-Effros (cf.\ \cite[Lemma 3.9]{choi-effros:lifting}): If
$\phi:A \to M_n(\mathbb{C})$ is c.c.p.\ then the unital linear map
$\tilde{\phi}:\tilde{A} \to M_n(\mathbb{C})$ defined by
$\tilde{\phi}(a + \lambda) = \phi(a) + \lambda 1$ is also
completely positive.  With this result in hand all six of the `if'
statements above are simple exercises (left to the reader).

The first four `only if' statements are completely trivial -- just
restrict the finite dimensional approximating maps on $\tilde{A}$
to $A$.  The locally finite dimensional case is only slightly
harder.  One must verify that the multiplicative domain of a map
defined on $\tilde{A}$ actually intersects $A$.  For general
subalgebras there is no reason that this should be true but
luckily $A$ is an ideal in $\tilde{A}$ and hence a routine
exercise shows that if $\mathfrak{F} \subset A$ is a given finite
set and $B \subset \tilde{A}$ is a subalgebra which nearly
contains $\mathfrak{F}$ then $B$ must intersect $A$ and this
intersection must almost contain $\mathfrak{F}$ as well.  Hence
the final two `only if' statements also follow by simply
restricting the maps which are given on $\tilde{A}$.
\end{proof}

Now that we know how to define amenable traces on nonunital
C$^*$-algebras we can prove a couple simple facts  which will be
needed later.

\begin{prop}[Restriction to an Ideal]
\label{thm:restrictiontoideal}
Let $I \triangleleft A$ be an
ideal, $\tau \in \TA$ and $\gamma = \frac{1}{\|\tau|_I\|}\tau \in
\mathrm{T}(I)$. If $\tau$ is amenable (resp.\ uniform amenable)
then $\gamma$ is also amenable (resp.\ uniform amenable).
\end{prop}

\begin{proof}  If the restriction of $\tau$ to $I$ happens to be a
state on $I$ then there is nothing to prove since we can just
restrict the finite dimensional approximating maps on $A$ to $I$
as in the last proposition.  However, that need not be the case so
a bit more care is required.

By uniqueness of GNS representations we may identify
$\pi_{\gamma}(I)^{\prime\prime}$ with the weak closure of
$\pi_{\tau}(I)$ inside
$\pi_{\tau}(A)^{\prime\prime}$.\footnote{Actually, since $I$ is an
ideal we get the decomposition $A^{**} = I^{**} \oplus (A/I)^{**}$
and hence we can identify $\pi_{\gamma}(I)^{\prime\prime}$ with a
direct summand of $\pi_{\tau}(A)^{\prime\prime}$.} In the uniform
amenable case the proof then follows from the fifth statement in
Theorem \ref{thm:mainthmUTAwafd} since hyperfiniteness passes to
all von Neumann subalgebras of a finite, hyperfinite von Neumann
algebra.

For the case of amenable traces one applies condition (5) from
Theorem \ref{thm:amenabletraces} as follows.  Let $A \subset B(H)$
be any faithful representation and $\Phi:A \to
\pi_{\tau}(A)^{\prime\prime}$ be a u.c.p.\ map such that $\Phi(a)
= \pi_{\tau}(a)$. Now decompose $\pi_{\tau}(A)^{\prime\prime}
\cong \pi_{\gamma}(I)^{\prime\prime} \oplus M$ and let $P\in
\pi_{\tau}(A)^{\prime\prime}$ be the (central) projection such
that $P\pi_{\tau}(A)^{\prime\prime} =
\pi_{\gamma}(I)^{\prime\prime}$.  Finally one defines the desired
u.c.p.\ map $\Psi: B(H) \to \pi_{\gamma}(I)^{\prime\prime}$ by
$\Psi(T) = P\Phi(T)$ and running back through the identifications
one sees that $\Psi(x) = \pi_{\gamma}(x)$ for all $x \in I$.
\end{proof}

\begin{prop}[Extending from an Ideal]
\label{thm:extensionfromideal} If $I \triangleleft A$ is an ideal
and $\tau \in \mathrm{T}(I)$ is an amenable (resp.\ uniform
amenable) trace then there is a (unique) trace $\gamma \in \TA$
which extends $\tau$ and it is amenable (resp.\ uniform amenable).
\end{prop}

\begin{proof} Since $I$ is an ideal there is a unique extension of
the GNS representation $\pi_{\tau}:I \to B(L^2(I,\tau))$ to a
$*$-homomorphism $\sigma:A \to \pi_{\tau}(I)^{\prime\prime}$.  Let
$\gamma$ be the unique tracial state extension of $\tau$ to $A$
(which is gotten by composing the GNS vector state on
$\pi_{\tau}(I)^{\prime\prime}$ with the map $\sigma$).  Since
$\sigma$ is unitarily equivalent to $\pi_{\gamma}:A\to
B(L^2(A,\gamma))$ it follows that $\gamma$ is uniform amenable
whenever $\tau$ is since hyperfiniteness is obviously passed to
$\pi_{\gamma}(A)^{\prime\prime}$.

Now assume that $\tau$ is only an amenable trace and fix a
faithful $*$-representation $A \subset B(H)$.  To see the essence
of the proof we will also further assume that $I$ is an essential
ideal in $A$ and hence belongs to the strict closure of $I$.  Now
if $\Phi:B(H) \to \pi_{\tau}(I)^{\prime\prime}$ is a u.c.p.\ map
such that $\Phi(x) = \pi_{\tau}(x)$ for all $x \in I$ then it
suffices to show that $\Phi(a) = \sigma(a)$ for all $a \in A$. But
for each $a \in A$ we can find a sequence $\{ x_n\} \subset I$
such that $yx_n \to ya$ and $x_n y \to ay$ (in norm) for all $y
\in I$.  Since $I$ falls in the multiplicative domain of $\Phi$ it
follows that $\pi_{\tau}(y)\pi_{\tau}(x_n) \to
\pi_{\tau}(y)\Phi(a)$ and $\pi_{\tau}(x_n)\pi_{\tau}(y) \to
\Phi(a)\pi_{\tau}(y)$ for all $y \in I$ and hence $\Phi(a)$ is the
strict limit of $\pi_{\tau}(x_n)$.  In other words, $\Phi$ is
automatically strictly continuous and hence must agree with
$\sigma$ (since $\sigma$ is also strictly continuous and agrees
with $\Phi$ on $I$).  This shows that $\gamma$ is an amenable
trace on $A$ as well.

The argument in the previous paragraph actually shows that
amenable traces always extend to amenable traces on the multiplier
algebra of a nonunital C$^*$-algebra.  Hence when $I \triangleleft
A$ is not an essential ideal we can uniquely extend an amenable
trace on $I$ to one on $A$ by first extending to the multiplier
algebra of $I$, denoted $M(I)$, and then composing with the
canonical $*$-homomorphism $A \to M(I)$.
\end{proof}

\chapter{Examples and special cases}

Having defined so many subsets of traces it may be
worthwhile to consider a number of examples.  Our hope is that a
detailed look at these examples will not only help cement the
definitions in the reader's mind but also begin to illustrate how
these traces are related to various other problems and
conjectures.

\section{Discrete groups}

As usual, groups provide some interesting and instructive examples. Since
we are sticking to unital algebras, we will only consider the discrete
case.

Following standard notation, we will let $C^*(\Gamma)$ and
$C^*_r(\Gamma)$ denote the full (i.e.\ universal) and reduced
C$^*$-algebra, respectively, of a discrete group $\Gamma$.  When
studying amenable traces it turns out that there is an enormous
difference between considering the full and reduced algebras.  For
$C^*_r(\Gamma)$ there is a sort of `all or nothing' principle
which asserts that either every trace on $C^*_r(\Gamma)$ is
amenable or no trace on $C^*_r(\Gamma)$ is amenable -- the first
case occurring if and only if $\Gamma$ is an amenable group.

\begin{prop}{\em (Reduced group C$^*$-algebras -- compare with
\cite[Corollary 2.11]{bedos:hypertraces}.)}
\label{thm:reducedgroupalgebras} For a discrete group $\Gamma$ the
following are equivalent:
\begin{enumerate}
\item $\Gamma$ is amenable.

\item $\mathrm{T}(C^*_r(\Gamma)) = \mathrm{AT}(C^*_r(\Gamma))$.

\item $\mathrm{AT}(C^*_r(\Gamma)) \neq \emptyset$.
\end{enumerate}
\end{prop}

\begin{proof}
(1) $\Longrightarrow$ (2).
A result of Lance \cite{lance:nuclear} states that $C^*_r(\Gamma)$ is
nuclear whenever $\Gamma$ is amenable.  Hence for every C$^*$-algebra
$B$ there is a unique C$^*$-norm on the algebraic tensor product
$C^*_r(\Gamma) \odot B$.  In particular, $C^*_r(\Gamma) \otimes B$
enjoys the universal property of the maximal tensor
product\footnote{Every pair of representations $\pi:C^*_r(\Gamma) \to
B(H)$ and $\sigma:B \to B(H)$ with commuting ranges induces a natural
$*$-homomorphism $C^*_r(\Gamma) \otimes B \to B(H)$.} and hence it is
clear from condition (4) of Theorem \ref{thm:amenabletraces} that
every trace on $C^*_r(\Gamma)$ is amenable.

(2) $\Longrightarrow$ (3) is immediate since discreteness implies
that $C^*_r(\Gamma)$ always has a tracial state.

(3) $\Longrightarrow$ (1).
Suppose that $\tau$ is an amenable trace on $C^*_r(\Gamma)$.  Then we
may extend $\tau$ to a state $\phi$ on $B(l^2(\Gamma))$ such that
$$\phi(\lambda_g T \lambda_{g}^*) = \phi(T)$$ for all $g \in \Gamma$
and $T \in B(l^2(\Gamma))$.  A well known calculation shows that if we
regard $l^{\infty}(\Gamma) \subset B(l^2(\Gamma))$ as multiplication
operators then the left translation action of $\Gamma$ on
$l^{\infty}(\Gamma)$ is just given by conjugation: $f \mapsto
\lambda_g f \lambda_{g}^*$.  Hence the restriction of the state $\phi$
to $l^{\infty}(\Gamma)$ defines a classical invariant mean which
implies $\Gamma$ is amenable.
\end{proof}

\begin{cor} $\Gamma$ is not amenable if and only if
$C^*_r(\Gamma)$ has no amenable traces.
\end{cor}

Jonathan Rosenberg first observed that if $C^*_r(\Gamma)$ is
quasidiagonal then $\Gamma$ must be amenable.  He has conjectured that
the converse holds as well.  While this conjecture remains elusive we
now point out that (a) Rosenberg's observation is easily deduced from
amenable trace considerations and (b) his conjecture would follow from a
similar statement about traces.

First suppose that $C^*_r(\Gamma)$ is quasidiagonal.  Then
$C^*_r(\Gamma)$ would have to have at least one amenable trace (cf.\
Proposition \ref{thm:QDcase}) and hence, by the last proposition,
$\Gamma$ would be amenable.  Note that this argument requires much
less than quasidiagonality since this was only used to get the
existence of an amenable trace.

Now to the tracial statement that would imply Rosenberg's conjecture.

\begin{prop}{\em (Rosenberg's Conjecture)}
\label{thm:rosenberg} If the canonical trace on $C^*_r(\Gamma)$ is
quasidiagonal then $C^*_r(\Gamma)$ is quasidiagonal.  In
particular, if this is the case for every amenable group then
Rosenberg's conjecture would follow.
\end{prop}

\begin{proof}
This is actually a special case of a more general result: If $\tau
\in \TAQD$ is faithful then $A$ is quasidiagonal.  Indeed, assume
that $\phi_n:A \to M_{k(n)}({\mathbb C})$ are u.c.p.\ maps which
are asymptotically multiplicative (in norm) and recover $\tau$.
Since these maps are almost multiplicative in norm we get a
$\tau$-preserving $*$-homomorphism into the quotient of the norm
bounded sequences by the ideal of sequences which tend to zero in
norm $$\frac{\Pi M_{k(n)}({\mathbb C})}{\oplus M_{k(n)}({\mathbb
C})}$$ similar to the proof of (1) $\Longrightarrow$ (2) from
Theorem \ref{thm:amenabletracesII}.  Since $\tau$ is faithful so
is the induced $*$-homomorphism $$A \to \frac{\Pi
M_{k(n)}({\mathbb C})}{\oplus M_{k(n)}({\mathbb C})}$$ and hence
$\|a\| = \limsup \|\phi_n(a)\|$ for every $a \in A$.  From
Voiculescu's abstract characterization of quasidiagonality it
follows that $A$ must be quasidiagonal.\footnote{The careful
reader may be worried about only having  $\|a\| = \limsup
\|\phi_n(a)\|$ instead of an honest limit. However it is a routine
exercise to replace the $\phi_n$'s with asymptotically
multiplicative maps satisfying the right limit condition by
passing to subsequences and taking direct sums.}
\end{proof}

Now we consider the full group C$^*$-algebras $C^*(\Gamma)$ where
things are totally different.  For example, this C$^*$-algebra
always has at least one amenable trace.  In fact, the trivial
representation can be regarded as a (uniform locally) finite
dimensional trace on $C^*(\Gamma)$.

There is another trace on $C^*(\Gamma)$ coming from the left
regular representation and, in contrast to the reduced
C$^*$-algebra case, this trace is also amenable for a large class
of non-amenable discrete groups. For the remainder of this section
we will let {\em $\tau$ denote the trace on $C^*(\Gamma)$ coming
from the left regular representation}.

\begin{prop}{\em (Residually Finite Groups)}
\label{thm:residuallyfinitegroups} If $\Gamma$ is residually
finite then $\tau$ is a locally finite dimensional trace on
$C^*(\Gamma)$. In particular, it is amenable.
\end{prop}

\begin{proof} We will show that there are finite dimensional
$*$-representations $\pi_n:C^*(\Gamma) \to M_{k(n)}({\mathbb C})$
such that for each non-trivial group element $g \in \Gamma$ we
have $$\mathrm{tr}_{k(n)}\circ\pi_n(\lambda_g) \to 0.$$  Since the
linear span of such unitaries, together with the unit, is dense in
$C^*(\Gamma)$ it will follow that $\tau$ is the weak-$*$ limit of
finite dimensional traces and the proof will be complete.

To construct such representations we let $\Gamma_1
\trianglerighteq \Gamma_2 \trianglerighteq \ldots$ be a descending
sequence of normal subgroups each of which has finite index in
$\Gamma$ and such that their intersection is the neutral element.
Let $\pi_n : C^*(\Gamma) \to B(l^2(\Gamma/\Gamma_n))$ be the
unitary representation induced by the left regular representation
of $\Gamma/\Gamma_n$.  Since $\Gamma/\Gamma_n$ is a finite group,
$B(l^2(\Gamma/\Gamma_n))$ is finite dimensional and a standard
calculation completes the proof.
\end{proof}

\begin{cor}{\em (Maximally Almost Periodic Groups)}
\label{thm:maximallyalmostperiodic} If $\Gamma$ is maximally
almost periodic\footnote{In other words, $\Gamma$ is isomorphic to
a subgroup of a compact group.  This class of groups contains all
residually finite groups.} then $\tau$ is a locally finite
dimensional trace on $C^*(\Gamma)$.
\end{cor}

\begin{proof} It is well known that every discrete maximally
almost periodic group is the increasing union of its residually
finite subgroups.\footnote{Since all finitely generated linear
groups are residually finite (cf.\ \cite{alperin}) and maximally
almost periodic groups have a separating family of homomorphisms
into linear groups it can be shown that every finitely generated
maximally almost periodic group is  residually finite.  Since
every discrete group is the union of its finitely generated
subgroups the claim follows.} If $\Lambda \subset \Gamma$ is a
subgroup then $C^*(\Lambda) \subset C^*(\Gamma)$ and the
restriction of the canonical trace on $C^*(\Gamma)$ is just the
canonical trace on $C^*(\Lambda)$. With these observations the
remainder of the proof amounts to a standard approximation
argument.  Indeed, if $\mathfrak{F} \subset C^*(\Gamma)$ is a
finite set then, after a small perturbation, we may assume that
$\mathfrak{F} \subset C^*(\Lambda) \subset C^*(\Gamma)$ where
$\Lambda$ is a residually finite subgroup of $\Gamma$.  By the
previous corollary we can find maps from $C^*(\Lambda)$ to
matrices whose multiplicative domains almost contain
$\mathfrak{F}$ and recover the canonical trace. We then extend
these maps to all of $C^*(\Gamma)$ by Arveson's extension theorem.
\end{proof}

We are now in a position to illustrate the difference between
being an amenable trace and a uniform amenable trace.

\begin{prop}
\label{thm:what?} If $\tau$ is a uniform amenable trace on
$C^*(\Gamma)$ then $\Gamma$ is amenable.
\end{prop}

\begin{proof}
If $\tau$ is a uniform amenable trace then statement (2) of
Proposition \ref{thm:passagetoquotient} implies that the {\em
reduced} group C$^*$-algebra also has a (uniform) amenable trace
and hence $\Gamma$ is amenable.
\end{proof}

\begin{rem} Propositions \ref{thm:residuallyfinitegroups} and
\ref{thm:what?} imply that the sets $\TAuim$, $\TAuQD$ and
$\TAulfd$ need not be weak-$*$ closed in general.  Indeed, if
$\Gamma$ is residually finite and non-amenable then $\tau$ arises
as the weak-$*$ limit of (uniform locally) finite dimensional
traces however it is {\em not} uniform amenable.
\end{rem}

Proposition \ref{thm:residuallyfinitegroups} and Corollary
\ref{thm:maximallyalmostperiodic} are well known (cf.\
\cite{kirchberg:propertyTgroups}, \cite{wassermann:exactbook}) but
usually formulated in terms of {\em Kirchberg's factorization
property}.

\begin{defn} A discrete group $\Gamma$ has {\em Kirchberg's
factorization property} if the natural $*$-representation
$$C^*(\Gamma) \odot C^*(\Gamma) \to B(l^2(\Gamma))$$ induced by
the left and right regular representations is continuous with
respect to the minimal tensor product norm.\footnote{In other
words, the $*$-representation defined on the maximal tensor
product factors through the minimal tensor product.}
\end{defn}

The following theorem is an immediate consequence of Theorem
\ref{thm:amenabletraces}.  Indeed, the `right regular'
representation used in that theorem is nothing but the obvious
generalization of the right regular representation of a group.

\begin{thm}{\em (Kirchberg)} $\Gamma$ has the factorization
property if and only if $\tau$ is an amenable trace on
$C^*(\Gamma)$.
\end{thm}

Evidently all amenable groups have the factorization property
(since $C^*(\Gamma) = C^*_r(\Gamma)$ is nuclear) but there are
plenty of non-amenable examples as well.  For example, all
residually finite groups (e.g.\ free groups, $SL(n,\mathbb{Z})$,
etc.) and, more generally, maximally almost periodic groups. With
this observation in hand we can now recover a result of Bekka
which states that Rosenberg's conjecture has an affirmative answer
in the case that $\Gamma$ is amenable and maximally almost
periodic. (See  \cite[Corollary 4]{dadarlat:QDapproximation} for
another proof of this result.)

\begin{cor}{\em (Bekka)} If $\Gamma$ is
amenable and maximally almost periodic then $C^*(\Gamma) =
C^*_r(\Gamma)$ is quasidiagonal.
\end{cor}

\begin{proof} Since $\Gamma$ is amenable the full and reduced
C$^*$-algebras are equal.  Hence Proposition
\ref{thm:maximallyalmostperiodic} implies the canonical trace on
$C^*_r(\Gamma)$ is quasidiagonal and so we can apply Proposition
\ref{thm:rosenberg}.
\end{proof}

The main result of \cite{kirchberg:propertyTgroups} asserts that
every discrete group which has both the factorization property and
Kazhdan's property T must be residually finite.  The main
technical tool in the proof is the following fact (see
\cite[Proposition 2.3]{kirchberg:propertyTgroups} for a
proof\footnote{The reader acquainted with Property T arguments may
find this a fun exercise. Basically one applies Stinespring to a
u.c.p.\ map $\phi:A \to M_n(\mathbb{C})$.  If $\phi$ is almost
multiplicative in 2-norm then the Stinespring projection gives a
Hilbert-Schmidt operator which is almost invariant under the
conjugation action.  One perturbs, by Property T, to a
Hilbert-Schmidt operator in the commutant of the Stinespring
representation and some functional calculus and other standard
estimates complete the proof.}).

\begin{lem}{\em (Kirchberg)} If $\Gamma$ has property T and
$\phi_n:C^*(\Gamma) \to M_{k(n)}(\mathbb{C})$ are u.c.p.\ maps
such that $\| \phi_n(ab) - \phi_n(a)\phi_n(b) \|_2 \to 0$ for all
$a,b \in C^*(\Gamma)$ then there exist $*$-homomorphisms
$\pi_n:C^*(\Gamma) \to M_{l(n)}(\mathbb{C})$ such that $$|
\mathrm{tr}_{k(n)}(\phi_n(a)) - \mathrm{tr}_{l(n)}(\pi_n(a))| \to
0$$ for all $a \in C^*(\Gamma)$.
\end{lem}

Given this lemma it is a simple matter to provide the proof of the
following fact.

\begin{prop} If $\Gamma$ is a property T group then
$$\mathrm{AT}(C^*(\Gamma)) = \mathrm{AT}(C^*(\Gamma))_{\mathrm{QD}} =
\mathrm{AT}(C^*(\Gamma))_{\mathrm{LFD}}.$$
\end{prop}

Since the factorization property is equivalent to knowing that
$\tau$ is an amenable trace on $C^*(\Gamma)$ we can combine
Kirchberg's result quoted above with some work of Olshanskii
\cite{olshanskii} to get the following connection between amenable
traces and one of the major open problems in geometric group
theory.

\begin{prop}
\label{thm:residuallyfinite} Every hyperbolic group is residually
finite if and only if $\tau$ is an amenable trace on $C^*(\Gamma)$
for every hyperbolic group $\Gamma$.
\end{prop}

\begin{proof} Since the canonical trace is always amenable
for residually finite groups we only need to observe the `if'
statement.  However in \cite{olshanskii} Olshanskii shows, among other
things, that every hyperbolic group can be embed into a hyperbolic
group with property T.  Since, for property T groups, residual
finiteness and the factorization property are equivalent the
proposition follows.
\end{proof}

Finally we consider the case of free groups.  As observed in
\cite{kirchberg:invent} understanding the space of amenable traces
on $C^*({\mathbb F}_{n})$ is of paramount importance.\footnote{The
main question, which turns out to be equivalent to Connes'
embedding problem, is whether or not every trace is amenable on
$C^*({\mathbb F}_{\infty})$.}  For now we just observe that, like
the case of property T groups, some of the sets of traces defined
here coincide in this case.

\begin{prop}
\label{thm:freegroupexample} For a free group we have
$$\mathrm{AT}(C^*({\mathbb F}_{n})) =
\mathrm{AT}(C^*({\mathbb F}_{n}))_{\mathrm{QD}} =
\mathrm{AT}(C^*({\mathbb F}_{n}))_{\mathrm{LFD}}.$$
\end{prop}

\begin{proof} This follows from a simple observation of Kirchberg
(see \cite[Lemma 4.5]{kirchberg:invent} or \cite{hadwin}). Namely,
since $R^{\omega}$ is a von Neumann algebra every unitary can be
lifted to a unitary in $l^{\infty}(R)$.  A standard approximation
argument yields the following improvement: If $\pi:C^*({\mathbb
F}_{n}) \to R^{\omega}$ is a $*$-homomorphism then there exist
$*$-homomorphisms $\pi_n:C^*({\mathbb F}_{n}) \to
M_{k(n)}(\mathbb{C})$ such that $$|\tau_{\omega}\circ\pi(a) -
\mathrm{tr}_{k(n)}\circ\pi_n(a) | \to 0$$ for all $a \in
C^*({\mathbb F}_{n})$.\footnote{First lift the generating
unitaries to unitaries in $l^{\infty}(R)$ and then approximate
them (in 2-norm) by unitaries in large finite dimensional
subfactors of $R$.} This evidently implies that every amenable
trace arises as the weak-$*$ limit of finite dimensional traces
and the proof is complete.
\end{proof}

\section{Nuclear and WEP C$^*$-algebras}

We now see how the theory of amenable traces develops for several
important classes of examples.  We begin with two where it turns out 
that every trace enjoys some sort of
amenability -- nuclear C$^*$-algebras and C$^*$-algebras with the
weak expectation property (WEP) of Lance (cf.\
\cite{lance:nuclear}).

In the proof of Proposition \ref{thm:reducedgroupalgebras} we saw
that every trace on a nuclear C$^*$-algebra is amenable.  This
result is a trivial consequence of the original tensor product
definition of nuclearity\footnote{$A$ is nuclear if and only if
there is a unique C$^*$-norm on $A \odot B$ for every
C$^*$-algebra $B$.} and Lance's tensor product trick from the
proof of Theorem \ref{thm:amenabletraces}. A far more substantial
result, however, is the fact that every trace on a nuclear
C$^*$-algebra is {\em uniform} amenable.

\begin{thm}{\em (Nuclear C$^*$-algebras)}
\label{thm:nuclear} If $A$ is nuclear then $$\TA = \TAim = \TAuim
.$$
\end{thm}

\begin{proof} This is an immediate consequence of Theorem
\ref{thm:mainthmUTAwafd} and the celebrated fact that an injective
von Neumann algebra is hyperfinite (cf.\
\cite{connes:classification}, \cite{haagerup:newproof},
\cite{popa:injectiveimplieshyperfiniteI}) since it is not hard to
show that every representation of a nuclear C$^*$-algebra yields
an injective von Neumann algebra (use Lance's trick from the proof
of Theorem \ref{thm:amenabletraces} to show the commutant is
injective).
\end{proof}

Recall that a C$^*$-algebra $A$ is said to have the {\em weak
expectation property} (WEP) if for every {\em faithful}
representation $\pi : A \to B(H)$ there exists a u.c.p.\ map $\Phi
: B(H) \to \pi(A)^{\prime\prime}$ such that $\Phi(\pi(a)) =
\pi(a)$, for all $a \in A$.  This is a large class of
C$^*$-algebras (including every injective C$^*$-algebra and every
nuclear C$^*$-algebra).  In fact, \cite[Corollary
3.5]{kirchberg:invent} states that every (separable) C$^*$-algebra
is contained in a (separable) simple C$^*$-algebra with the WEP
(in particular, the WEP class contains lots of separable
non-nuclear C$^*$-algebras).

\begin{prop}{\em (C$^*$-algebras with the WEP)}
\label{thm:WEPcase}
If $A$ has the WEP then $\TA = \TAim$.
\end{prop}

\begin{proof} Let $\tau \in \TA$ be given and assume that $A
\subset A^{**} \subset B(H)$ is the universal
representation.\footnote{This violates our standing separability
assumption, but otherwise causes no harm.}  Applying the
definition of the weak expectation property to this representation
we get a u.c.p.\ map $\Phi : B(H) \to A^{**}$ such that $\Phi(a) =
a$ for all $a \in A$.  In particular, $A$ falls in the
multiplicative domain of $\Phi$ and hence, letting $\tau^{**}$
denote the normal extension of $\tau$ to $A^{**}$, we have that
$\tau^{**} \circ \Phi$ is a state on $B(H)$ which extends $\tau$
and contains $A$ in its centralizer.
\end{proof}

Though the previous proposition is a simple consequence of the
definitions it turns out to have some interesting consequences.
Indeed, a common phenomenon in mathematics is that certain
problems can be easily answered once the right point of view is
found.  Amenable traces turn out to be the ``right" concept to
consider for a number of operator algebraic questions as we will
see in this paper.

Our first application is related to Kirchberg's remarkable
embedding theorem: A C$^*$-algebra is exact if and only if it
embeds into the Cuntz algebra ${\mathcal O}_2$. It is natural to
wonder if $A$ is exact and has additional properties whether or
not $A$ can be embed into a nuclear C$^*$-algebra with the same
additional properties.  For example, it is not known whether or
not every exact, quasidiagonal C$^*$-algebra can be embed into a
nuclear, quasidiagonal C$^*$-algebra. (See, however,
\cite{ozawa:AFembed} where cones and suspensions over exact
C$^*$-algebras are shown to be embeddable into AF algebras!) On
the other hand, some experts have already observed that the stably
finite, exact C$^*$-algebra $C^*_r ({\mathbb F}_2)$ can't be embed
into a stably finite nuclear C$^*$-algebra.\footnote{Whether or
not every exact, stably finite C$^*$-algebra can be embed into a
nuclear, stably finite C$^*$-algebra was a problem pondered by
some experts a few years ago.  This problem naturally arose from
an important question of Blackadar-Kirchberg which asks whether or
not every nuclear, stably finite C$^*$-algebra is quasidiagonal. A
natural way to provide counterexamples to the quasidiagonal
question is to try to embed a (necessarily exact and stably
finite) non-quasidiagonal C$^*$-algebra into a stably finite,
nuclear C$^*$-algebra. Unfortunately this strategy is doomed to
failure for all of the examples of exact, stably finite,
non-quasidiagonal C$^*$-algebras that we currently know -- see
Proposition \ref{thm:BFDIII}.}  The proof typically quoted depends
on the uniqueness of the trace on $C^*_r ({\mathbb F}_2)$,
Haagerup's remarkable result that every unital, stably finite,
exact C$^*$-algebra has a tracial state
(\cite{haagerup:quasitraces}, \cite{haagerup-thorbjornsen})
together with Connes' theorem that injective implies hyperfinite
and hence is highly non-trivial. On the other hand, amenable
traces show that this is a general fact about non-amenable
discrete groups. (See also Bed{\o}s' observations along these
lines \cite[Corollary 2.11]{bedos:hypertraces}.)

\begin{cor}
\label{thm:nohomomorphisms}
Let $\Gamma$ be a discrete group. One can find a 
$*$-homomorphism $\pi : C^*_r (\Gamma) \to A$, where $A$ has the WEP
and at least one tracial state, if and only if $\Gamma$ is amenable.
\end{cor}

\begin{proof} If $\Gamma$ is amenable then $C^*_r (\Gamma)$ is nuclear
(hence WEP) and has a trace.  On the other hand, since every trace
on $A$ is amenable, the existence of $\pi$ would imply that $C^*_r
(\Gamma)$ has an amenable trace (just restrict the trace on $A$ to
$\pi(C^*_r (\Gamma))$) and hence, by Proposition
\ref{thm:reducedgroupalgebras}, $\Gamma$ would be amenable.
\end{proof}

\begin{cor}
\label{thm:nohomomorphisms2}
If $\Gamma$ is a discrete, non-amenable group then $C^*_r (\Gamma)$
can't be embed into a stably finite, nuclear C$^*$-algebra or into a
finite, hyperfinite von Neumann algebra.
\end{cor}

\begin{proof}
Hyperfinite, finite von Neumann algebras are injective (hence WEP)
and have a tracial state.  Every unital, stably finite, nuclear
C$^*$-algebra has a tracial state (cf.\ \cite{haagerup-thorbjornsen}).
\end{proof}

Note, of course, that not only are embeddings impossible but one can't
even find a non-zero $*$-homomorphism from $C^*_r (\Gamma)$ to a
stably finite, nuclear C$^*$-algebra when $\Gamma$ is non-amenable.
This result is in contrast to the following theorem of Kirchberg:
Every exact C$^*$-algebra admits a complete order embedding into the
CAR algebra (cf.\ \cite[Theorem 9.1]{wassermann:exactbook}).  All
together these results say that if $\Gamma$ is exact and non-amenable
then $*$-homomorphic embeddings of $C^*_r (\Gamma)$ into purely
infinite, nuclear algebras always exist; complete order embeddings
into stably finite, nuclear algebras always exist; but, non-zero
$*$-homomorphisms to stably finite, nuclear algebras {\em never}
exist.

We now use amenable traces to give a simple proof of a
generalization of \cite[Theorem 6.8]{bekka}.  If $\pi : G \to
B(H)$ is a unitary representation of a locally compact group $G$
then Bekka defines the {\em representation} to be amenable if
there exists a state $\varphi$ on $B(H)$ such that $\varphi(\pi(g)
T \pi(g^{-1})) = \varphi(T)$ for all $T \in B(H)$ and for all $g
\in G$.  In the language of this paper $\pi$ is amenable if and
only if $A = C^*( \{ \pi(g) : g \in G \} ) = span \{ \pi(g) : g
\in G \}^{-}$ (norm closure, since $G$ is a group) has an amenable
trace (it is a simple exercise to show that this is equivalent to
Bekka's definition). In \cite[Theorem 6.8]{bekka} Bekka proves the
following result in the case that $A$ is nuclear.

\begin{cor}
Let $\pi : G \to B(H)$ be a strongly continuous, unitary
representation of a locally compact group $G$ and $A = C^*( \{ \pi(g)
: g \in G \} )$.  If $A$ has the WEP then $\pi$ is an amenable
representation if and only if $A$ has a tracial state.
\end{cor}

\begin{proof} Strictly speaking one might be concerned about
separability issues here, but we only assume separability for
convenience and this is not necessary for virtually everything
discussed in this paper.  If $\pi$ is amenable then, by
definition, $A$ has an amenable trace.  If $A$ has the WEP and a
tracial state then this trace is amenable and hence $\pi$ is
amenable.
\end{proof}

\begin{rem}  In the paragraph preceding \cite[Theorem
6.8]{bekka}, Bekka asked if there is any general relationship
between the amenability of $\pi$ and the nuclearity of $A = C^*(
\{ \pi(g) : g \in G \} )$. The answer is that -- in general --
there is absolutely no relationship.  That is, it is possible that
$A$ be nuclear while $\pi$ is not amenable (let $A$ be a Cuntz
algebra, for example, and $\Gamma$ any dense subgroup of the
unitary group of $A$) and that $A$ be non-nuclear while $\pi$ is
amenable (take $A$ to be any non-nuclear Popa algebra and then $A$
must have at least one quasidiagonal trace since it is a
quasidiagonal C$^*$-algebra).
\end{rem}

\section{Locally reflexive, exact and quasidiagonal C$^*$-algebras}

We regard amenable traces as corresponding to nice GNS
representations (`nice' is quite ambiguous at the moment).
One natural question to ask is whether or not being `nice' can be
characterized in terms of the von Neumann algebra one gets in the
GNS representation. Unfortunately this is not the case, in
general, since the canonical traces on $C^*({\mathbb F}_{2})$ and
$C^*_r({\mathbb F}_{2})$ yield the same free group factor while
the trace on $C^*({\mathbb F}_{2})$ is amenable (by residual
finiteness) and the trace on $C^*_r({\mathbb F}_{2})$ is not
amenable (by non-amenability of free groups).  Thus it is the {\em
representation} (including the domain), as opposed to the {\em
image} of the representation, which matters in general. However,
if one assumes some additional structure on the algebra then the
{\em image} of the representation does in fact characterize
`niceness'. We will see another example of this later on (cf.\
Proposition \ref{thm:LLPexample}): If $A$ has the so-called local
lifting property then a trace on $A$ is amenable if and only if
the corresponding GNS von Neumann algebra admits a trace
preserving embedding into the ultraproduct of the hyperfinite
II$_1$-factor. In this section we show that assuming local
reflexivity gives a much sharper statement (see Corollary
\ref{cor:locallyreflexive}).

The notion of local reflexivity in the context of C$^*$-algebras
is due to Effros and Haagerup \cite{effros-haagerup}.  A
remarkable result of Kirchberg asserts that every exact
C$^*$-algebra has this property
\cite{kirchberg:commutantsofunitaries} and he has conjectured that
the converse holds as well.

\begin{defn} $A$ is {\em locally reflexive} if for each finite
dimensional operator system $X \subset A^{**}$ there exists a net
of u.c.p.\ maps $\phi_{\lambda}:X \to A$ such that
$$\eta(\phi_{\lambda}(x)) \to \eta(x)$$ for all $x \in X$ and
$\eta \in A^*$. (i.e.\ $\phi_{\lambda} \to id_{X}$ in the point
weak-$*$ topology.)
\end{defn}

The following improvement is a standard Hahn-Banach type argument
but we will soon find it quite useful.

\begin{lem} Assume $A$ is locally reflexive and $\mathfrak{F}
\subset A$ is a finite set.  Then for any finite dimensional
operator system $\mathfrak{F} \subset X \subset A^{**}$ there
exists a net of u.c.p.\ maps $\phi_{\lambda}:X \to A$ such that
$$\eta(\phi_{\lambda}(x)) \to \eta(x)$$ for all $x \in X$ and
$\eta \in A^*$ and such that $\|a - \phi_{\lambda}(a) \| \to 0$
for each $a \in \mathfrak{F}$.
\end{lem}

\begin{proof} Assume first that $\mathfrak{F}$ consists of a single
element $a \in A$. Applying the Hahn-Banach theorem to the convex
hull of the u.c.p.\ maps provided by the definition of locally
reflexivity one can find a net converging in norm.  For general
finite sets one employs the usual trick of taking direct sums and
reducing to the singleton case (as we already saw in the proof of
(1) $\Longrightarrow$ (2) from Theorem \ref{thm:amenabletraces},
for example).
\end{proof}

The main result concerning amenable traces on locally reflexive
C$^*$-algebras is as follows.

\begin{thm}
\label{thm:locallyreflexive} Assume that  $A$ is locally
reflexive. Then it is always the case that $\TAim = \TAuim$ and
$\TAQD = \TAuQD$.\footnote{As previously mentioned, we don't know
whether $\TAlfd = \TAulfd$ for all locally reflexive $A$ -- this
is a very important open problem (see Section
\ref{thm:classificationprogram}).}
\end{thm}

\begin{proof}
We only give the proof of $\TAQD = \TAuQD$ as it will be clear that
essentially the same proof gives the other equality.

So let $\tau \in \TAQD$ be arbitrary.  Evidently it suffices to
prove that if $\mathfrak{F} \subset A$ is an arbitrary finite set
and $\varepsilon > 0$ then there exists a u.c.p.\ map $\varphi : A
\to B$, where $B$ is a finite dimensional C$^*$-algebra, such that
$\| \varphi(xy) - \varphi(x)\varphi(y) \| < \varepsilon$, for all
$x,y \in \mathfrak{F}$ and such that there exists a trace $\gamma$
on $B$ such that $\| \tau - \gamma\circ\varphi \|_{A^*} <
\varepsilon$ (cf.\ Lemma \ref{thm:tracepreservingembedding}).

In order to do this, we will show that for each finite dimensional
operator system $X \subset A^{**}$ containing both the set
$\mathfrak{F}$ and $\{ ab : a,b \in \mathfrak{F} \}$, there exists
a sequence of normal, u.c.p.\ maps $\psi_n : A^{**} \to M_{s(n)}
({\mathbb C})$ such that $\mathrm{tr}_{s(n)} \circ \psi_n (x) \to
\tau^{**} (x)$, for all $x \in X$ and {\em each} $\psi_n$ is
$\varepsilon$-multiplicative on $\mathfrak{F}$.  If we are able to
do this then one can construct a net of normal, u.c.p.\ maps
$\varphi_{\lambda} : A^{**} \to M_{k(\lambda)} ({\mathbb C})$ with
the property that $\mathrm{tr}_{k(\lambda)} \circ
\varphi_{\lambda} \in (A^{**})_* = A^*$ for all $\lambda$, each
$\varphi_{\lambda}$ is $\varepsilon$-multiplicative on
$\mathfrak{F}$ and (here is the key) $\mathrm{tr}_{k(\lambda)}
\circ \varphi_{\lambda} \to \tau^{**}$ in the weak topology coming
from $A^{**}$.  Hence, by the Hahn-Banach theorem, $\tau^{**}$
belongs to the {\em norm} closure of the convex hull of $\{
\mathrm{tr}_{k(\lambda)} \circ \varphi_{\lambda} \} \subset A^*$.
Then one would be able to choose a finite set $\lambda_1, \ldots,
\lambda_p$ and positive real numbers $\theta_1, \ldots, \theta_p$
such that $\sum \theta_i = 1$ and $\| \tau^{**} - \sum \theta_i
\mathrm{tr}_{k(\lambda_i)} \circ \varphi_{\lambda_i} \|_{A^*} <
\varepsilon$.  Finally one would define $B =
M_{k(\lambda_1)}({\mathbb C}) \oplus \cdots \oplus
M_{k(\lambda_p)} ({\mathbb C})$, $\varphi = \varphi_{\lambda_1}
\oplus \cdots \oplus \varphi_{\lambda_p}$ and $\gamma = \sum
\theta_i \mathrm{tr}_{k(\lambda_i)}$.  Since we arranged that
$\varphi_{\lambda}$ is $\varepsilon$-multiplicative on
$\mathfrak{F}$ {\em for every} $\lambda$, it is clear that
$\varphi$ will also be close to multiplicative on $\mathfrak{F}$.

So let $X \subset A^{**}$ be any finite dimensional operator
system containing the sets $\mathfrak{F}$ and $\{ ab : a,b \in
\mathfrak{F} \}$.  Since $\tau \in \TAQD$ we can find a sequence
of u.c.p.\ maps $\varphi_m : A \to M_{k(m)} ({\mathbb C})$ which
are asymptotically multiplicative (in norm) and which recapture
$\tau$ (as a weak-$*$ limit) after composing with the traces on
$M_{k(m)} ({\mathbb C})$. Note that by passing to a subsequence,
if necessary, we may further assume that each $\varphi_m$ is as
close to multiplicative as one likes on the set $\mathfrak{F}$.
Since $A$ is locally reflexive, we can find a net of u.c.p.\ maps
$\beta_t : X \to A$ such that $\beta_t (x) \to x$ in the weak-$*$
topology (coming from $A^*$) for all $x \in X$. By the lemma
preceding this theorem we may further assume that  $\| a - \beta_t
(a) \| \to 0$ for all $a \in \mathfrak{F} \cup \{ ab : a,b \in
\mathfrak{F} \}$. Passing to a subnet, if necessary, we may also
assume that $\| \beta_t (a) - a \| < \varepsilon$ for all $a \in
\mathfrak{F} \cup \{ ab : a,b \in \mathfrak{F} \}$ and for all
$t$. In particular, this implies that $\varphi_m \circ \beta_t$ is
nearly multiplicative on $\mathfrak{F}$ for all $t,m$.

We almost have the desired maps $\psi_n$. Since $X$ is finite
dimensional we can choose a linear basis $\{ x_1, \ldots, x_q \}$.
For each $n \in {\mathbb N}$ first choose $t(n)$ such that $|
\tau^{**} (x_i) - \tau(\beta_{t(n)}(x_i)) | < 1/n$ for $1 \leq i
\leq q$.  Next, choose $m_n$ such that $$| \tau(\beta_{t(n)}(x_i)) -
\mathrm{tr}_{k(m_n)} \circ \varphi_{m_n} (\beta_{t(n)}(x_i)) | <
1/n$$ for $1 \leq i \leq q$.  Then defining $\tilde{\psi}_n =
\varphi_{m_n} \circ \beta_{t(n)} : X \to M_{k(m_n)} ({\mathbb C})$
we have that $\mathrm{tr}_{k(m_n)} \circ \tilde{\psi}_n (x) \to
\tau^{**}(x)$ for all $x \in X$.

By Arveson's extension theorem, we may assume that each
$\tilde{\psi}_n$ is actually defined on all of $A^{**}$.  The only
problem is that we can't be sure that the Arveson extensions are
normal on $A^{**}$.  However a tiny perturbation of the
$\tilde{\psi}_n$ will yield normal maps $\psi_n$ with all the
right properties. (Use the fact that c.p.\ maps to matrix algebras
are nothing but positive linear functionals on matrices over the
given algebra and that the set of normal linear functionals on a
von Neumann algebra is dense in the dual space.)
\end{proof}

Theorem \ref{thm:mainthmUTAwafd} together with the result above gives
our first corollary.

\begin{cor}
\label{cor:locallyreflexive} Let $A$ be locally reflexive and
$\tau$ be a tracial state on $A$. Then $\tau$ is amenable if and
only if $\pi_{\tau}(A)^{\prime\prime}$ is hyperfinite.
\end{cor}

The next corollary was proved in \cite[Theorem
7.5]{kirchberg:invent}. Our proof is not any simpler than
Kirchberg's but it emphasizes the role of amenable traces and
representation theory.  Kirchberg, in fact, also used (classical)
invariant means but his emphasis is on exactness and tensor
products (cf.\ \cite[Proposition 7.1]{kirchberg:invent}).

\begin{cor}
Let $\Gamma$ be a discrete group with the factorization property.
Then $\Gamma$ is amenable if and only if $C^*(\Gamma)$ is locally
reflexive. In particular, if $\Gamma$ is any residually finite,
non-amenable group then $C^*(\Gamma)$ is not locally reflexive (hence
not exact).
\end{cor}

\begin{proof} If $\Gamma$ is amenable then $C^*(\Gamma)$ is nuclear
and hence locally reflexive.  If $C^*(\Gamma)$ is locally
reflexive and $\tau$ is an amenable trace then $\pi_{\tau}
(C^*(\Gamma))^{\prime\prime}$ is a hyperfinite von Neumann
algebra. Hence the canonical trace on the {\em reduced} algebra
$C^*_r(\Gamma)$ is also amenable and thus $\Gamma$ is amenable.
\end{proof}

We won't actually need condition (3) in the next corollary however, it
may be of independent interest.  Pisier gave the first proof of (3)
$\Longrightarrow$ (4), though Ozawa later observed that exactness was
an unnecessary assumption (see Theorem \ref{thm:mainthmUTAwafd}).

\begin{cor}
\label{prop:mainpropsection4}
Let $A$ be an exact C$^*$-algebra and $\tau \in \TA$. Then the
following are equivalent:

\begin{enumerate}
\item $\tau \in \TAim$.

\item $\tau \in \TAuim$.

\item The GNS representation is a nuclear map into $\pi_{\tau}
(A)^{\prime\prime}$.  That is, there exist u.c.p.\ maps $\phi_n : A \to
M_{k(n)} ({\mathbb C})$ and $\psi_n : M_{k(n)} ({\mathbb C}) \to
\pi_{\tau} (A)^{\prime\prime}$ such that $\| \pi_{\tau} (a) - \psi_n
\circ \phi_n (a) \| \to 0$ for all $a \in A$.

\item $\pi_{\tau} (A)^{\prime\prime}$ is hyperfinite.
\end{enumerate}
\end{cor}

\begin{proof} Since exact C$^*$-algebras are locally reflexive
(cf.\ \cite[pg.\ 71]{kirchberg:commutantsofunitaries}) we have the
equivalence of (1), (2) and (4).  Since nuclearity obviously implies
weak nuclearity we also get the implication (3) $\Longrightarrow$ (4)
from Theorem \ref{thm:mainthmUTAwafd}.  Hence we are only left to
prove (1) $\Longrightarrow$ (3).

Let $A \subset B(H)$ be any faithful representation of $A$. Then
the inclusion $A \hookrightarrow B(H)$ is a nuclear map (cf.\
\cite[Theorem 7.3]{wassermann:exactbook}).  From part (5) of
Theorem \ref{thm:amenabletraces} there exists a u.c.p.\ map $\Phi
: B(H) \to \pi_{\tau} (A)^{\prime\prime}$ such that $\Phi(a) =
\pi_{\tau} (a)$, for all $a \in A$.  It follows that $\pi_{\tau}$
(which can be identified with $\Phi|_A$ composed with $A
\hookrightarrow B(H)$) is a nuclear map.
\end{proof}

It is an important open problem to determine which quasidiagonal
C$^*$-algebras satisfy the equation $\TA = \TAQD$.  Indeed, if $A
= C^*({\mathbb F}_{\infty})$ then the equation $\TA = \TAQD$ is
equivalent to Connes' embedding problem (see Section
\ref{thm:connesproblem}). If $A$ is nuclear and quasidiagonal then
this equation is predicted by Elliott's conjecture (see
Proposition \ref{thm:elliott}). When we first began studying
amenable traces we thought that it may be possible to prove
Connes' problem by verifying the equation $\TA = \TAQD$ for all
quasidiagonal C$^*$-algebras.  Support for this strategy is
provided by the following observation.

\begin{prop}
If $A$ is quasidiagonal and $\tau \in \TA$ then there exists a sequence of
u.c.p.\ maps $\phi_n : A \to M_{k(n)}({\mathbb C})$, which are
asymptotically multiplicative in norm, and a sequence of {\em states}
$\sigma_n \in S(M_{k(n)}({\mathbb C}))$ such that $\sigma_n \circ
\phi_n (a) \to \tau(a)$ for all $a \in A$.
\end{prop}

\begin{proof} In fact, this result holds for all states on $A$.
Indeed, let $\phi_n : A \to M_{k(n)}({\mathbb C})$ be any sequence
of u.c.p.\ maps such that $\| a \| = \lim \| \phi_n (a) \|$ for
all $a \in A$.  Consider the following set of states on $A$:
$${\mathcal S} = \{ \sigma \in S(A): \exists \sigma_n \in
S(M_{k(n)}({\mathbb C})) \ \mathrm{such} \ \mathrm{that} \ \sigma
= \sigma_n \circ \phi_n \}.$$ Since $\| a \| = \lim \| \phi_n (a)
\|$ for all $a \in A$ it is easy to see that for every
self-adjoint $a \in A$ we have $\| a \| = \sup \{ |\sigma(a)| :
\sigma \in {\mathcal S} \}$. From the Hahn-Banach theorem it
follows that the weak-$*$ closure of ${\mathcal S}$ is the entire
state space of $A$ (cf.\ Lemma \ref{thm:HB}).  An argument similar
to that given in the proof of Proposition \ref{thm:convexity}
completes the proof.
\end{proof}

Hence we wondered if one could always replace the {\em states} in the
lemma above with traces by passing to centralizers or some other
general argument.  Unfortunately this is impossible, in general,
though the nuclear case and the case of $C^*({\mathbb F}_{\infty})$
are still open.

\begin{cor}
\label{thm:exactQD} There exists an exact, residually finite
dimensional\footnote{$A$ is residually finite dimensional if it
has a separating family of finite dimensional representations.
Another way of saying this is that $A$ embeds into $\Pi
M_{k(n)}(\mathbb{C})$ for some sequence of integers $k(n)$.}
(hence quasidiagonal) C$^*$-algebra $A$ such that $\TA \neq \TAim$
(hence $\TA \neq \TAQD$).
\end{cor}

\begin{proof}
Since $C^*_r ({\mathbb F}_2)$ is exact and every exact
C$^*$-algebra is a quotient of an exact, residually finite
dimensional algebra (cf.\ \cite{brown:QDsurvey}) it follows that
we can find an exact, residually finite dimensional C$^*$-algebra
$A$ and a trace $\tau \in \TA$ such that
$\pi_{\tau}(A)^{\prime\prime}$ is a free group factor. Since free
group factors are not hyperfinite, $\tau$ is not an amenable trace
on $A$.
\end{proof}

\section{Type I C$^*$-algebras}

Recall from Theorem \ref{thm:nuclear} that {\em every trace on a
nuclear C$^*$-algebra is a uniform amenable trace.} It may be that
this is the strongest approximation property enjoyed by the entire
class of nuclear C$^*$-algebras.  However, in this section we
will show that for type I C$^*$-algebras a much stronger
approximation property always holds -- every trace is uniform
locally finite dimensional.

Our first lemma is rather technical but gives a useful
characterization of these traces.

\begin{lem}
\label{thm:localfinitedimensional} Let $\tau \in \TA$ be given.
Then $\tau \in \TAulfd$ if and only if for every finite set
${\mathfrak F} \subset A$ and $\epsilon > 0$ there exists a
C$^*$-subalgebra $B \subset A$ and a u.c.p.\ map $\phi : B \to
M_n({\mathbb C})$ such that $d({\mathfrak F}, B_{\phi}) <
\epsilon$ and $\| \tau|_B - \mathrm{tr}_n\circ\phi \|_{B^*} <
\epsilon$.
\end{lem}

\begin{proof}
Evidently we need only show the `if' part.  The obvious approach
to proving this lemma would be to simply extend, via Arveson's
theorem, the map on $B$ to all of $A$.  However, it does not seem
obvious that if $\tilde{\phi}:A \to M_n({\mathbb C})$ is a u.c.p.\
extension then we still have $\| \tau -
\mathrm{tr}_n\circ\tilde{\phi} \|_{A^*} < \epsilon$ and hence our
proof will take another route.

So, assume ${\mathfrak F} \subset A$ and $\epsilon > 0$ are given
and choose $B \subset A$ and $\phi : B \to M_n({\mathbb C})$ such
that $d({\mathfrak F}, B_{\phi}) < \epsilon$ and $\| \tau|_B -
\mathrm{tr}_n\circ\phi \|_{B^*} < \epsilon$.

Let $P_{\tau} \in A^{**}$ (resp.\ $P_{\mathrm{tr}_n\circ\phi} \in
B_{\phi}^{**} \subset B^{**} \subset A^{**}$) be the central cover
(cf.\ \cite{pedersen:book}) of the GNS representation of $\tau$
(resp.\ the GNS representation of
$\mathrm{tr}_n\circ\phi|_{B_{\phi}}$). Let $Q =
P_{\tau}P_{\mathrm{tr}_n\circ\phi} \in A^{**}$.  Then $Q$ commutes
with $B_{\phi} \ (\subset A^{**})$ and $QB_{\phi}$ is a finite
dimensional subalgebra of $QA^{**}Q \subset P_{\tau}A^{**} \cong
\pi_{\tau}(A)^{\prime\prime}$ (since
$P_{\mathrm{tr}_n\circ\phi}B_{\phi}$ is finite dimensional).
Letting $\tau^{**}$ denote the normal extension of $\tau$ to
$A^{**}$, note that $\tau^{**}(Q) =
\tau^{**}(P_{\mathrm{tr}_n\circ\phi}) =
\tau^{**}|_{B_{\phi}^{**}}(P_{\mathrm{tr}_n\circ\phi}) =
(\tau|_{B_{\phi}})^{**}(P_{\mathrm{tr}_n\circ\phi}) > 1 -
\epsilon$ since $\| \tau|_{B_{\phi}} -
\mathrm{tr}_n\circ\phi|_{B_{\phi}} \|_{B_{\phi}^*} < \epsilon$.
Let $E : QA^{**}Q \to QB_{\phi}$ be a conditional expectation
which preserves the tracial state $\theta$ on $QA^{**}Q$ defined
by $QxQ \mapsto \frac{1}{\tau^{**}(Q)}\tau^{**}(QxQ)$.  Then
define $\Phi : A \to QB_{\phi}$ by $\Phi(a) = E(QaQ)$.  Note that
$B_{\phi} \subset A_{\Phi}$ and hence the proof is complete once
one observes that
$$|\tau(a) - \theta\circ\Phi(a)| = |\tau^{**}(Q^{\perp}a) +
\tau^{**}(Qa) - \frac{1}{\tau^{**}(Q)}\tau^{**}(Qa)| \leq
(\epsilon^{1/2} + \frac{\epsilon}{1 - \epsilon}) \| a \|,$$ and
then applies Lemma \ref{thm:tracepreservingembedding} to get from
$QB_{\phi}$ to a full matrix algebra.
\end{proof}

The following lemma of Huaxin Lin provides the key first step to
proving that all traces on type I C$^*$-algebras are uniform
locally finite dimensional.

\begin{lem}(cf.\ \cite[Lemma 4.7]{lin:ACtraces}) Let $B$ be a
subhomogeneous C$^*$-algebra (i.e.\ assume the dimension of every
irreducible representation of $B$ is uniformly bounded).  Then,
for every trace $\tau \in \TB$, finite set ${\mathfrak F} \subset B$
and $\epsilon > 0$ there exists an ideal $I \subset B$, a trace
$\gamma \in {\rm T(B/I)}$ and a unital, finite dimensional subalgebra $C
\subset B/I$ such that (1) $\| \tau - \gamma\circ \sigma \|_{B^*} <
\epsilon$ and (2) $d(\sigma({\mathfrak F}), C) < \epsilon$, where $\sigma :
B \to B/I$ is the canonical quotient mapping.
\end{lem}

\begin{lem}
\label{thm:whocares}
Let $\tau \in \TA$ be a trace and assume that the weak closure of
the associated GNS representation of $A$ is a type I von Neumann algebra.
Then $\tau \in \TAulfd$.
\end{lem}

\begin{proof} Fix a finite set ${\mathfrak F} \subset A$ and $\epsilon > 0$.
Let $\pi_{\tau} : A \to B(L^2(A,\tau))$ be the GNS representation
arising from $\tau$.  Since $\pi_{\tau}(A)^{\prime\prime}$ is a
type I von Neumann algebra we have two cases to consider.

{\noindent Case 1:} Assume $\pi_{\tau}(A)^{\prime\prime}$ is
isomorphic to $$\prod_{i = 1}^{k} M_{n(i)} ({\mathbb C}) \otimes
L^{\infty}(X_i, \mu_i)$$ for some probability spaces $(X_i,
\mu_i)$ and some natural number $k$.  Let $B = \pi_{\tau}(A)$ and
$\tilde{{\mathfrak F}} = \pi_{\tau}({\mathfrak F}) \subset B$.
Then $B$ is subhomogeneous and by Lin's lemma we can find an ideal
$I \subset B$, trace $\gamma \in {\rm T(B/I)}$ and unital, finite
dimensional subalgebra $C \subset B/I$ such that $\| \tau -
\gamma\circ\sigma \|_{B^*} < \epsilon$ and
$d(\sigma(\tilde{{\mathfrak F}}), C) < \epsilon$. Note that if $D
= (\sigma\circ\pi_{\tau})^{-1}(C) \subset A$ is the pullback of
$C$ then $d({\mathfrak F},D) < \epsilon$. Note also that $\|
\tau|_D - \gamma\circ\sigma\circ\pi_{\tau}|_D \|_{D^*} < \epsilon$
in this situation.  Applying Lemma
\ref{thm:localfinitedimensional} to the subalgebra $D \subset A$
and $*$-homomorphism $\sigma\circ\pi_{\tau}: D \to C$ we get the
result in this case. (Strictly speaking one should apply Lemma
\ref{thm:tracepreservingembedding} to get from $C$ to a full
matrix algebra, but we have seen this type of argument numerous
times already.)

{\noindent Case 2:} The only other possibility is that the weak
closure of $\pi_{\tau} (A)$ is isomorphic to $$\prod_{i =
1}^{\infty} M_{n(i)} ({\mathbb C}) \otimes L^{\infty}(X_i,
\mu_i).$$ However, in this case $$\tau = \sum_{i=1}^{\infty}
\theta_i \tau_i$$ where $\theta_i$ are positive real numbers which
sum up to one and $\tau_i$ is gotten by restricting $\tau$ to the
$i^{th}$ summand of $\prod_{i = 1}^{\infty} M_{n(i)} ({\mathbb C})
\otimes L^{\infty}(X_i, \mu_i)$ (and renormalizing to get a
state).\footnote{This is because of normality of the induced
vector trace in the GNS representation.}  However, each of the
tracial functionals $$\eta_k = \sum_{i=1}^{k} \theta_i \tau_i$$
falls under Case 1 and the $\eta_k$'s converge in norm to $\tau$.
In other words, $\tau$ is a norm limit of elements of $\TAulfd$
and since this set is norm closed it follows that $\tau \in
\TAulfd$ as well.
\end{proof}

\begin{cor}
\label{thm:typeI} If $A$ is an inductive limit of type I
C$^*$-algebras then $\TA = \TAulfd$. In particular, every trace on
a type I C$^*$-algebra is uniform locally finite dimensional.
\end{cor}

\begin{proof}
Let $\tau \in \TA$ be given and assume $A$ is the norm closure of
type I subalgebras $\{A_i\}$.  By the previous lemma we have that
$\tau|_{A_i}$ is uniform locally finite dimensional and hence from
Lemma \ref{thm:localfinitedimensional} it follows that $\tau \in
\TAulfd$.
\end{proof}

Note, of course, that one really does not need an inductive limit
decomposition in the previous corollary.  It suffices to be
``locally type I" in the sense that every finite set is almost
contained in a type I subalgebra (though the type I subalgebras
need not be increasing).

\section{Tracially AF C$^*$-algebras}

This section is devoted to Huaxin Lin's class of tracially AF
algebras \cite{lin:TAF}.  This class of C$^*$-algebras is similar
to the class of Popa algebras, however Lin's work on classifying
such algebras has been a breakthrough in the classification
program and proving that particular C$^*$-algebras are tracially
AF now immediately yields new classification results.

Our first goal is to show that tracial approximation properties give a
very simple characterization of these algebras in many instances.
This was also recognized by Lin in \cite{lin:ACtraces} where he
introduced a notion of `approximately AC trace' (cf.\ \cite[Definition
3.1]{lin:ACtraces}) and used this to characterize certain tracially AF
algebras.  The present approach, however, has the advantage of being
conceptually and technically easier to digest. Moreover, our
definitions do not presuppose an inductive limit type structure as is
required in \cite[Definition 3.1]{lin:ACtraces}.

The following theorem summarizes several of Lin's results and
explains the `tracially AF' terminology (cf.\ \cite[Theorems 6.9, 6.11
and 6.13]{lin:tracialtopologicalrank}).

\begin{thm}
\label{thm:Linscharacterization} Let $A$ be a simple
C$^*$-algebra.  Then $A$ is tracially AF if and only if $A$ has
real rank zero, stable rank one, weakly unperforated
K-theory\footnote{Recall that real rank zero means the
self-adjoint elements with finite spectrum are dense in the set of
all self-adjoints, stable rank one means the invertible elements
are dense in the whole algebra and weakly unperforated K-theory
means that if $x \in K_0(A)$ and $\tau(x) > 0$ for all $\tau \in
\TA$ then $x > 0$.} and for every finite subset ${\mathfrak F}
\subset A$ and $\epsilon > 0$ there exists a finite dimensional
subalgebra $B \subset A$ with unit $e$ such that:

\begin{enumerate}
\item $\| [x,e] \| < \epsilon$, for all $x \in {\mathfrak F}$.

\item $d(e{\mathfrak F}e, B) < \epsilon$.

\item $\tau(e) > 1 - \epsilon$ for all $\tau \in \TA$.
\end{enumerate}
\end{thm}

As mentioned above, we are primarily interested in knowing which
C$^*$-algebras are tracially AF.  We remind the reader that Popa's
work shows that every simple, quasidiagonal C$^*$-algebra with real
rank zero is `almost' tracially AF (more precisely, satisfies the
approximation property with (1) and (2), but not necessarily (3),
above).  The point of this section is that if one happens to know that
the traces on $A$ are sufficiently well behaved then it can be shown
that $A$ is tracially AF (i.e.\ we can also arrange condition (3)).
The difficulty is that condition (3) is global (i.e.\ must hold for
all traces) and hence if one only knows approximation properties of
particular traces then one is forced to add restrictions on the size
of $\TA$ so that a constructive procedure can be carried out.  The
following lemma, which for all intents and purposes is due to Lin,
makes this more precise.

\begin{lem}
\label{thm:traciallylocal}
Assume that $\TA$ is $\| \cdot \|_{A^*}$-separable and
every hereditary subalgebra $A_0 \subset A$ has the following
(tracially local) approximation property: For every trace $\tau \in
\TA$, finite subset ${\mathfrak F} \subset A_0$ and $\epsilon > 0$
there exists a finite dimensional subalgebra $B \subset A_0$ with unit
$e$ such that:

\begin{enumerate}
\item $\| [x,e] \| < \epsilon$, for all $x \in {\mathfrak F}$.

\item $d(e{\mathfrak F}e, B) < \epsilon$.

\item $\tau(e) > (1 - \epsilon)\| \tau|_{A_0} \|$.
\end{enumerate}

Then it follows that for every finite subset ${\mathfrak F} \subset A$
and $\epsilon > 0$ there exists a finite dimensional subalgebra $B
\subset A$ with unit $e$ such that:

\begin{enumerate}
\item[($1^{\prime}$)] $\| [x,e] \| < \epsilon$, for all $x \in {\mathfrak F}$.

\item[($2^{\prime}$)] $d(e{\mathfrak F}e, B) < \epsilon$.

\item[($3^{\prime}$)] $\tau(e) > 1 - \epsilon$ for all $\tau \in \TA$.
\end{enumerate}
In other words, if one assumes that $\TA$ is $\| \cdot
\|_{A^*}$-separable and happens to know that for each trace one
can find a ``large" subalgebra with the right tracially AF
properties then, in fact, it is possible to find a subalgebra
which is large in all traces simultaneously.\footnote{It would be
very nice if one could remove the assumption of $\| \cdot
\|_{A^*}$-separability from this result -- not just for aesthetic
reasons but also because it would immediately imply significant
new classification results.  While we believe this should be
possible it has resisted our best efforts and those of a few other
colleagues we have spoken to.}
\end{lem}

\begin{proof} This lemma is basically an abstraction of the proof of
\cite[Theorem 4.13]{lin:ACtraces}.  Hence we will only sketch the main
idea.

Fix a finite subset ${\mathfrak F} \subset A$ and $\epsilon > 0$.  Let
$\{ \tau_i \} \subset \TA$ be a {\em norm} dense sequence.  Applying
the assumed approximation property to $({\mathfrak F}, \tau_1,
\epsilon/2)$ we can find $B_1$ with unit $e_1$ such that

\begin{enumerate}
\item[(a)] $\| [x,e_1] \| < \epsilon/2$, for all $x \in {\mathfrak F}$,

\item[(b)] $d(e_1{\mathfrak F}e_1, B_1) < \epsilon/2$,

\item[(c)] $\tau_1(e_1) > 1 - \epsilon/2$.
\end{enumerate}

Now we look at the hereditary subalgebra which is orthogonal to $B_1$
(i.e.\ $e_1^{\perp} A e_1^{\perp}$) and apply the assumed
approximation property to $(e_1^{\perp}{\mathfrak F}e_1^{\perp},
\tau_2, \epsilon/8)$.  We get a finite dimensional subalgebra $B_2
\subset e_1^{\perp} A e_1^{\perp}$ with unit $e_2$ satisfying the
stated conditions.  The key observation is that the finite dimensional
subalgebra $B_1 \oplus B_2$ has the following property:

\begin{enumerate}
\item[(d)] $\| [x,e_1\oplus e_2] \| < \epsilon/2 + \epsilon/8$, for
all $x \in {\mathfrak F}$,

\item[(e)] $d([e_1\oplus e_2]{\mathfrak F}[e_1\oplus e_2], B_1
\oplus B_2) < \epsilon/2 + \epsilon/4$,

\item[(f)] $\tau_2(e_1\oplus e_2) > 1 - \epsilon/8$.
\end{enumerate}

By induction we can produce a sequence of finite dimensional
subalgebras $C_n$\footnote{$C_1 = B_1$, $C_2 = B_1 \oplus B_2,
\ldots$.} with units $f_n$ such that $f_1 \leq f_2 \leq f_3 \leq
\cdots$ and

\begin{enumerate}
\item[(g)] $\| [x,f_n] \| < \epsilon$, for
all $x \in {\mathfrak F}$,

\item[(h)] $d(f_n{\mathfrak F}f_n, C_n) < \epsilon$,

\item[(i)] $\tau_n(f_n) > 1 - \frac{\epsilon}{2^n}$.
\end{enumerate}

The only thing left to verify is that $\tau(f_n) \to 1$ for every
$\tau \in \TA$ because if this is the case then Dini's theorem applies
(identifying $\{f_n \}$ with an {\em increasing} sequence of functions
on the compact metric space $\TA$) and we conclude that the
convergence is uniform; i.e.\ there exists some large $n$ such that
$\tau(f_n) > 1 - \epsilon$ for all $\tau \in \TA$.  This is where the
assumption of $\| \cdot \|_{A^*}$-separability of $\TA$ comes in.
Indeed, we may assume that the sequence $\{ \tau_i \} \subset \TA$ has
infinite multiplicity (i.e.\ every trace in the sequence appears
infinitely many times in the sequence) and thus condition (i) ensures
that $\tau(f_n) \to 1$ on a {\em norm} dense subset of $\TA$.  However
this implies pointwise convergence to one on all of $\TA$ and we are
done.
\end{proof}

\begin{lem} Let $\tau \in \TAulfd$ and a projection $p \in A$ be
given.  Then the tracial state $\frac{1}{\tau(p)} \tau|_{pAp}$ on
$pAp$ is also uniform locally finite dimensional.
\end{lem}

\begin{proof} Let $p \in {\mathfrak F} \subset pAp$ and
$\epsilon > 0$ be given.  Let $\delta > 0$ be much smaller (to be
determined later) than $\epsilon$ and let $\phi : A \to M_k
({\mathbb C})$ be a u.c.p.\ map such that $\| \mathrm{tr}_k \circ
\phi - \tau \| < \tau(p) \delta$ and ${\mathfrak F}$ is
$\delta$-contained in the multiplicative domain $A_{\phi}$ of
$\phi$.  Since $p$ is a projection it follows from some standard
K-theoretic results that we can find a projection $q \in
A_{\phi}$, which is very close to $p$, and a unitary $u \in A$,
which is very close to the identity, such that $upu^* = q$. Hence
defining $\psi : pAp \to \phi(q) M_k ({\mathbb C})\phi(q) \cong
M_l ({\mathbb C})$ (for some $l \leq k$) by the formula $\psi(pap)
= \phi(upapu^*)$ we see that ${\mathfrak F}$ is nearly contained
in the multiplicative domain of $\psi$ (since $q \in A_{\phi}$ and
$u$ is close in norm to the identity, $u^* qA_{\phi}q u$ is in the
multiplicative domain of $\psi$ and almost contains ${\mathfrak
F}$) and $\mathrm{tr}_l \circ \psi$ is very close in norm (within
$\delta$ actually) to $\frac{1}{\tau(p)} \tau|_{pAp}$.
\end{proof}

We will need one more lemma before coming to our characterization
of tracially AF algebras in terms of tracial approximation
properties.  The following result is taken from the work of Huaxin
Lin (cf.\ \cite[proof of Theorem 5.3]{lin:TAF}) but we include
those ideas in the proof which are no so easy to extract from
\cite{lin:TAF}.

\begin{lem}
\label{thm:fdsplitting} Assume that $A$ has real rank zero, $1_A
\in D \subset A$ is a C$^*$-subalgebra which contains an ideal $I
\triangleleft D$ such that $I$ is a hereditary subalgebra of $A$
and $D/I = C$ is finite dimensional.  Then the exact sequence
$$0 \to I \to D \stackrel{\pi}{\to} C \to 0$$ has a $*$-homomorphic
splitting -- i.e.\ there is a $*$-homomorphism $\sigma:C \to D$
such that $\pi\circ\sigma = id_C$.  Moreover, there is an
approximate unit of projections $\{q_n \} \subset I$ such that
$[q_n, \sigma(C)] = 0$ for all $n \in \mathbb{N}$.
\end{lem}

\begin{proof} Once one knows that projections from $C$ can always
be lifted to $D$ then the existence of a $*$-homomorphic splitting
follows exactly as in the proof of \cite[Lemma III.6.2]{davidson}.
However, a nice proof that one can always lift projections can be
found in \cite[Lemma 5.2]{lin:TAF} and hence our map $\sigma$
always exists.

That an approximate unit exists which both consists of projections
and commutes with $\sigma(C)$ is a bit technical but only relies
on fairly standard matrix manipulations.  The main ingredients are
already contained in the special case that $C \cong
M_2(\mathbb{C})$ (where notation becomes much simpler) so we give
all the details in this setting and leave the general case to the
reader.

So let $B = \sigma(C) \subset A$ and denote by $e_1$ and $e_2$ a
pair of orthogonal minimal projections in $B$.  Let $u \in B$ be a
unitary such $ue_1 u^* = e_2$.  Consider the ideal $e_1 I e_1
\triangleleft e_1 D e_1$. Since real rank zero passes to
hereditary subalgebras it follows that $e_1 I e_1$ has real rank
zero and hence we can find an approximate unit of projections $\{
p_n \} \subset e_1 I e_1$.  Evidently $\{ up_n u^* \} \subset e_2
I e_2$ is also an approximate unit of projections and hence we can
define projections $$q_n = r_n + p_n + up_n u^* \in I,$$ where $\{
r_n \} \subset (e_1 + e_2)^{\perp} I (e_1 + e_2)^{\perp}$ is also
an approximate unit of projections. We must verify that these
projections form an approximate unit for $I$ but the commutation
part is now simple since one easily checks that each $q_n$
commutes with all four of the operators $e_1, e_2, ue_1$ and
$e_1u^*$ (and these four operators are a set of matrix units for
$B$).

In general, if $I \triangleleft D$ is an ideal, $\{ s_n \} \subset
dId$ is an approximate unit ($d \in D$ is some fixed element) and
$x \in I$ is arbitrary then a straightforward calculation using
the C$^*$-identity shows that $$\| (xd)s_n - xd \| \to 0.$$  Now,
our goal is to show that $q_n = r_n + p_n + up_n u^*$ is an
approximate unit of $I$ and this essentially reduces to the remark
in the preceding sentence.  Indeed, we can write each $x \in I$ as
a 3$\times$3 matrix with respect to the decomposition $1 = (e_1 +
e_2)^{\perp} + e_1 + e_2$.  It then suffices to consider the ``off
diagonal elements" and show, for example, $$\| e_1 x (e_1 +
e_2)^{\perp} q_n - e_1 x (e_1 + e_2)^{\perp} \| = \| e_1\big( x
(e_1 + e_2)^{\perp} r_n - x (e_1 + e_2)^{\perp}\big) \| \to 0.$$
This follows from the general remarks above and hence we leave the
remaining details to the reader.

For more general finite dimensional algebras $B$ we first
decompose $B$ as a finite direct sum of full matrix algebras and
then apply the procedure above underneath each of the summands. It
is quite messy to properly write down so we won't attempt it.
\end{proof}

\begin{prop}
\label{thm:characterization}
Let $A$ be a simple C$^*$-algebra with real rank zero,
stable rank one, weakly unperforated K-theory and assume that $\TA$ is
norm separable.  Then $A$ is tracially AF if and only if $\TA =
\TAulfd$.
\end{prop}

\begin{proof} We begin with the `only if' statement.  So assume that
$A$ is tracially AF and let $\tau \in \TA$, a finite set ${\mathfrak
F} \subset A$ and $\epsilon > 0$ be given.  Choose a finite
dimensional subalgebra $B \subset A$ as in Lin's characterization
above and note that there is a conditional expectation $eAe \to B$
which preserves the tracial state $\frac{1}{\tau(e)}\tau|_{eAe}$.
Composing this map with the u.c.p.\ map $A \to eAe$, $a \mapsto eae$
we get a u.c.p.\ map from $A$ to $B$ and ${\mathfrak F}$ is nearly
contained in the multiplicative domain of this map since condition (1)
in Theorem \ref{thm:Linscharacterization} above implies that $a
\approx eae \oplus e^{\perp}ae^{\perp}$ for all $a \in {\mathfrak F}$
(and clearly $e^{\perp}Ae^{\perp}$ belongs to the multiplicative
domain).  Finally, it is easily seen that $\tau$ is close (in norm) to
$\frac{1}{\tau (e)} \tau(e \cdot e)$ since $\tau(e)$ is close to one.

For the converse, assume that $\TA = \TAulfd$ and we will show
that $A$ has the `tracially local' approximation property from
Lemma \ref{thm:traciallylocal}. Note that from the proof of that
lemma we really only need to consider corners of $A$ (instead of
general hereditary subalgebras).  Moreover, by the previous lemma,
restrictions of traces on $A$ will again be  uniform locally
finite dimensional and hence we may assume that the hereditary
subalgebra is $A$ itself (since the proof will carry over verbatim
to corners). Hence we need to show that if ${\mathfrak F} \subset
A$ is a finite set, $\tau \in \TA$ is arbitrary and $\epsilon > 0$
is given then we can find a finite dimensional subalgebra $B
\subset A$ with unit $e$ such that $\tau(e)$ is large, $e$ almost
commutes with ${\mathfrak F}$ and cutting down by $e$ almost
pushes ${\mathfrak F}$ into $B$.

By assumption we can find a u.c.p.\ map $\phi : A \to M_k
({\mathbb C})$ such that $\| \mathrm{tr}_k \circ \phi - \tau
\|_{A^*} < \epsilon$ and ${\mathfrak F}$ is nearly contained
(within $\epsilon$) in the multiplicative domain $A_{\phi}$ of
$\phi$. Let $J \subset A_{\phi}$ be the kernel of the
$*$-homomorphism $\phi|_{A_{\phi}}$.  We claim that $J$ is a
hereditary subalgebra of $A$ and hence inherits real rank zero
from $A$.  To see this we assume $0 \leq x \leq y$ and $y \in J$.
By positivity we have that $\phi(x) = 0$.  The Cauchy-Schwartz
inequality for c.p.\ maps tells us that $0 \leq \phi(\sqrt{x})^2
\leq \phi((\sqrt{x})^2) = 0$ and hence $\phi(\sqrt{x}) = 0$ as
well.   This shows that $\sqrt{x} \in A_{\phi}$ and thus $x$
belongs to the multiplicative domain too.  It follows that $x \in
J$ which shows that $J$ is a hereditary subalgebra of $A$.

Applying Lemma \ref{thm:fdsplitting} to the short exact sequence
$$0 \to J \to A_{\phi} \to C \to 0,$$ where $C = \phi(A_{\phi})
\subset M_k ({\mathbb C})$, we can find a finite dimensional
C$^*$-subalgebra $B \subset A_{\phi}$, with unit $e$, such that
$\phi$ is an isomorphism from $B$ onto $C$.  Moreover, we can find
an approximate unit of projections $\{q_n\} \subset J$ which
commutes with $B$ and hence is quasicentral in $A_{\phi}$.
Borrowing a few more techniques from the proof of \cite[Theorem
5.3]{lin:TAF} will complete the proof.  The main new ingredient,
however, we wish to explain from the outset.  Namely, note that
$$\tau(\big(1 - q_n\big) e) > 1 - \epsilon$$ for every $n$ since $\|
\mathrm{tr}_k \circ \phi - \tau \|_{A^*} < \epsilon$, $\big(1 -
q_n\big) e$ belongs to the unit ball of $A$ and, finally,
$\phi(\big(1 - q_n\big) e) = 1$ since $q_n \in J$ and $e$ is a
lift of the unit of $C$.

Hence the proof will be complete as soon as we verify:
\begin{enumerate}
\item $\limsup \| [x, \big(1 - q_n\big) e ] \| < \epsilon$ for all
$x \in \mathfrak{F}$.

\item $\limsup d(\big(1 - q_n\big) e x \big(1 - q_n\big) e, \big(1
- q_n\big) e B) < \epsilon$ for all $x \in \mathfrak{F}$.
\end{enumerate}
Indeed, if we are successful in showing this then we can take our
desired finite dimensional algebra to be $\big(1 - q_n\big) B$ for
any sufficiently large $n$.  Since $\mathfrak{F}$ is almost
contained in $A_{\phi}$ it will be sufficient to replace $x$ above
by an element in the multiplicative domain and show that the
$\limsup$'s are actually zero in this case.   But if $x \in
A_{\phi}$ then $$\| x \big(1 - q_n\big) e - \big(1 - q_n\big) e
x\| \stackrel{\epsilon}{\approx} \|(1 - q_n)(xe - ex) \| \to
\|\phi(xe - ex)\| = 0$$ where we used the fact that $\{q_n\}$ is
quasicentral to get the first approximation and the fact that $e$
is a lift of $1$ (which is central in the quotient, of course) for
the second part. For the assertion $d(\big(1 - q_n\big) e x \big(1
- q_n\big) e, \big(1 - q_n\big) e B) \to 0$ we first pick a lift
$b \in B$ of $\phi(x) \in C$.  A similar argument to the one given
above shows that $d(\big(1 - q_n\big) e x \big(1 - q_n\big) e,
\big(1 - q_n\big) e b) \to 0$ and hence the proof is complete.
\end{proof}

\begin{rem} This proposition is the analogue of Theorems 3.8
and 4.13 from \cite{lin:ACtraces}.  Though our proof certainly
uses some key ideas from \cite{lin:ACtraces} it is also true that
the present approach avoids a number of technical difficulties
which must be dealt with in Lin's approach.  Note that the `only
if' part holds without the separability assumption on the tracial
space. Also, the `if' statement actually generalizes \cite[Theorem
4.13]{lin:ACtraces} as it can be shown that every `approximately
AC trace' is uniform locally finite dimensional. Indeed, it is not
too hard to show that every `AC trace' is uniform locally finite
dimensional and so Lemma \ref{thm:localfinitedimensional} implies
that approximately AC traces have this property as well.
\end{rem}

Though it won't be needed in what follows we wish to point out that
the (annoying) assumption that $\TA$ is norm separable is basically
necessary if one wishes to carry out a constructive procedure as in
the proof of Lemma \ref{thm:traciallylocal}.  Moreover, we remind the
reader that a constructive procedure is also basically necessary in
the category of C$^*$-algebras (i.e.\ one can't use Zorn's Lemma type
maximality arguments -- which work so beautifully in von Neumann
algebras -- since an infinite sum of projections doesn't make sense in
a C$^*$-algebra).

\begin{prop} Let $A$ be a simple, nuclear C$^*$-algebra with real
rank zero, stable rank one and weakly unperforated K-theory.  Then the
following statements are equivalent:

\begin{enumerate}
\item $\TA$ is weakly (i.e.\ $\sigma(A^*, A^{**})$) separable.

\item $\TA$ is norm separable.

\item $\TA$ has countably many extreme points.

\item There exists a countable subset ${\mathcal S} \subset \TA$ with
the following property: For every increasing sequence of projections
in $A$, say $p_1 \leq p_2 \leq \cdots$, if one knows that $\tau(p_n)
\to 1$ for all $\tau \in {\mathcal S}$ then it follows that $\tau(p_n)
\to 1$ for all $\tau \in \TA$.
\end{enumerate}
\end{prop}

\begin{proof} The equivalence of (1) and (2) follows from the
Hahn-Banach Theorem and the convexity of $\TA$.

(3) $\Longrightarrow$ (2).  By Choquet theory, for every $\tau \in
    \TA$ there exist non-negative real numbers $\{ a_i \}$ such that
    $\sum a_i =1$ and if we list the extreme points of $\TA$ as $\{
    \tau_1, \tau_2, \ldots \}$ then $$\tau = \sum a_i \tau_i.$$
    Evidently this implies $\TA$ is norm separable.

(2) $\Longrightarrow$ (3). We prove the contrapositive.  Note that
    since extreme points of $\TA$ yield factors in their GNS
    representations and, in general, factor representations are either
    disjoint or equivalent (cf.\ \cite[3.8.13]{pedersen:book}) it
    follows (using the fact that finite factors have a unique trace)
    that for two extreme traces $\tau_1, \ \tau_2$ we either have that
    $\tau_1 = \tau_2$ or $\| \tau_1 - \tau_2 \|_{A^*} = 2$.  Hence if
    $\TA$ has uncountably many extreme points we see that it can't be
    norm separable.

(2) $\Longrightarrow$ (4). Any norm dense subset has the property
    described in condition (4).

(4) $\Longrightarrow$ (2). We again prove the contrapositive.  So
    assume $\TA$ is not separable in norm and let $\pi = \oplus_{\tau
    \in \TA} \pi_{\tau} : A \to B(\oplus_{\tau \in \TA} L^2(A,\tau))$
    be the direct sum of all GNS representations arising from traces
    on $A$.  Let $M = \pi(A)^{\prime\prime}$ and note that $M$ is a
    finite von Neumann algebra and every trace on $A$ naturally
    extends to a weakly continuous trace on $M$.  Hence the predual of
    $M$ is not separable (since we assume $\TA$ is not separable).
    Thus $M$ has no faithful representation on a separable Hilbert
    space.

Now let ${\mathcal S} \subset \TA$ be any countable subset and we will
show that there exists an increasing sequence of projections $p_1 \leq
p_2 \leq \cdots$ such that $\tau(p_n) \to 1$ for all $\tau \in
{\mathcal S}$ but such that there exists some $\gamma \in \TA$ for
which $\gamma(p_n)$ does not tend to one.  The first step in
constructing such projections is to construct an increasing sequence
of positive elements from the unit ball of $A$ with this
property. This will follow from Pedersen's Up-Down Theorem (cf.\
\cite[Theorem 2.4.3]{pedersen:book}).

Let $\sigma_{{\mathcal S}} : M \to B(\oplus_{\tau \in {\mathcal S} }
L^2(A,\tau))$ be the natural restriction of $M$ to the sum of the GNS
spaces coming from traces in ${\mathcal S}$.  Note that $\oplus_{\tau
\in {\mathcal S} } L^2(A,\tau)$ is a separable Hilbert space and hence
$\sigma$ can't be faithful.  Thus there exists a central projection $q
\in M$ such that $qM$ is naturally isomorphic to $\sigma(M)$.
Evidently we can then find a trace $\gamma \in \TA$ such that
$\gamma(q) = 0$ ($\gamma$ also denoting the natural extension to $M$).
Note also that $\tau(q) = 1$ for all $\tau \in {\mathcal S}$.  Now let
$\sigma_{\tilde{{\mathcal S}}} : M \to B(\oplus_{\tau \in
\tilde{{\mathcal S}}} L^2(A,\tau))$.  We will identify $qM$ with a
(nontrivial) direct summand of $\sigma_{\tilde{{\mathcal S}}}(M)$.
Note that since $A$ is simple, $\sigma_{\tilde{{\mathcal S}}}(A) \cong
A$.  At this point what we have arranged is: There exists a
representation of $A$, $\sigma : A \to B(H)$, on a separable Hilbert
space, such that every trace in $\tilde{{\mathcal S}}$ naturally
extends to a weakly continuous trace on $\sigma(A)^{\prime\prime}$ and
there exists a projection $q \in \sigma(A)^{\prime\prime}$ such that
$\gamma(q) = 0$ while $\tau(q) = 1$ for all $\tau \in {\mathcal S}$
(these are the abstract properties we really need).

Since $A \cong \sigma(A)$ is weakly dense in
$\sigma(A)^{\prime\prime}$ we can apply Pedersen's Up-Down Theorem to
find a sequence of positive elements from the unit ball of $A$, say
$\{ a_i \}$, such that (a) $0 < a_1 < a_2 < \cdots$ and (b)
$a_i \to Q \in \sigma(A)^{\prime\prime}$ (weakly) (c) $Q \geq q$ and
(d) $\gamma(Q) < 1$ (actually, as close to zero as you like).  From
(c) and the fact that $\| Q \| = 1$ it follows that $\tau(Q) = 1$ for
all $\tau \in {\mathcal S}$.  Hence, from (b) and (d), we see that
$\tau(a_n) \to 1$ for all $\tau \in {\mathcal S}$ while $\lim
\gamma(a_n) < 1$.

Our final task is to replace the sequence of positive elements $\{ a_i
\}$ above with an increasing sequence of projections which have the
same property.  To do this we need to recall the following K-theoretic
facts about simple, nuclear C$^*$-algebras with real rank zero, stable
rank one and weakly unperforated K-theory (cf.\ \cite[Corollary 6.9.2
and Theorem 6.9.3]{blackadar:book}; we also need Haagerup's result
that quasitraces on nuclear algebras are actually traces
\cite{haagerup-thorbjornsen}): (1) For two projections $p$ and $q$ in
(matrices over) $A$, $p$ is equivalent to a proper subprojection of
$q$ if and only if $\tau(p) < \tau(q)$ for all $\tau \in \TA$; (2) If
$f : \TA \to [0,1]$ is any continuous, affine function and $\epsilon >
0$ is given then there exists a projection $p \in A$ such that $|
f(\tau) - \tau(p) | < \epsilon$ for all $\tau \in \TA$ (this is
essentially a special case of the fact that $K_0 (A)$ is dense in
$Aff(\TA)$ together with the strict ordering fact from (1)).  With
these results in hand it is a simple matter to produce a sequence of
projections $\{ p_n \}$ in $A$ such that $\tau(p_n) \to 1$ for all
$\tau \in {\mathcal S}$ while $\lim \gamma(p_n) < 1$ and, moreover,
such that $0 < [p_1] < [p_2] < \cdots$ (in $K_0(A)$) since each of the
positive elements $a_i$ naturally defines a positive affine function on
$\TA$ and there is a uniform gap between each of these functions
(where we can stick a projection by (2)) since $A$ is simple (hence
every trace is faithful).  The final step is to recall that in a
C$^*$-algebra with stable rank one, two projections are (Murray-von
Neumann) equivalent if and only if they are unitarily equivalent.
Hence we can find unitaries $u_n \in A$ such that $p_1 < u_2 p_2 u_2^*
< u_3 p_3 u_3^* < \cdots$.
\end{proof}

\begin{rem} Note that the equivalence of (1), (2) and (3) above holds
for arbitrary C$^*$-algebras.
\end{rem}

\begin{cor}
Let $A$ be a (separable) simple, nuclear C$^*$-algebra with real rank
zero, stable rank one, weakly unperforated K-theory and such that
$\TA$ is not norm separable.  Then for every countable, weak-$*$ dense
subset ${\mathcal S} \subset \TA$ there exists a sequence of
projections $p_1 \leq p_2 \leq \cdots$ such that $\tau(p_n) \to 1$ for
all $\tau \in {\mathcal S}$ but $\lim \gamma(p_n) < 1$ for some
$\gamma \in \TA$.
\end{cor}

\chapter{Finite representations}

Now that we have spent some time discussing approximation
properties of traces we will turn to another aspect of finite
representation theory.  A typical goal in the representation
theory of groups is to classify all {\em irreducible}
representations of a given group -- i.e.\ factor representations
of type I. Our goal in these notes is quite different, however, as
we wish to study what sort of II$_1$-factors can arise as GNS
representations of a given C$^*$-algebra or class of
C$^*$-algebras.  We believe it is premature to ask for
classification results in the realm of type II
representations\footnote{Indeed, this is probably hopeless except
in trivial cases.} and hence will stick to the more modest goal of
constructing as many different type II$_1$ representations as
possible.

\section{II$_1$-factor representations of some universal C$^*$-algebras}

This section will be concerned with the II$_1$-factor
representation theory of a few examples of universal
C$^*$-algebras.  More precisely, we will consider one of the
infinite (universal) free product C$^*$-algebras
$$\bigast_1^{\infty} C(\mathbb{T}), \bigast_1^{\infty} C([-1,1])
\ \mathrm{or} \bigast_1^{\infty} C([0,1])\footnote{In other words,
the universal C$^*$-algebra generated by a countable number of
unitaries, contractive self-adjoints or, respectively, contractive
positive operators.}$$ and study their II$_1$-factor
representations.  The arguments involved are the same in all three
cases so {\em throughout this section we will let $A$ denote
any one of the C$^*$-algebras above}.\footnote{In fact, all the
results presented here are valid for any C$^*$-algebra $A$ with
the following properties: (1) $A$ is isomorphic to
$\bigast_1^{\infty} A$, (2) $A$ is residually finite dimensional
and (3) every (separable, unital) C$^*$-algebra arises as a
quotient of $A$.}

Since every C$^*$-algebra arises as a quotient of $A$ it is not
hard to see that for every von Neumann algebra $M \subset B(H)$
there is a representation $\pi:A \to M \subset B(H)$ such that
$\pi(A)^{\prime\prime} = M$.  A less obvious fact is that in many
cases $\pi$ can be taken to be {\em injective}.

\begin{prop}
\label{thm:universalfaithful} Let $M$ be a II$_1$-factor.  There
exists a {\em faithful} trace $\tau$ on $A$ such that $\pi_{\tau}
(A)^{\prime\prime} \cong M$.
\end{prop}

\begin{proof}  It suffices to construct a $*$-monomorphism 
$A \hookrightarrow M$ with weakly dense range.  Indeed, uniqueness 
of GNS representations implies that the (unique, faithful) trace on $M$, restricted 
to $A \subset M$, is the $\tau$ we are after.

The construction, though a bit technical, really amounts to some
universal trickery and is not particularly deep.  We first need to
write $A$ as an inductive limit of free products of itself.  That
is, we define $$A_1 = A, A_2 = A_1 * A, \ldots, A_n = A_{n-1} * A,
\ldots,$$ where $*$ denotes the full (i.e.\ universal) free
product (with amalgamation over the scalars).

Letting $B$ denote the inductive limit of the sequence $A_1 \to
A_2 \to \cdots$ it is easy to see (by universal considerations)
that $A \cong B$.\footnote{$B$ is also the universal C$^*$-algebra
generated by countably many operators of the same type as those
that generate $A$.} Since $A$ is residually finite dimensional
(see \cite{davidson} for the case of $\bigast_1^{\infty}
C(\mathbb{T}) \cong C^*(\mathbb{F}_{\infty})$ -- the other two
cases have a similar, but easier, proof) we can find a sequence of
integers $\{ k(n) \}$ and a unital $*$-monomorphism $\sigma : B
\hookrightarrow \Pi M_{k(n)} ({\mathbb C})$.  Note that we may
naturally identify each $A_i$ with a subalgebra of $B$ and hence,
restricting $\sigma$ to this copy of $A_i$, get an injection of
$A_i$ into $\Pi M_{k(n)} ({\mathbb C})$.

To construct the desired embedding of $B$ into $M$, it suffices to
prove the existence of a sequence of unital $*$-homomorphisms
$\rho_i : A_i \to M$ with the following properties:

\begin{enumerate}
\item Each $\rho_i$ is injective.

\item $\rho_{i+1}|_{A_i} = \rho_i$, where we identify $A_i$ with
the `left side' of $A_i * A = A_{i+1}$.

\item The (increasing) union of $\{ \rho_i (A_i) \}$ is weakly
dense in $M$.
\end{enumerate}

To this end, we first choose an increasing sequence of projections
$p_1 \leq p_2 \leq \cdots $ from $M$ such that $\tau_M (p_i) \to
1$. Then define orthogonal projections $q_2 = p_2 - p_1, q_3 = p_3
- p_2, \ldots$ and consider the II$_1$-factors $Q_j = q_j M q_j$
for $j = 2, 3, \ldots$.  As is well known and not hard to
construct, we can, for each $j \geq 2$, find a unital embedding
$\Pi M_{k(n)} ({\mathbb C}) \hookrightarrow Q_j \subset M$ and
thus we get a sequence of (orthogonal) embeddings $B
\hookrightarrow \Pi M_{k(n)} ({\mathbb C}) \hookrightarrow Q_j
\subset M$ which will be denoted by $\sigma_j$.

We are almost ready to construct the $\rho_i$'s. Indeed, for each
$i \in {\mathbb N}$ let $\pi_i : A \to p_i M p_i$ be a (not
necessarily injective!) $*$-homomorphism with weakly dense range.
We then define $\rho_1$ as $$\rho_1 = \pi_1 \oplus \bigg(
\bigoplus_{j \geq 2} \sigma_j|_{A_1} \bigg) : A_1 \hookrightarrow
p_1 M p_1 \oplus \bigg( \Pi_{j \geq 2} Q_j \bigg) \subset M.$$
Note that this is a unital $*$-monomorphism from $A_1$ into $M$
(since each $\sigma_j$ is already faithful on all of $B$).  Now
define a $*$-homomorphism $\theta_2 : A_2 = A_1 * A \to p_2 M p_2$
as the free product of the $*$-homomorphisms $A_1 \to p_2 M p_2$,
$x \mapsto p_2 \rho_1 (x) p_2$, and $\pi_2 : A \to p_2 M p_2$. We
then put $$\rho_2 = \theta_2 \oplus \bigg( \bigoplus_{j \geq 3}
\sigma_j|_{A_2} \bigg) : A_2 \hookrightarrow p_2 M p_2 \oplus
\bigg( \Pi_{j \geq 3} Q_j \bigg) \subset M.$$ Note that
$\rho_2|_{A_1} = \rho_1$.  Hopefully it is now clear how to
proceed.  In general, we construct a map (whose range is dense in
$p_{n+1} M p_{n+1}$) $\theta_{n+1}: A_{n+1} = A_n * A \to p_{n+1}
M p_{n+1}$ as the free product of the cutdown (by $p_{n+1}$) of
$\rho_n$ and $\pi_{n+1}$. This map need not be faithful and hence
we take a direct sum with $\oplus_{j \geq n+2}
\sigma_j|_{A_{n+1}}$ to remedy this deficiency.  It is then easy
to see that these maps have all the required properties and hence
the proof is complete.
\end{proof}

\begin{rem}  While it is not true that every infinite
dimensional von Neumann algebra contains a weakly dense copy of
$A$ (e.g.\ abelian or, more generally, homogeneous von Neumann
algebras of finite type can't contain $A$) it is true that every
infinite dimensional factor contains a weakly dense copy of $A$. A
careful inspection of the proof shows that one really only needs
smaller and smaller corners each of which contains a copy of $\Pi
M_{k(n)} ({\mathbb C})$ for the proof above to go through.
\end{rem}

\section{Elliott's intertwining argument for II$_1$-factors}

This section contains a fairly straightforward adaptation of
Elliott's celebrated approximate intertwining argument which has
been so successful in the classification program.  It has nothing
to do with finite representation theory per se but will be needed
in the next section when we study factor representations of
Popa algebras. Though usually done in the setting of
C$^*$-algebras we will need Elliott's argument in the setting of
von Neumann algebras. While not the most general possible form,
the following version is more than sufficient for our purposes.
The set-up is as follows.

Assume that $M \subset B(L^2 (M,\tau))$ and $N \subset
B(L^2(N,\gamma))$ are von Neumann algebras acting standardly, with
faithful, normal tracial states $\tau$ and, respectively,
$\gamma$. Let $X_1 \subset X_2 \subset \ldots \subset M$ and $Y_1
\subset Y_2 \subset \ldots \subset N$ be (not necessarily unital)
C$^*$-subalgebras such that $\cup X_i$ is weakly dense in $M$ and
$\cup Y_i$ is weakly dense in N.  Further assume that we have c.p.
maps $\alpha_n : Y_n \to X_n$, $\beta_n : X_n \to Y_{n + 1}$,
which are contractive both with respect to the operator norms and
the 2-norms coming from $\tau$ and $\gamma$, and finite subsets
$\Lambda_i \subset X_i$ and $\Omega_i \subset Y_i$ with the
following properties:

\begin{enumerate}
\item $\Lambda_i \subset \Lambda_{i + 1}$, $\Omega_i \subset
\Omega_{i + 1}$, for all $i \in {\mathbb N}$, and the linear spans
of $\cup \Lambda_i$ and $\cup \Omega_i$ are norm dense in $\cup
X_i$ and $\cup Y_i$, respectively, and hence weakly dense in $M$
and $N$, respectively.  To simplify things, we will also assume
that $x_1, x_2 \in \Lambda_i \Longrightarrow x_1x_2 \in
\Lambda_{i+1}$ and, similarly, that $\Omega_{i+1}$ contains the
product of any pair of elements from $\Omega_i$.

\item $\alpha_i (\Omega_i) \subset \Lambda_i$ and $\beta_i
(\Lambda_i) \subset \Omega_{i + 1}$ for all $i \in {\mathbb N}$.

\item Both $\{ \alpha_i \}$ and $\{ \beta_i \}$ are weakly
asymptotically multiplicative. That is, $\| \alpha_i (y_1 y_2) -
\alpha_i (y_1) \alpha_i (y_2) \|_{2,\tau} \to 0$, as $i \to
\infty$, for all $y_1, y_2 \in \cup Y_i$ and similarly for $\{
\beta_i \}$.
\end{enumerate}

\begin{thm}(Elliott's Intertwining)
In the setting described above, if it happens that $\| x -
\alpha_{n + 1} \circ \beta_n (x) \|_{2,\tau} < 1/2^n$ and $\| y -
\beta_n \circ \alpha_n (y) \|_{2,\gamma} < 1/2^n$ for all $x \in
\Lambda_n$, $y \in \Omega_n$ and all $n \in {\mathbb N}$, then $M
\cong N$.
\end{thm}

\begin{proof} For the proof we will need the following two standard
facts: $i)$ every norm bounded sequence which is Cauchy in 2-norm
converges (in 2-norm) and $ii)$ on norm bounded subsets, the
2-norm topology is the same as the strong operator topology (since
our von Neumann algebras are acting standardly on the $L^2$ spaces
coming from their traces).

Since the $\alpha_n$'s and $\beta_n$'s are 2-norm contractive, we
first claim that for each $y \in \cup \Omega_i$, the sequence $\{
\alpha_n (y) \}$ is Cauchy in 2-norm (and similarly for each $x
\in \cup \Lambda_i$, $\{ \beta_n(x) \}$ is Cauchy in 2-norm).  To
see this we first fix $m \in \mathbb{N}$ and note that $$\|
\alpha_m(y) - \alpha_{m+1}(y) \|_{2,\tau} \leq \|\alpha_m(y) -
\alpha_{m+1}\circ\beta_m(\alpha_m(y)) \|_{2,\tau} + \|
\alpha_{m+1}\big( \beta_m\circ\alpha_m(y) - y \big)\|_{2,\tau}$$ 
which, in turn, is bounded above by $\frac{1}{2^{m-1}}.$  Repeated applications of this
inequality shows that for $m \leq n$ $$\| \alpha_m(y) -
\alpha_n(y) \|_{2,\tau} \leq \sum_{i=m-1}^{n-2} \frac{1}{2^i}$$
and hence the sequence is Cauchy as claimed.  Evidently the
sequences $\{ \alpha_n (y) \}$ and $\{ \beta_n(x) \}$ are still
Cauchy when $y$ and $x$ are taken from the linear spans of $\cup
\Omega_i$ and $\cup \Lambda_i$.

Since the $\alpha_n$'s and $\beta_n$'s are norm contractive, it
follows that for each $y$ in the linear span of $\cup \Omega_i$,
the sequence $\{ \alpha_n (y) \}$ is convergent in $M$ (and
similarly for all $x \in span(\cup \Lambda_i$)).  Hence we can
define linear maps $\Phi : span(\cup \Lambda_i) \to N$ and $\Psi :
span(\cup \Omega_i) \to M$ by $\Phi(x) = \lim \beta_n(x)$ and
$\Psi (y) = \lim \alpha_n(y)$.  Note that $\Phi$ and $\Psi$ are
contractive with respect to both the operator norms and the
2-norms.  This implies that $\Phi$ and $\Psi$ can be (uniquely)
extended to the norm closures of $span(\cup \Lambda_i)$ and
$span(\cup \Omega_i)$ (which are weakly dense C$^*$-algebras, by
condition (1) above) and, moreover, that these extensions are
2-norm contractive as well.  Now one uses Kaplansky's density
theorem (and the fact that our extensions are still norm and
2-norm contractive) to extend beyond these weakly dense
C$^*$-subalgebras to all of $M$ and $N$. (i.e.\ For each $x \in M$
we take a norm bounded sequence $\{ x_n \}$ from the norm closure
of $span(\cup \Lambda_i)$ which converges to $x$ in 2-norm. The
image of $\{ x_n \}$ in $N$ is then a norm bounded sequence which
is Cauchy in 2-norm and hence we map $x$ to the (2-norm) limit of
this sequence.)

To save notation, we will also let $\Phi : M \to N$ and $\Psi : N
\to M$ denote the maps constructed in the previous paragraph. Note
that these maps are 2-norm contractive and linear. They are also
$*$-preserving for if $x_n \to x$ (in 2-norm) then it follows that
$x_n^* \to x^*$ as well (i.e.\ the strong and strong-$*$
topologies agree on bounded subsets of a tracial von Neumann
algebra since $\| x \|_{2,\tau} = \| x^* \|_{2,\tau}$).  Since the
maps $\alpha_n$ and $\beta_n$ preserve adjoints by assumption it
follows that the maps $\Phi$ and $\Psi$ preserve adjoints as well.

It is easy to check that $\Phi$ and $\Psi$ are mutual inverses on
the spans of $\cup \Lambda_i$ and $\cup \Omega_i$.  By 2-norm
contractivity it follows that they are mutual inverses on all of
$M$ and $N$. Hence we only have to observe that both $\Phi$ and
$\Psi$ are multiplicative on $M$ and $N$.  Since multiplication is
continuous, on bounded sets, in the 2-norm, it is not hard to
check that $\Phi$ and, respectively, $\Psi$ are multiplicative on
the norm closures of the linear spans of $\cup \Lambda_i$ and,
respectively, $\cup \Omega_i$.  Finally, another application of
Kaplansky's density theorem and a standard interpolation argument
allow one to deduce multiplicativity on all of $M$ and $N$.
\end{proof}

\section{II$_1$-factor representations of Popa Algebras}
\label{sec:Popareps}

The main motivation for Popa's work in \cite{popa:simpleQD} was to
try to understand the relationship between quasidiagonality and
nuclearity. Indeed, in \cite[pg.\ 157]{popa:simpleQD} Popa asked
whether every Popa algebra with unique trace is necessarily
nuclear.  More generally, he asked in \cite[Remark
3.4.2]{popa:simpleQD} whether the hyperfinite II$_1$-factor $R$
was the {\em only} II$_1$-factor which could arise from a GNS
representation of a Popa algebra. However counterexamples to
Popa's first question were constructed by Dadarlat in
\cite{dadarlat:nonnuclearsubalgebras}. Interestingly enough,
though, in \cite{dadarlat:nonnuclearsubalgebras} Dadarlat
constructs nonnuclear tracially AF algebras and hence all of their
II$_1$-factor representations are hyperfinite (even though the
algebras are not nuclear).  Support in favor of Popa's second
question was provided in \cite[Remark 3.4.2]{popa:simpleQD} where
Popa proved that if a factorial trace is in $\TAuQD$ then it gives
the hyperfinite II$_1$-factor (compare with Theorem
\ref{thm:mainthmUTAwafd} -- see also the Introduction).

In this section we will show that Popa's second question also
has a negative answer.  In fact, we will show that there is a
universal Popa algebra for the class of McDuff II$_1$-factors
(meaning there exists a Popa algebra $A$ with the property that
every II$_1$-factor of the form $R\overline{\otimes}M$ appears as
the GNS representation of $A$ with respect to some trace -- see
Theorem \ref{thm:arbitraryMcDuff}). On the other hand, we also
observe that our results on tracial approximation properties yield
an affirmative answer to Popa's question in one non-trivial case;
if $A$ is a locally reflexive Popa algebra with unique trace
$\tau$ then $\pi_{\tau}(A)^{\prime\prime} \cong R$ (see Theorem
\ref{thm:locallyreflexiveuniquetrace}).

We now turn to the main technical result which will allow us to
construct a wide variety of II$_1$-factor representations of Popa
algebras. The informed reader will note that every aspect of this
result can be traced back to the classification program. Indeed,
we will adapt the inductive limit techniques of Dadarlat
\cite{dadarlat:nonnuclearsubalgebras} to construct new Popa
algebras and use Elliott's intertwining argument to understand
their GNS representations.

In the following theorem $\mathfrak{C}$ will denote some
collection of C$^*$-algebras which is closed under increasing
unions (i.e.\ inductive limits with injective connecting maps) and
tensoring with finite dimensional matrix algebras.  A
C$^*$-algebra $E$ is called {\em residually finite dimensional} if
$E$ has a separating family of finite dimensional representations
(i.e.\ for every $0 \neq x \in E$ there exists a $*$-homomorphism
$\pi : A\to M_n ({\mathbb C})$ such that $\pi(x) \neq 0$).

\begin{thm}
\label{thm:basicconstruction} Let $E \in \mathfrak{C}$ be a
residually finite dimensional C$^*$-algebra.  Then there exists a
Popa algebra $A \in \mathfrak{C}$ such that for every $\varepsilon
> 0$ we can find a $*$-monomorphism $\rho : E \hookrightarrow A$
with the property that for each trace $\tau \in \mathrm{T}(E)$,
there exists a trace $\gamma \in \TA$ such that

\begin{enumerate}
\item $|\gamma\circ\rho (x) - \tau(x) | < \varepsilon \| x \|$ for
all $x \in E$ and,

\item $\pi_{\gamma} (A)^{\prime\prime} \cong R \bar{\otimes}
\pi_{\tau} (E)^{\prime\prime}.$
\end{enumerate}
\end{thm}

The proof of this result becomes much more transparent once the
main idea is understood.  Hence we think it is worthwhile to give
the idea first and leave the details to the end.

So suppose that $E$ is a residually finite dimensional
C$^*$-algebra and $\tau \in \mathrm{T}(E)$. Let $\mathcal{U}$ be
some UHF algebra. Then the canonical, unital inclusion $E
\hookrightarrow E\otimes \mathcal{U}$, $x \mapsto x\otimes 1$ is
honestly trace preserving (in fact, yields an isomorphism of
tracial spaces) and the weak closure in any GNS representation is
obviously of the form $R \bar{\otimes} \pi_{\tau}
(E)^{\prime\prime}.$ The problem, of course, is that $E\otimes
\mathcal{U}$ is not a Popa algebra. So the idea is that we will
use an inductive limit construction to get a sequence $$ E \to
E\otimes M_{k(1)} ({\mathbb C}) \to E\otimes M_{k(1)} ({\mathbb
C}) \otimes M_{k(2)} ({\mathbb C}) \to \cdots,$$ such that the
limit is a Popa algebra, but the connecting maps above will be
chosen so that when one applies a trace it will (approximately)
look like the sequence which yields $E \otimes \mathcal{U}$.

We now describe the basic construction which will be needed to get
our Popa algebras.  Let $\pi : E \to M_k ({\mathbb C})$ be a
representation, $\tau \in \mathrm{T}(E)$ and $\epsilon > 0$.
Choose $n \in {\mathbb N}$ very large and consider the map $\rho :
E \to E \otimes M_n({\mathbb C})$ given by
$$x \mapsto 1_E \otimes diag(0_{n-k}, \pi(x)) + x \otimes
diag(1_{n-k}, 0_k),$$ where $diag(0_{n-k}, \pi(x))$ is the block
diagonal element in $M_n({\mathbb C})$ whose first $n - k$ entries
down the diagonal are zero and the bottom block is given by
$\pi(x)$, while $diag(1_{n-k}, 0_k) \in M_n({\mathbb C})$ has $n -
k$ $1$'s down the diagonal followed by $k$ zeros. The key remarks
about this choice of connecting map are:

\begin{enumerate}
\item If $\frac{n - k}{n} > 1 - \epsilon$ then $|\tau \otimes
\mathrm{tr}_n (\rho(x)) - \tau(x) | < 2\epsilon \|x\|$ for all $x
\in E$.  That is, in trace the connecting map $\rho$ is almost the
same as the map $x \mapsto x \otimes 1_{M_n}$ (which would be the
natural connecting maps to use if we were trying to construct $E
\otimes \mathcal{U}$).

\item If $I \subset E$ is an ideal and there exists an element $x
\in I$ such that $\pi(x) \neq 0$ then the ideal generated by
$\rho(I)$ is all of $E\otimes M_n({\mathbb C})$.  This follows
from the definition of $\rho$ and the simplicity of $M_n({\mathbb
C})$.  It is this fact that will allow us to deduce simplicity of
our inductive limits.

\item There exists a finite dimensional C$^*$-algebra $B \subset
E\otimes M_n({\mathbb C})$ with unit $e$ such that $e\rho(x) -
\rho(x)e = 0$ and $e\rho(x)e \in B$, for all $x \in E$.  (Let $B =
diag(0,\ldots,0,\pi(E)) \subset M_n({\mathbb C})$.) This remark
will immediately imply that our inductive limits satisfy the
finite dimensional approximation property which defines Popa
algebras.

\item The representation $\pi \otimes \id : E \otimes M_n({\mathbb
C}) \to M_k \otimes M_n({\mathbb C})$ is again a finite
dimensional representation and hence this whole procedure can be
reapplied to the algebra $E\otimes M_n({\mathbb C})$ (thus
yielding an inductive system).
\end{enumerate}

We now enter the gory details.  So let $E$ be a residually finite
dimensional C$^*$-algebra and $\pi_i : E \to M_{k(i)} ({\mathbb
C})$ be a separating sequence of representations.  In fact, we
will assume that for every $x \in E$, $\| x \| = \lim_i \| \pi_i
(x) \|$ (taking direct sums, it is not hard to see that every
residually finite dimensional C$^*$-algebra has such a sequence).
Note that for every $n \in {\mathbb N}$, $\pi_i \otimes id : E
\otimes M_n({\mathbb C}) \to M_{k(i)} \otimes M_n({\mathbb C})$ is
a separating sequence of the same type.

Now choose natural numbers $1 = n(0) \leq n(1) \leq n(2) \leq
\ldots$ such that $$\frac{n(0)n(1)\cdots n(j-1)k(j)}{n(j)} <
2^{-j},$$ for all $j \in {\mathbb N}$.  One then defines algebras
$E = E_0, E_1 = E_0 \otimes M_{n(1)}, E_2 = E_1 \otimes M_{n(2)},
E_3 = E_2 \otimes M_{n(3)}, \ldots $ and inclusions $\rho_i : E_i
\hookrightarrow E_{i + 1}$ as in the basic construction described
above where the inclusion $\rho_i$ uses the finite dimensional
representation $\pi_{i+1} \otimes id \otimes \cdots \otimes id : E
\otimes M_{n(1)} \otimes \cdots \otimes M_{n(i)} \to M_{k(i+1)}
\otimes M_{n(1)} \otimes \cdots \otimes M_{n(i)} $ in the lower
right hand corner.  Letting $\Phi_{j,i} : E_i \to E_j$, $i \leq
j$, be defined by $\Phi_{j,i} = \rho_{j-1}\circ \cdots \circ
\rho_i$ we get an inductive system $\{ E_i, \Phi_{j,i} \}$.

There are some projections in the above inductive system which we
will need.  Let $P_i \in E_i$ be the projection $$P_i =
1_{E_{i-1}} \otimes diag(1_{n(i) - n(1)\cdots n(i-1)k(i)},
0_{n(1)\cdots n(i-1)k(i)}).$$ Note that $P_{i+1}$ commutes with
all of $\rho_i (E_i)$ (and, in particular, with $\rho_i (P_i)$).
Note also that if we write $E_{i+1} = E_i \otimes M_{n(i+1)}$ then
$$P_{i+1}\rho_i (P_i) = P_i \otimes diag(1_{n(i+1) - n(1)\cdots
n(i)k(i+1)}, 0_{n(1)\cdots n(i)k(i+1)}).$$

Letting $A$ be the inductive limit of the inductive system above,
we only have to show that $A$ is the Popa algebra we are after.

\vspace{2mm}

{\noindent\bf Proof of Theorem \ref{thm:basicconstruction}:} We
keep all the notation above.  We leave it to the reader to verify
that $A$ is a Popa algebra as this follows from our remarks above
and the construction of $A$. (That $A$ is unital and satisfies the
right finite dimensional approximation property is obvious while
simplicity follows from the remark that any ideal in A must
eventually intersect some $E_i$ (cf.\ \cite[Lemma
III.4.1]{davidson}).)  Note also that $A$ was constructed as an
inductive limit of matrices over $E$ and hence belongs to the
class $\mathfrak{C}$ when $E$ does.

Now observe that given a trace $\tau \in T(E)$ we can define
traces $\tau_j \in T(E_j)$ by $\tau_j = \tau\otimes
\mathrm{tr}_{n(1)} \otimes \cdots \otimes \mathrm{tr}_{n(j)}$.
Then the embedding $\rho_j : E_j \to E_{j+1}$ almost intertwines
$\tau_j$ and $\tau_{j+1}$. More precisely, a straightforward (but
rather unpleasant) calculation shows that for $i < j$, $$\tau_j
(\Phi_{j,i} (x)) = \frac{\prod\limits_{s = i}^{j-1} (n(s+1) -
n(1)n(2)\cdots n(s)k(s+1))}{\prod\limits_{s = i}^{j-1} n(s+1)}
\tau_i (x) + \lambda_{i,j} \eta_{i,j}(x),$$ where $\lambda_{i,j} =
1 - \frac{\prod\limits_{s = i}^{j-1} (n(s+1) - n(1)n(2)\cdots
n(s)k(s+1))}{\prod\limits_{s = i}^{j-1} n(s+1)}$ and $\eta_{i,j}$
is some tracial state on $E_i$. Hence we get the estimate
$$|\tau_j (\Phi_{j,i}(x)) - \tau_i (x) | \leq 2\lambda_{i,j} \| x
\|,$$ for all $x \in E_i$.  But, it can be shown by induction that
$\prod\limits_{s=i}^{j-1} (1 - \frac{n(1)n(2)\cdots
n(s)k(s+1)}{n(s+1)}) \geq \prod\limits_{s=i}^{j-1} (1 - 2^{-s-1})
\geq 1 - 2^{-i} + 2^{-j} \geq 1 - 2^{-i},$ for all $i < j \in
{\mathbb N}$. Hence we get that $$|\lambda_{i,j}| = | 1 -
\prod\limits_{s=i}^{j-1} (1 - \frac{n(1)n(2)\cdots
n(s)k(s+1)}{n(s+1)}) | \leq 2^{-i}.$$

We have almost established part (1) in Theorem
\ref{thm:basicconstruction}. For each $i \in {\mathbb N}$, extend
$\tau_i$ to a state on $A$ (after identifying $E_i$ with it's
image in $A$).  It is clear that if we take any weak-$*$ cluster
point, $\gamma$, of this sequence then we will get a trace on $A$.
Moreover, by the estimates above, we have that for each $x \in
E_i$,
$$|\gamma(x) - \tau_i (x)| \leq 2^{-i} \| x \|.$$ Since we always have
$\tau$-preserving embeddings of $E$ into $E_i$, it should be clear
how to construct the embedding $\rho$ in the statement of the
theorem.

Our last task is to prove that $\pi_{\gamma} (A)^{\prime\prime}
\cong \pi_{\tau} (E)^{\prime\prime} \bar{\otimes} R$.  To do this,
we will need to study the projections $P_j \in E_j$ defined above.
The idea is that we will use the $P_j$'s to construct different
projections $Q^{(i)} \in \pi_{\gamma} (A)^{\prime\prime}$ with the
following properties:

\begin{enumerate}
\item $Q^{(i)} = P_i Q^{(i+1)} = Q^{(i+1)} P_i$ and hence $Q^{(i)}
\leq Q^{(i+1)}$ for all $i \in {\mathbb N}$.

\item $Q^{(i+1)} \in \pi_{\gamma} (E_i)^{\prime}$, for all $i \in
{\mathbb N}$ (where we have identified $E_i$ with it's image in
$A$).

\item $\gamma(Q^{(i)}) \geq 1 - 2^{-i}$.

\item For each $i \in {\mathbb N}$,
$\frac{\gamma(Q^{(i+1)}x)}{\gamma(Q^{(i+1)})} = \tau_i (x)$, for
all $x \in E_i$.

\item The natural inclusion of the weak closure of $Q^{(i)}
\pi_{\gamma} (E_{i-1}) Q^{(i)}$ into the weak closure of
$Q^{(i+1)}\pi_{\gamma} (E_{i})Q^{(i+1)}$ (which is a natural
inclusion by (1) above) is isomorphic to the (non-unital)
inclusion $$\pi_{\tau_{i-1}} (E_{i-1})^{\prime\prime}
\hookrightarrow \pi_{\tau_{i}}(E_{i})^{\prime\prime} \cong
\pi_{\tau_{i-1}} (E_{i-1})^{\prime\prime} \otimes M_{n(i)}$$ given
by $$x \mapsto x\otimes diag(1_{n(i) - n(1)n(2)\cdots n(i-1)k(i)},
0_{n(1)n(2)\cdots n(i-1)k(i)}).$$
\end{enumerate}

We claim that the construction of such $Q^{(i)}$'s will complete
the proof.  Indeed, if we can do this then one uses part (5) and
Elliott's approximate intertwining argument to compare the
(non-unital) inclusions $Q^{(1)} \pi_{\gamma} (E_0) \subset
Q^{(2)} \pi_{\gamma} (E_1) \subset \ldots$ to the natural (unital)
inclusions $E_0 \subset E_0 \otimes M_{n(1)} \subset E_0 \otimes
M_{n(1)} \otimes M_{n(2)} \subset \ldots$.  Part (3) ensures that
the former sequence recaptures $\pi_{\gamma} (A)^{\prime\prime}$
while the latter sequence gives $E_0 \otimes \mathcal{U}$ in the
limit, where $\mathcal{U}$ is a UHF algebra, and hence the weak
closure will be as desired and the proof will be complete.

The construction of the $Q^{(i)}$'s is fairly simple.  For each $i
\in {\mathbb N}$ we define projections $Q^{(i)}_n =
\pi_{\gamma}(P_i P_{i+1} \cdots P_{i+n}) \in
\pi_{\gamma}(A)^{\prime\prime}$.  Since the $Q^{(i)}_n$'s are
decreasing (as $n \to \infty$), there exists a strong operator
topology limit.  Define $Q^{(i)} = {\rm sot}-\lim_{n\to \infty}
Q^{(i)}_n$.  Then $Q^{(i)}$ is a projection and it is
straightforward to verify conditions (1) and (2) above. (Recall
that $P_j$ commutes with $\Phi_{j,i} (E_i)$ whenever $i < j$.)
Thus we are left to verify the last three conditions.

Proof of (3). It suffices to show that $\gamma (Q^{(i)}_n) \geq 1
-
   2^{i} - 2^{-i-n}$ for all n.  But we may identify $Q^{(i)}_n$ with
   a projection in $E_{i+n}$ and so using the first part of the proof
   of this theorem (and using the identification) we get $|
   \gamma(Q^{(i)}_n) - \tau_{i+n}(Q^{(i)}_n)| < 2^{-i-n}$.  However,
   it follows from the construction of $Q^{(i)}_n$ that $$\tau_{i+n}
   (Q^{(i)}_n) = \prod\limits_{s = 1}^{n} \mathrm{tr}_{n(i+s)} (diag(1_{n(i+s)
   - n(1)\cdots n(i+s-1)k(i+s)}, 0_{ n(1)\cdots n(i+s-1)k(i+s)})).$$
   Thus, by the calculations given in the first part of the proof of
   this theorem, we see that $\tau_{i+n} (Q^{(i)}_n) \geq 1 - 2^{-i}$
   and hence $\gamma(Q^{(i)}_n) \geq 1 - 2^{i} - 2^{-i-n}$.

Proof of (4). We must show that if $x \in E_{i-1}$ then
$$\tau_{i-1}
    (x)\gamma(Q^{(i)}) = \lim_{n\to \infty} \gamma(Q^{(i)}_n x).$$  In
    order to show this, it suffices to prove that $$\tau_{i-1} (x) =
    \frac{\tau_{i+n}(Q^{(i)}_n x)}{\tau_{i+n}(Q^{(i)}_n)},$$ for all
    $n \in {\mathbb N}$.  This last equality, however, is evident from
    the construction.

Proof of (5). It suffices to show (essentially due to the
uniqueness, up to unitary equivalence, of GNS representations
together with part (4)) that there exists a surjective
$*$-homomorphism $\eta : E_{i-1} \otimes M_{n(i)} \to
Q^{(i+1)}\pi_{\gamma} (E_i)$ such that for every $x \in E_{i-1}$,
$$\eta(x\otimes diag(1_{n(i) - n(1)\cdots n(i-1)k(i)},
0_{n(1)\cdots n(i-1)k(i)})) = Q^{(i)}\pi_{\gamma}(x) =
Q^{(i+1)}P_i \pi_{\gamma}(x).$$ But since
$$P_i \rho_{i-1}(x) = x\otimes diag(1_{n(i) - n(1)\cdots n(i-1)k(i)},
0_{n(1)\cdots n(i-1)k(i)}) \in E_i = E_{i-1} \otimes M_{n(i)},$$
we get the desired homomorphism by identifying $E_{i-1} \otimes
M_{n(i)}$ with it's image in $A$, passing to the GNS construction
$\pi_{\gamma}$ and then cutting down by $Q^{(i+1)}$. \ \ $\square$

With the technical preliminaries now out of the way we can address
Popa's question concerning the set of II$_1$-factors arising from
representations of Popa algebras.  Our previous work on amenable
traces proves relevant in at least one case.

\begin{thm}
\label{thm:locallyreflexiveuniquetrace} Assume $A$ is a locally
reflexive Popa algebra with unique trace $\tau$.  Then
$\pi_{\tau}(A)^{\prime\prime} \cong R$.
\end{thm}

\begin{proof}  Since $A$ is quasidiagonal $\tau$ must be a
quasidiagonal trace.  In particular, it must be an amenable trace
and hence, by Corollary \ref{cor:locallyreflexive}, produce the
hyperfinite II$_1$-factor in the GNS construction.
\end{proof}

Note that we really didn't need $A$ to be a Popa algebra in the
previous result as this was only used in order to deduce
quasidiagonality.  Indeed, if $A$ is any locally reflexive,
quasidiagonal C$^*$-algebra then $A$ must have at least one trace
whose GNS representation gives either a matrix algebra or
$R$\footnote{Use the fact that $\TAim$ is not empty, by
quasidiagonality, and is a face in $\TA$, hence contains an
extreme point of $\TA$.  It follows that $A$ must have at least
one factorial amenable trace and this will do thanks to Corollary
\ref{cor:locallyreflexive}.} and hence in the unique trace case we
are forced to get something hyperfinite.

We will eventually see that the theorem above is false without the
assumption of local reflexivity.  However, we now wish to observe
just how rich the II$_1$-factor representation theory of a Popa
algebra can be in general. Recall that a II$_1$-factor is called
McDuff if it is isomorphic to something of the form
$N\bar{\otimes}R$.  Also recall that such a factor is hyperfinite
if and only if $N \cong R$ (and hence McDuff factors are almost
never hyperfinite).

\begin{thm}
\label{thm:arbitraryMcDuff} There exists a Popa algebra $A$ with
the property that for each McDuff factor, $M$, there exists a
trace $\tau_M \in \TA$ such that $\pi_{\tau} (A)^{\prime\prime}
\cong M$.
\end{thm}

\begin{proof}
Since every II$_1$-factor arises as the weak closure of a GNS
representation of $C^* ({\mathbb F}_{\infty})$, and $C^* ({\mathbb
F}_{\infty})$ is residually finite dimensional, this theorem
follows immediately from Theorem \ref{thm:basicconstruction}.
\end{proof}

One cute consequence of this result is that every McDuff factor
has some approximation properties on a dense subalgebra which are
close to some of the various characterizations of $R$ (though most
McDuff factors are not hyperfinite). Compare with Theorem
\ref{thm:R} to see the analogous statements which characterize
$R$.

\begin{cor}
\label{thm:BFDV}
If $M \subset B(L^2(M))$ is a McDuff factor then
there exists a weakly dense C$^*$-subalgebra $A \subset M$ such
that:

\begin{enumerate}
\item There exist finite rank projections, $P_1 \leq P_2 \leq
\ldots$, such that $\| [P_n, a] \| \to 0$ for all $a \in A$.

\item There exists a state $\phi$ on $B(L^2(M))$ such that $A
\subset B(L^2(M))_{\phi} = \{ T \in B(L^2(M)): \phi(TS) =
\phi(ST), S \in B(L^2(M))\}.$

\item There exists a sequence of II$_1$-factors, $R_n \subset
B(L^2(M))$, such that $R_n \cong R$ for all $n$ and for each $a
\in A$ we can find $x_n \in R_n$ such that $\| a - x_n \| \to 0$.

\item There exists a sequence of normal, u.c.p.\ maps $\varphi_n :
M \to M_{k(n)} ({\mathbb C})$ such that $\| \varphi_n (ab) -
\varphi_n (a) \varphi_n (b) \|_2 \to 0$ for all $a,b \in A$.

\item There exists a completely positively liftable u.c.p. map
$\Phi : M \to R^{\omega}$ such that $\Phi|_A$ is a
$*$-monomorphism.
\end{enumerate}
\end{cor}

\begin{proof}
Let $A$ be the universal Popa algebra from the previous theorem.
With the exception of statement (3), all assertions above follow
from the quasidiagonality of $A$ (identified with a weakly dense
subalgebra of $M$). Indeed, (1) follows from the definition of
quasidiagonality plus Voiculescu's theorem (since a II$_1$-factor
never contains any non-zero compact operators). To see property
(4) we first note that such maps exist on $A$ and hence we can
extend them to u.c.p.\ maps on all of $M$. Of course, the
extensions need not be normal but, as we have seen before, we can
always perturb to normal maps without destroying the asymptotic
multiplicativity on $A$.  Note that (5) follows from (4) via an
argument similar to the one we saw in the proof of (1)
$\Longrightarrow$ (2) from Theorem \ref{thm:amenabletracesII}.
Finally, statement (2) is just asserting that $A$ is contained in
the centralizer of some state on $B(L^2(M))$ which is to say that
$A$ has at least one amenable trace.  However, $A$ has at least
one quasidiagonal trace (actually it has tons of amenable traces)
which is enough to show (2).

Hence we are left to demonstrate (3).  However, this follows from
Voiculescu's theorem since we can also find a faithful
representation of $A$ into $R$ and this representation must be
approximately unitarily equivalent to the inclusion $A \subset M
\subset B(L^2(M))$.
\end{proof}

When we first discovered part (5) in the corollary above, we
thought that it may be useful in showing that every II$_1$-factor
embeds into $R^{\omega}$.  However, if the map $\Phi$ is normal --
which would imply that it is a $*$-homomorphism since it is so on
a weakly dense subalgebra -- then $M \cong R$ and hence there is
no hope of proving normality in general.

Of course, if a Popa algebra (or any other C$^*$-algebra) is
nuclear then all of its representations give hyperfinite von
Neumann algebras.  The next larger class of algebras is the class
of exact C$^*$-algebras and our work shows that as soon as one
leaves the class where hyperfiniteness is automatic one finds
counterexamples to Popa's question concerning the II$_1$-factor
representation theory of Popa algebras.

\begin{thm}
There exists an exact, Popa algebra, $A$, with non-hyperfinite
II$_1$-factor representations.
\end{thm}

\begin{proof}  According to the proof of Corollary
\ref{thm:exactQD} there is an exact residually finite dimensional
C$^*$-algebra $B$ with a trace $\tau$ such that
$\pi_{\tau}(B)^{\prime\prime}$ is a free group factor. Applying
Theorem \ref{thm:basicconstruction} to this example we get an
exact Popa algebra with a trace whose GNS representation yields
$\pi_{\tau}(B)^{\prime\prime} \overline{\otimes}R$ (which is not
hyperfinite).
\end{proof}

Another curious consequence of this work is that McDuff factors
which are generated by exact C$^*$-algebras always have `norm
microstates' on a dense subalgebra.

\begin{thm}
If $M \subset B(L^2(M))$ is McDuff and contains a weakly dense,
exact C$^*$-subalgebra then there exists a weakly dense
C$^*$-subalgebra $A \subset M$ and finite dimensional matrix
subalgebras $M_n \subset B(L^2(M))$ such that for each $a \in A$
there exists a sequence $a_n \in M_n$ such that $\| a - a_n \| \to
0$.  (Hence, for every noncommutative polynomial $P$ in $k$
variables and finite set $\{ a^{(i)} \}_{i = 1}^{k} \subset A$ we
have $\| P(a^{(1)}, \ldots, a^{(k)}) - P(a^{(1)}_n, \ldots,
a^{(k)})_n) \| \to 0$ as $n \to \infty$.)
\end{thm}

\begin{proof}
Since every exact C$^*$-algebra is the quotient of an exact,
residually finite dimensional C$^*$-algebra (cf.\ \cite[Corollary
5.3]{brown:QDsurvey}), it follows from Theorem
\ref{thm:basicconstruction} that $M$ contains a weakly dense,
exact Popa algebra.  In particular, $M$ contains a weakly dense,
exact, quasidiagonal C$^*$-algebra.  Since $M$ is a factor, it
can't contain any nonzero compact operators and hence the result
now follows from \cite{dadarlat:QDapproximation} (see also
\cite{brown:herrero} for the general case).
\end{proof}

Note that the preceding theorem covers many group von Neumann
algebras (e.g.\ $\Gamma = G_1 \times G_2$ where $G_1$ is discrete,
amenable and i.c.c. while $G_2$ is discrete, exact and i.c.c.). We
are not, however, claiming that this result implies $R^{\omega}$
embeddability for such group von Neumann algebras.  Indeed, it is
not at all clear that the existence of norm microstates implies
the existence of `weak' microstates (in the sense of Voiculescu)
since there does not appear to be any way of understanding how the
traces behave on the norm approximations.

\chapter{Applications and connections with other areas}

In the final chapter of these notes we wish to consider a number
of applications and connections with problems which, on the
surface, seem to be quite far removed from finite representation
theory.

\section{Elliott's classification program}
\label{thm:classificationprogram}

The classification program has been an ambitious attempt to find
complete isomorphism invariants for the class of nuclear
C$^*$-algebras.  At this level of generality it is really
impossible to imagine that such an invariant, if it exists, could
ever be `computable' in general.  For example, there is a
one-to-one correspondence between compact metric spaces and
(unital, separable) abelian C$^*$-algebras and hence such a
classification would necessarily include a classification, up to
homeomorphism, of all compact metric spaces.  On the other hand,
there are results in geometry which state, under suitable
`rigidity' assumptions, that manifolds which are only assumed
homotopy equivalent must in fact be diffeomorphic. Also group
theorists had a complete classification of the finite {\em simple}
groups two decades ago but no such result is expected, as far as
we know, for non-simple finite groups.  Our point is that though a
very general classification theorem for nuclear C$^*$-algebras may
exist in principle it is unlikely to be very applicable due to the
inherit difficulty of computing whatever invariant is capable of
completely classifying such a broad class of algebras. However, if
we look to other fields of mathematics for inspiration, then it
does not seem unreasonable to expect that if one further assumes
some sort of `rigidity', in addition to nuclearity, then
classification results with computable invariants might be
possible.  Exactly what `rigid' means in this context is not yet
clear but most operator algebraists would agree that requiring
{\em simplicity}, like in the case of finite groups, is a
reasonable start.

The invariant which has been most successful, and hence attracted
the most attention, is known as the Elliott invariant which we now
describe.\footnote{Many thanks to M. R{\o}rdam and H. Lin for some
helpful discussions regarding various issues discussed in this
section.} For a nuclear C$^*$-algebra $A$, the Elliott
invariant is the triple $(K_0 (A), K_1 (A), \TA)$, where $\TA$ is
the set of tracial states on $A$, together with the natural
pairing $P_A : K_0 (A) \times \TA \to {\mathbb R}$.\footnote{It is
quite possible that the tracial state space be empty and then, of
course, there is no pairing with $K_0(A)$.} Given two algebras $A$
and $B$, we say that their Elliott invariants are isomorphic if
$K_1 (A) \cong K_1(B)$ and there exist a scaled, ordered group
isomorphism $\Phi : K_0 (A) \to K_0 (B)$ and an affine
homeomorphism $T : \TA \to \mathrm{T}(B)$ such that $P_A(x,\tau) =
P_B (\Phi(x),T(\tau))$, for all $(x,\tau) \in K_0 (A) \times \TA$.

For some time it was felt that the Elliott invariant may be a
complete invariant for the class of {\em simple} nuclear
C$^*$-algebras.  However examples of R{\o}rdam definitively showed
that this was not the case (cf.\
\cite{rordam:infinitenotpurelyinfinite}).  But, his examples were
not stably finite and one could still hope that the Elliott
invariant may be complete for simple, stably finite, nuclear
C$^*$-algebras. Unfortunately, some recent work of Toms (cf.\
\cite{toms}) has shown that this is also not possible as he has
constructed examples of non-isomorphic simple ASH algebras with
the same Elliott invariant.

On the other hand, some experts have suggested that the proper
class of `rigid' algebras to expect general classification results
is for the simple, nuclear C$^*$-algebras of real rank zero or
stable rank one (or both).  Since this still contains many of the
`naturally occurring' examples of finite simple nuclear
C$^*$-algebras (e.g.\ irrational rotation algebras, crossed
products of the Cantor set by minimal homeomorphisms, the crossed
product of the CAR algebra by the noncommutative Bernoulli shift
and others) there is still good reason to pursue a general
classification theorem in this setting.

If we restrict to the class of algebras with stable rank one
(adding real rank zero later on) then one can formulate Elliott's
conjecture as follows.

\vspace{2mm}

(Special case of) {\bf Elliott's Conjecture:} If two simple,
nuclear C$^*$-algebras of stable rank one have isomorphic Elliott
invariants (as described above) then they are isomorphic.

\vspace{2mm}

To the untrained eye this special case may look very special
indeed. However, it is still unimaginably general and it seems to
the present author that we are a long way from the resolution of
this `special' case of Elliott's conjecture.  However, it is our
hope that if we add a few more hypotheses, such as real rank zero,
quasidiagonality, weakly unperforated invariant\footnote{See the
section on tracially AF algebras for the definitions of these
things.} and unique trace then classification results may now be
in sight. Moreover, our feeling is that the key to unlocking
Elliott's conjecture (at least in the real rank zero, stable rank
one case) is approximation properties of traces.

As evidence to support this point of view we now (a) present a
number of classification results where uniform locally finite
dimensional traces play an essential role, (b) point out which
tracial approximation question needs to be resolved in order to
complete the general unique trace case mentioned above and (c)
recall some facts which indicate there is reasonable hope of
carrying out part (b). Finally we will end this section with a
number of predictions of Elliott's conjecture as they may be of
independent interest.

So let's have a look at some classification results where uniform
locally finite dimensional traces play an important role.  As
discussed in the introduction to this paper our feeling is that
Huaxin Lin's classification theorem for tracially AF algebras is
the `right' classification theorem (in the real rank zero, stable
rank one, weakly unperforated case) and hence we have no intention
of presenting more general classification theorems.  Rather our
goal will be to show that certain examples are already covered by
Lin's theorem and the way to do this is to use tracial
approximation properties to show that these examples are in fact
tracially AF. We will begin with some `rational' classification
results.  It will be convenient to summarize (special cases of)
some known facts.

\begin{thm}[cf.\ \cite{rordam:tensorUHFI},
\cite{rordam:tensorUHFII},
\cite{blackadar-kumjian-rordam}]\footnote{If one assumes that $A$
is also approximately divisible in the sense of
\cite{blackadar-kumjian-rordam} then one need not tensor with the
UHF algebra in this theorem.} \label{thm:rordamUHF} Assume $A$ is
simple, nuclear and has a unique tracial state. Then for any UHF
algebra $\mathcal{U}$ the (simple, nuclear) tensor product algebra
$A\otimes \mathcal{U}$ has a unique trace, real rank zero, stable
rank one, Blackadar's fundamental comparison
property\footnote{This means that if $p,q \in A$ are projections
and $\tau(p) < \tau(q)$ for all $\tau \in \TA$ then $p$ is
Murray-von Neumann equivalent to a subprojection of $q$. Note that
this implies weak unperforation.} and the Riesz
property.\footnote{In other words, if $x, y_1, y_2$ are positive
elements of $K_0(A)$ and $x \leq y_1 + y_2$ then there exist $z_1,
z_2 \in K_0(A)$ such that $x = z_1 + z_2$ and $z_i \leq y_i$, $i
=1,2$.}
\end{thm}

\begin{thm} If $A$ is simple, nuclear and has a unique trace
$\tau$ which happens to be uniform locally finite dimensional then
for any UHF algebra $\mathcal{U}$ we have that $A\otimes
\mathcal{U}$ is tracially AF.
\end{thm}

\begin{proof}  Since $A\otimes \mathcal{U}$ has real rank zero,
stable rank one and weakly unperforated K-theory this result
follows from our tracial characterization of tracially AF algebras
(cf.\ Proposition \ref{thm:characterization}) together with the
fact that the unique trace on $A\otimes \mathcal{U}$ is also
uniform locally finite dimensional by Proposition
\ref{thm:tensorproduct}.
\end{proof}

\begin{cor}\footnote{Though we haven't seen it written down
anywhere, this corollary may be known to the experts as it can
also be deduced from the results of \cite{lin:ACtraces}, for
example.}  Let $M$ be a compact manifold and $h:M \to M$ a minimal
diffeomorphism.  If $C(M) \rtimes_{h} \mathbb{Z}$ has a unique
trace then for each UHF algebra $\mathcal{U}$ it follows that
$(C(M) \rtimes_{h} \mathbb{Z})\otimes \mathcal{U}$ is isomorphic
to one  of the AH algebras considered in \cite{elliott-gong}.
\end{cor}

\begin{proof}  Letting $A = C(M) \rtimes_{h} \mathbb{Z}$ in the
previous theorem, simplicity and nuclearity are well known and we
have assumed a unique tracial state.  Hence we only have to know
that $A$ satisfies the Universal Coefficient Theorem (UCT) of
Rosenberg-Schochet  -- which it does by results of
Rosenberg-Schochet \cite{rosenberg-schochet} -- and that the
unique trace is uniform locally finite dimensional for then the
previous theorem will imply that $(C(M) \rtimes_{h}
\mathbb{Z})\otimes \mathcal{U}$ is a nuclear, tracially AF algebra
satisfying the UCT and hence it must be isomorphic to one  of the
AH algebras considered in \cite{elliott-gong} by Lin's
classification theorem \cite{lin:TAFclassification}.  However, the
unique trace is locally finite dimensional since the remarkable
work of Q.\ Lin and Phillips (cf.\ \cite{lin-phillips}) implies
that $C(M) \rtimes_{h} \mathbb{Z}$ is an inductive limit of type I
C$^*$-algebras (actually very special `recursive' subhomogeneous
algebras) and hence Corollary \ref{thm:typeI} applies.
\end{proof}

Note that a similar result can be formulated for simple inductive
limits of type I algebras which have unique trace (i.e.\ all of
these become AH algebras after tensoring with a UHF algebra).

Our assumption of a unique trace above was essential to ensure
that the rationalization had real rank zero.  Since a unique trace
does not imply real rank zero in general we will have to add this
assumption (along with stable rank one and weak unperforation
which also need not be automatic) if we wish to drop the
rationalization procedure.  However, once these assumptions are
added then one can relax the unique trace hypothesis to allow a
norm separable tracial space as in Proposition
\ref{thm:characterization} and still obtain classification
results.  But since many people seem to find the unique trace
formulation more aesthetically pleasing we will continue to use
this assumption (though it can be relaxed in all of the
classification results below).

\begin{thm}
\label{thm:ulfdimpliesTAF} Let $A$ be a simple, nuclear
C$^*$-algebra with real rank zero, stable rank one, weakly
unperforated invariant and satisfying the UCT
(\cite{rosenberg-schochet}).  If $A$ has a unique trace $\tau$
which happens to be uniform locally finite dimensional then $A$ is
isomorphic to one of the AH algebras considered in
\cite{elliott-gong}).
\end{thm}

\begin{proof} According to Proposition \ref{thm:characterization}
we don't even need to know nuclearity or the UCT in order to
deduce that $A$ is tracially AF.  However, these two assumptions
are certainly needed to apply Lin's classification theorem
\cite{lin:TAFclassification}.
\end{proof}

\begin{defn}
Given $\tau \in \TA$, we will say that $A$ is {\em
$\tau$-tracially type I} if for each finite subset $\mathfrak{F}
\subset A$ and $\varepsilon > 0$ there exists a type I subalgebra
$B \subset A$ with unit $e$ such that $\| [x,e] \| < \varepsilon$
for all $x \in \mathfrak{F}$, $e\mathfrak{F} e
\subset^{\varepsilon} B$ and $\tau(e)
> 1 - \varepsilon$.
\end{defn}

The definition above is admittedly artificial.  However, it allows
us to treat both inductive limits of type I C$^*$-algebras (even
``locally type I'' algebras) and C$^*$-algebras with finite
tracial topological rank (in the sense of either \cite[Definition
3.1]{lin:tracialtopologicalrank} or \cite[Definition
3.4]{lin:tracialtopologicalrank}) at the same time as both of
these classes are easily seen to be $\tau$-tracially type I with
respect to every $\tau \in \TA$.

\begin{cor}
\label{thm:traciallytypeI} Assume that $A$ is nuclear, has unique
trace $\tau$, is $\tau$-tracially type I, simple, real rank zero,
stable rank one, weakly unperforated K-theory and satisfies the
UCT.  Then $A$ is isomorphic to an AH-algebra.
\end{cor}

\begin{proof} By the previous theorem, we only need to show that
$\tau$ is a uniform locally finite dimensional trace.  The proof
of this is similar to the proof of Proposition
\ref{thm:characterization} but also requires Lemma
\ref{thm:localfinitedimensional} and Lemma \ref{thm:whocares}.

Indeed, according to Lemma \ref{thm:localfinitedimensional} it
suffices to show that for each finite set $\mathfrak{F} \subset A$
and $\epsilon > 0$ there exists a C$^*$-subalgebra $C \subset A$
and a u.c.p.\ map $\phi:C \to M_n(\mathbb{C})$ such that
$d(\mathfrak{F}, C_{\phi}) < \epsilon$ and $\| \tau|_C -
\mathrm{tr}_n\circ\phi \|_{C^*} < \epsilon$.  But, by definition
of $\tau$-tracially type I we can find a type I subalgebra $B
\subset A$ with unit $e$ satisfying the three conditions stated
above.  If we let $C = e^{\perp}Ae^{\perp} \oplus B$ then $C$
nearly contains $\mathfrak{F}$ in norm and if we knew that the
trace $\tau|_C$ was uniform locally finite dimensional then we
would be done.  While that may or may not be true, it is a fact
that $\tau|_C$ is very close in norm to a uniform locally finite
dimensional trace (and this is good enough to complete the proof).
Indeed, since $B$ is a type I algebra the tracial state on $C$
defined by composing the quotient map $C \to B$ with
$\frac{1}{\tau(e)} \tau(\cdot)$ yields a type I von Neumann
algebra in its GNS representation and hence, by Lemma
\ref{thm:whocares}, is uniform locally finite dimensional.  But
this trace is also close in norm to $\tau|_C$ (since $\tau(e) > 1
- \epsilon$) and hence the proof is complete.
\end{proof}

\begin{rem} The point of the previous corollary is that basically
anything which is built out of type I C$^*$-algebras (either via
classical constructions like inductive limits or Lin's more recent
``measurable'' approximations) satisfies the equation $\TA =
\TAulfd$. Hence in the presence of nuclearity, simplicity, real
rank zero, stable rank one, weak unperforation and a natural
separability condition on the size of the tracial space one
immediately gets classification results for such algebras.

We should also mention that if $A$ has real rank zero and finite
decomposition rank in the sense of Kirchberg-Winter
\cite{kirchberg-winter} then it can be shown that $\TA = \TAulfd$.
Hence this class of algebras also fits into the tracial
approximation picture.  However, for the purposes of
classification there is currently little reason to take this point
of view as Wilhelm Winter has shown that the structural
assumptions involved in finite decomposition rank allow one to
prove much nicer classification theorems than we can presently
achieve only from assumptions on traces (cf.\
\cite{winter:classification}).
\end{rem}

The next three results were first proved by Huaxin Lin, but they
are also simple consequences of Corollary
\ref{thm:traciallytypeI}.  Note that in all of these corollaries
the role played by uniform locally finite dimensional traces
appears to be crucial.  Indeed we are not aware of an argument,
other than the tracial approximation approach,  which will provide
the large finite dimensional algebras required to deduce tracially
AF in any of the cases considered below.

\begin{cor}\cite[Theorem 5.16]{lin:ACtraces} Every simple, real rank
zero, stable rank one, weakly unperforated, inductive limit of
type I C$^*$-algebras with unique trace is isomorphic to an
AH-algebra as in \cite{elliott-gong}.
\end{cor}

\begin{cor}\cite[Theorem 7.7]{lin:tracialtopologicalrank} Every
simple, nuclear, C$^*$-algebra with finite tracial topological
rank, satisfying the UCT and having a unique trace is isomorphic
to an AH-algebra as in \cite{elliott-gong}.
\end{cor}

\begin{proof} In \cite{lin:tracialtopologicalrank} Lin shows that
these hypotheses imply real rank zero, stable rank one and weakly
unperforated K-theory.
\end{proof}

\begin{cor}\cite[Corollary 5.5]{lin:ACtraces}
\label{thm:BFDIV} Let $h : M \to M$ be a minimal diffeomorphism of
a compact manifold and $C(X) \rtimes_{h} {\mathbb Z}$ the
corresponding crossed product.  If $C(X) \rtimes_{h} {\mathbb Z}$
has a unique trace $\tau$ and $\tau_* (K_0(C(X) \rtimes_{h}
{\mathbb Z}))$ is dense in $\mathbb{R}$ then it is isomorphic to
an AH-algebra as in \cite{elliott-gong}.
\end{cor}

\begin{proof} By \cite{lin-phillips}, such crossed products are
inductive limits of `recursive subhomogeneous' algebras (which, in
particular, are type I), they have stable rank one and weakly
unperforated K-theory by \cite{phillips:stablerankone} and,
finally, real rank zero by \cite{phillips:realrankzero}.
\end{proof}

\begin{rem}  Note that the first two results above hold true
under the weaker assumption of norm separability of the tracial
space. However, Corollary \ref{thm:BFDIV} does not generalize
quite so easily.  As proved by Phillips in
\cite{phillips:realrankzero}, if one knows that the canonical map
$K_0(C(X) \rtimes_{\phi} {\mathbb Z}) \to Aff(\mathrm{T}(C(X)
\rtimes_{\phi} {\mathbb Z}))$ has dense range then real rank zero
follows and then, together with norm separability of the tracial
space, one deduces that $C(X) \rtimes_{\phi} {\mathbb Z}$ is
isomorphic to an AH algebra.
\end{rem}

We now discuss what needs to be done in order to complete the case
of Elliott's conjecture mentioned at the beginning of this
section as well as point out why there is reason to hope that
this may be possible.

\begin{cor} Let $A$ be a simple, nuclear C$^*$-algebra with
real rank zero, stable rank one, weakly unperforated invariant,
unique trace $\tau$ and satisfying the UCT.  Then the following
are equivalent:
\begin{enumerate}
\item $A$ is isomorphic to one of the AH algebras considered in
\cite{elliott-gong}.

\item $A$ is an inductive limit of type I subalgebras.

\item $A$ is an inductive limit of subalgebras $A_i$ such that
$\pi_{\tau}(A_i)^{\prime\prime}$ is a type I von Neumann algebra
for each $i \in \mathbb{N}$.

\item $\tau$ is uniform locally finite dimensional.
\end{enumerate}
\end{cor}

\begin{proof} Evidently the last statement is the weakest and
Theorem \ref{thm:ulfdimpliesTAF} says that (4) $\Longrightarrow$
(1) so we are done.
\end{proof}

Of course, the corollary above is just a reformulation of things
we have already seen but we are trying to drive home the point
that {\em the real rank zero, stable rank one, weakly
unperforated, UCT, unique trace case of Elliott's conjecture is
equivalent to proving that, under the same assumptions, the unique
trace is always a uniform locally finite dimensional trace.} Now,
we don't mean to give the impression that verifying this tracial
approximation property will be trivial.  However, the problem has
now been isolated and the next two results show that if {\em
either} of the conditions defining uniform locally finite
dimensional traces are relaxed -- and the assumption of
quasidiagonality is added -- then the resulting approximation
property does, in fact, always hold.  In other words, we believe
there is evidence in favor of using uniform locally finite
dimensional traces to complete the case of Elliott's conjecture
described above.

\begin{thm}[Special case of Theorem \ref{thm:locallyreflexive}]
Let $A$ be a nuclear, quasidiagonal C$^*$-algebra with unique
trace $\tau$. Then $\tau$ is a uniform quasidiagonal trace.
\end{thm}

\begin{thm}
\label{thm:what} Let $A$ be a simple, quasidiagonal C$^*$-algebra
with unique trace $\tau$.  Then $\tau$ is  locally finite
dimensional.
\end{thm}

\begin{proof} Assume $A \subset B(H)$ is a faithful, irreducible
representation (which exists by simplicity).  Since $A$ is
quasidiagonal we can find finite rank projections $P_1 \leq P_2
\leq \cdots$ converging strongly to the identity and such that $\|
[ a, P_n ] \| \to 0$ for all $a \in A$ (note that $A$ is
necessarily in general position -- i.e.\ contains no compacts --
since $A$ is simple and unital). Define u.c.p.\ maps $\phi_n : A
\to P_nB(H)P_n$ by $\phi_n(a) = P_n a P_n$ and note that, by
uniqueness of $\tau$, $\mathrm{tr}_{rank(P_n)} \circ \phi_n (a)
\to \tau(a)$ for all $a \in A$.  Since the representation is
irreducible, it can be shown (cf.\ \cite[Corollary
3.5]{blackadar-kirchberg:II}) that for every $a \in A$, the
distance from $a$ to the multiplicative domain of $\phi_n$ tends
to zero and hence $\tau$ is locally finite dimensional.
\end{proof}

\begin{cor}  If $A$ is simple, nuclear, quasidiagonal and has
unique trace $\tau$ then $\tau \in \TAuQD \cap \TAlfd$.
\end{cor}

Note that Theorem \ref{thm:what} has absolutely nothing to do with
nuclearity.  Thus a general strategy naturally presents itself;
take the maps provided by Theorem \ref{thm:what} and use
nuclearity (or, perhaps, just local reflexivity as in the proof of
Theorem \ref{thm:locallyreflexive}?) to average them in a way
which preserves multiplicative domains while forcing norm
convergence to the trace $\tau$. We remark, however, that
nuclearity (or local reflexivity) {\em must} play a significant
role in the proof as we will show in the next section that
there exists a (non-locally reflexive) Popa algebra with unique
trace which is not {\em uniform} locally finite dimensional
(though it is, by Theorem \ref{thm:what}, locally finite
dimensional).

Finally, we should mention that Popa's work together with Theorem
\ref{thm:locallyreflexive} provides another approach to attacking
this problem.   Indeed, we have the following result which is
obviously relevant to this discussion as it shows that one can
always find a finite family of finite dimensional subalgebras (as
opposed to a single finite dimensional subalgebra as is required
for tracially AF) such that a linear combination of their units is
large in trace.

\begin{thm}
\label{thm:idontknow} Let $A$ be a locally reflexive Popa algebra
with real rank zero and unique trace $\tau$. Then for every finite
set ${\mathcal F} \subset A$ and $\varepsilon > 0$ there exists $n
\in {\mathbb N}$, subalgebras $Q_1, \ldots, Q_m \subset A$ each of
which is isomorphic to $M_n ({\mathbb C})$ and with units $e_1,
\ldots, e_m$ such that $\| [e_i, x] \| < \varepsilon$ for $1 \leq
i \leq m$ and all $x \in {\mathcal F}$, $e_i {\mathcal F} e_i
\subset^{\varepsilon} Q_i$ for $1 \leq i \leq m$ and, finally, $$
\| \frac{1}{n} \sum_{k = 1}^{m} e_k - 1_A \|_{\tau,1} <
\varepsilon,$$ where $\| x \|_{\tau,1} = \tau(|x|)$.
\end{thm}

\begin{proof}
By \cite[Theorem 3.3]{popa:simpleQD} it suffices to show that
$\tau \in \TAuQD$.  However, as above this follows from Theorem
\ref{thm:locallyreflexive} and hence we are done.
\end{proof}

We now want to spend some time looking at predictions of Elliott's
conjecture  -- i.e.\ necessary conditions.  It is quite
interesting that various statements about quasidiagonality and
approximation of traces naturally appear as necessary conditions
for Elliott's conjecture to hold.  We begin with a few simple
facts which are known to the classification experts.

\begin{prop}
\label{thm:elliottpredictsASH} Elliott's conjecture predicts that
if $A$ is a stably finite, simple, nuclear C$^*$-algebra and
$\mathcal{M}_{2^{\infty}}$ denotes the CAR algebra then $A \otimes
\mathcal{M}_{2^{\infty}}$ is an inductive limit of subhomogeneous
algebras (i.e.\ ASH).
\end{prop}

\begin{proof}  To deduce stable rank one and Blackadar's comparison
property for the algebra $A \otimes \mathcal{M}_{2^{\infty}}$ we
really only need to assume stable finiteness of $A$ in Theorem
\ref{thm:rordamUHF}.  Hence $A \otimes \mathcal{M}_{2^{\infty}}$
is a simple, nuclear C$^*$-algebra with stable rank one and weakly
unperforated K-theory.  But Elliott has shown that every possible
weakly unperforated Elliott invariant arises from an ASH algebra
(see, for example, the appendix of \cite{elliott-villadsen}) and
hence we can find a simple ASH algebra $B$ whose Elliott invariant
is isomorphic to that of $A\otimes \mathcal{M}_{2^{\infty}}$. If
$B$ does not have stable rank one then we can replace $B$ with $B
\otimes \mathcal{M}_{2^{\infty}}$ to get a simple ASH algebra with
stable rank one and Elliott invariant isomorphic to that of
$A\otimes \mathcal{M}_{2^{\infty}}$ (since $A\otimes
\mathcal{M}_{2^{\infty}} \otimes \mathcal{M}_{2^{\infty}} \cong
A\otimes \mathcal{M}_{2^{\infty}}$).  Thus our particular
formulation of Elliott's conjecture indeed predicts that $A
\otimes \mathcal{M}_{2^{\infty}}$ is ASH.
\end{proof}

\begin{cor} Elliott's conjecture predicts that every
simple, stably finite, nuclear C$^*$-algebra is quasidiagonal.
\end{cor}

\begin{proof} $A \otimes \mathcal{M}_{2^{\infty}}$ is quasidiagonal
since it is an inductive limit of subhomogeneous algebras and
quasidiagonality passes to subalgebras (i.e.\ $A \cong A\otimes 1
\subset A \otimes \mathcal{M}_{2^{\infty}}$).
\end{proof}

\begin{prop}
\label{thm:elliottpredictsTAF} Elliott's conjecture predicts that
if $A$ is simple, nuclear, stable rank one, real rank zero and has
weakly unperforated invariant then $A$ is tracially AF.
\end{prop}

\begin{proof} We first remark that in the real rank zero case the
tracial simplex is no longer relevant and hence Elliott's
invariant reduces to K-theory alone.  Indeed, if both $A$ and $B$
are C$^*$-algebras of real rank zero and $\Phi : K_0 (A) \to K_0
(B)$ is a scaled, ordered group isomorphism such that $\Phi([1_A])
= [1_B]$ then $\Phi$ induces an affine homeomorphism $\TA \to
\mathrm{T}(B)$ since we may (affinely, homeomorphically) identify
$\TA$ (resp.\ $\mathrm{T}(B)$) with the states in ${\rm
Hom}(K_0(A), {\mathbb R})$ (resp.\ ${\rm Hom}(K_0(B), {\mathbb
R})$). To see that this is true, we first note that the obvious
map $\TA \to {\rm Hom}(K_0(A), {\mathbb R})$ is affine and
injective since $A$ has real rank zero.  It is also onto the
states in ${\rm Hom}(K_0(A), {\mathbb R})$ since every state on
$K_0 (A)$ comes from a trace on $A$ when $A$ is unital and exact
(cf.\ \cite{haagerup-thorbjornsen}).  Finally, it is easy to check
(again using real rank zero) that a sequence of traces $\tau_n \in
\TA$ converges to $\tau \in \TA$ in the weak-$*$ topology if and
only if their images in ${\rm Hom}(K_0(A), {\mathbb R})$ converge
in the topology of pointwise convergence and hence our
identification is also a homeomorphism.

In \cite{elliott-gong} it is shown how to construct simple AH
algebras with real rank zero and with arbitrary unperforated
K-theory and Riesz interpolation property.  Our assumptions on $A$
also imply the Riesz property (cf.\ \cite{eilers-elliott}) and
hence Elliott's conjecture predicts that $A$ is isomorphic to one
of the AH algebras constructed in \cite{elliott-gong}.  However,
as observed by Lin \cite[Proposition 2.6]{lin:TAF}, the
Elliott-Gong construction always yields tracially AF algebras and
the proof is complete.
\end{proof}

Our next goal is to show that Elliott's conjecture predicts that
every trace on every (not necessarily simple) nuclear,
quasidiagonal C$^*$-algebra is uniform quasidiagonal. (Recall that
the corresponding statement for exact, quasidiagonal
C$^*$-algebras is not true -- see Corollary \ref{thm:exactQD}.)
The passage from the simple to non-simple cases is provided by the
next lemma.

\begin{lem}
\label{thm:lemmasection4} If $\mathfrak{C}$ is a collection of
C$^*$-algebras which contains $\mathbb C$ and is closed under $i)$
increasing unions (i.e.\ inductive limits with injective
connecting maps), $ii)$ quasidiagonal, semi-split extensions
(i.e.\ if $0 \to I \to E \to B \to 0$ is a semi-split (cf.\
\cite{blackadar:book}), short exact sequence, $I$ contains an
approximate unit {\em of projections} which is quasicentral in $E$
and both $I, B \in \mathfrak{C}$ then $E \in \mathfrak{C}$) and
$iii)$ tensoring with finite dimensional matrix algebras then the
following are equivalent:

\begin{enumerate}
\item $\TA = \TAQD$ (resp.\ $\TA = \TAuQD$) for every
quasidiagonal $A \in \mathfrak{C}$.

\item $\TA = \TAQD$ (resp.\ $\TA = \TAuQD$) for every residually
finite dimensional $A \in \mathfrak{C}$.

\item $\TA = \TAQD$ (resp.\ $\TA = \TAuQD$) for every Popa algebra
$A \in \mathfrak{C}$.
\end{enumerate}

If, moreover, the class $\mathfrak{C}$ is closed under tensor
products with (non-unital) abelian algebras (it actually suffices
to know $A \in \mathfrak{C} \Longrightarrow A\otimes C_0 ((0,1])
\in \mathfrak{C}$) then the following are equivalent:

\begin{enumerate}
\item[(4)] $\TA = \TAim$ (resp.\ $\TA = \TAuim$) for every $A \in
\mathfrak{C}$.

\item[(5)] $\TA = \TAim$ (resp.\ $\TA = \TAuim$) for every
quasidiagonal $A \in \mathfrak{C}$.

\item[(6)] $\TA = \TAim$ (resp.\ $\TA = \TAuim$) for every
residually finite dimensional $A \in \mathfrak{C}$.

\item[(7)] $\TA = \TAim$ (resp.\ $\TA = \TAuim$) for every Popa
algebra $A \in \mathfrak{C}$.
\end{enumerate}
\end{lem}

\begin{proof} We first prove the equivalence of (1) - (3)
 and then indicate the changes necessary to prove the second part.
 The proofs are the same whether dealing with $\TAQD$ or $\TAuQD$ and
 hence we just treat the uniform quasidiagonal
 case.

(1) $\Longrightarrow$ (3) is immediate.  (3) $\Longrightarrow$ (2)
    follows from Theorem \ref{thm:basicconstruction}.

(2) $\Longrightarrow$ (1).  Let $A \in \mathfrak{C}$, $\tau \in
\TA$ and ${\mathfrak{F}} \subset A$ be an arbitrary finite set.
Since $A$ is quasidiagonal we can find a sequence of u.c.p.\ maps
$\varphi_{n} : A \to M_{k(n)} ({\mathbb C})$ which are
asymptotically multiplicative and asymptotically isometric (i.e.\
$\| a \| = \lim \| \varphi_{n} (a) \|$, for all $a \in A$).
Passing to a subsequence, if necessary, we may assume that
$\varphi_{1}$ (and all the other $\varphi_n$'s) is as close to
multiplicative on $\mathfrak{F}$ as we like.  Put $\Phi = \oplus_n
\varphi_{n} : A \to \Pi_n M_{k(n)} ({\mathbb C})$, and let $E$ be
the C$^*$-algebra generated by $\Phi(A)$.  Note that $\Phi : A \to
E$ is as close to multiplicative on $\mathfrak{F}$ as we like, by
construction.

Now observe that we have a semi-split, quasidiagonal, short exact
sequence: $$0 \to \oplus M_{k(n)} ({\mathbb C}) \to E + \oplus
M_{k(n)} ({\mathbb C}) \to A \to 0.$$ Since $\mathfrak{C}$ is
closed under all of the operations used, it follows that $E +
\oplus M_{k(n)} ({\mathbb C}) \in {\mathfrak{C}}$ and it is clear
that $ E + \oplus M_{k(n)} ({\mathbb C})$ is residually finite
dimensional.  Hence every trace on $E + \oplus M_{k(n)} ({\mathbb
C})$ is  uniform quasidiagonal by (2).  In particular, the trace
$\tau \in \TA$ induces such a trace on $E + \oplus M_{k(n)}
({\mathbb C})$ by composing with the quotient map. Since the
splitting $\Phi : A \to E \subset E + \oplus M_{k(n)} ({\mathbb
C})$ is almost multiplicative, it follows that we can construct a
u.c.p.\ map on $A$ (by composing maps on $E$ with $\Phi$) which is
almost multiplicative on $\mathfrak{F}$ and which approximately
recaptures the trace $\tau$. This completes the proof of (2)
$\Longrightarrow$ (1).

For the equivalence of (4) - (7), we really only need to show the
implication (5) $\Longrightarrow$ (4) as the arguments above go
through without change for the other implications.  To prove this
we will need Proposition \ref{thm:passagetoquotient} which says
that amenable traces pass to quotients whenever there is a c.p.\
splitting for the quotient map and that uniform amenable traces
always pass to quotients.

With those results in hand, (5) $\Longrightarrow$ (4) becomes very
simple. Indeed, let $B$ be the unitization of the cone over $A$
(i.e.\ the unitization of $C_0 ((0,1]) \otimes A$).  Then $B$ is
quasidiagonal (cf.\ \cite{dvv:QDhomotopy}) and belongs to
$\mathfrak{C}$.  Moreover, there is a natural surjective
$*$-homomorphism $B \to A \oplus {\mathbb C} \to A$.  A
(non-unital) c.p. splitting for this quotient map is given by $a
\mapsto e \otimes a$, where $e \in C_0 ((0,1])$ is any
non-negative function such that $e(1) = 1$.  Hence if every trace
on $B$ is (uniform) amenable then every trace on $A$ enjoys the
same property.
\end{proof}

The assumptions on the class $\mathfrak{C}$ may seem unusual, but
note that any one of the following classes of C$^*$-algebras is
closed under the operations needed in the lemma above: nuclear
C$^*$-algebras, exact C$^*$-algebras (cf.\ \cite[Section
7]{kirchberg:commutantsofunitaries}), real rank zero
C$^*$-algebras (cf.\ \cite[2.10, 3.1, 3.14]{brown-pedersen}, it is
easy to prove that if an extension is semi-split and quasidiagonal
then every projection in the quotient lifts to a projection in the
middle algebra) - though these algebras are not closed under
tensoring with $C_0 ((0,1])$.

\begin{prop}
\label{thm:elliott} If Elliott's conjecture holds for all nuclear,
simple, quasidiagonal C$^*$-algebras with stable rank one and
weakly unperforated K-theory then $\TA = \TAuQD$ for every
nuclear, quasidiagonal C$^*$-algebra $A$.
\end{prop}

\begin{proof}
We will apply the previous lemma to the set $\mathfrak{C}$ of all
 nuclear C$^*$-algebras. We remark that extensions of nuclear
 C$^*$-algebras are again nuclear by \cite[Corollary
 3.3]{choi-effros:nuclearityandinjectivity}.

So assume that Elliott's conjecture holds for all nuclear, simple,
quasidiagonal C$^*$-algebras with stable rank one and weakly
unperforated K-theory and let $A$ be quasidiagonal and nuclear. By
the previous lemma, we may assume that $A$ is simple. Let
$\mathcal{M}_{2^{\infty}}$ be the CAR algebra and note that every
trace on $A$ extends (in fact, uniquely) to a trace on $A \otimes
\mathcal{M}_{2^{\infty}}$. Hence it suffices to show that
$\mathrm{T}(A\otimes \mathcal{M}_{2^{\infty}}) =
\mathrm{UAT}(A\otimes \mathcal{M}_{2^{\infty}})_{{\rm LFD}}$. (It
is not clear that this will imply $\TA = \TAulfd$, but it will
certainly imply that $\TA = \TAuQD$.)

But just as in the proof of Proposition
\ref{thm:elliottpredictsASH}, it follows that $A\otimes B$ is an
ASH algebra. In particular, it is an inductive limit of type I
algebras and so from Corollary \ref{thm:typeI} we deduce that
$\mathrm{T}(A\otimes B) = \mathrm{UAT}(A\otimes B)_{{\rm LFD}}$.
\end{proof}

Note that if every stably finite nuclear C$^*$-algebra turns out
to be quasidiagonal (recall that this question was asked by
Blackadar and Kirchberg \cite{blackadar-kirchberg}) then Elliott's
conjecture would predict that every trace on every nuclear
C$^*$-algebra is uniform quasidiagonal.  Whether or not this
happens in the non-simple case, we still have the following
corollary.

\begin{cor}
If Elliott's conjecture holds for all nuclear, simple,
quasidiagonal C$^*$-algebras with stable rank one and unperforated
K-theory then $\TA = \TAuQD$ for every simple, nuclear
C$^*$-algebra $A$.
\end{cor}

\begin{proof}
If $A$ is simple and not stably finite then $A$ has no tracial
states. On the other hand, if $A$ is simple and stably finite then
Elliott's conjecture predicts that $A$ is also quasidiagonal.
\end{proof}

Note that the corollary above begs the following question: If $A$
is nuclear and $\tau \in \TA$ does there exist a {\em simple},
nuclear C$^*$-algebra $B$ with trace $\gamma \in \TB$ and a
$*$-homomorphism $\pi : A \to B$ such that $\tau =
\gamma\circ\pi$?  If this is the case then Elliott's conjecture
predicts that every trace on every nuclear C$^*$-algebra is
uniform quasidiagonal.

We hope that the reader is now convinced that tracial
approximation properties arise naturally as both necessary
conditions (cf.\ Proposition \ref{thm:elliott}) and sufficient
conditions (cf.\ Theorem \ref{thm:ulfdimpliesTAF}) for various
cases of Elliott's conjecture.  We wish to end this section
with a possible approach to the general real rank zero, stable
rank one, weakly unperforated case of Elliott's conjecture. Of
course our focus is on tracial approximation so the following
strategy is by no means the unique path to resolving this case of
the conjecture.

Consider the following problems for a simple, nuclear,
C$^*$-algebra $A$ with real rank zero, stable rank one, weakly
unperforated K-theory and satisfying the UCT.
\begin{enumerate}
\item Is it true that $\TAQD \neq \emptyset$? (Recall that $\TA =
\TAuim \neq \emptyset$ by nuclearity and stable finiteness.)

\item If $A$ is quasidiagonal then do we have $\TAim = \TAQD =
\TAlfd$? (Recall that this is true in the unique trace case.)

\item Is it true that $\TAuQD \cap \TAlfd = \TAulfd$?

\item Is it possible to remove the hypothesis of $\|\cdot
\|_{A^*}$-separability from Proposition \ref{thm:characterization}
in this setting?
\end{enumerate}

\begin{prop}
\label{thm:BFD}
The real rank zero, stable rank one, weakly
unperforated, UCT case of Elliott's conjecture is equivalent to
affirmative answers to all four of the problems posed above.
\end{prop}

\begin{proof}  Assume first that Elliott's conjecture holds in the
case described above.  It follows that if $A$ is any C$^*$-algebra
with all of those properties then it must be tracially AF (cf.\
Proposition \ref{thm:elliottpredictsTAF}).  Since $\TA = \TAulfd$
for every tracially AF algebra it would follow that all four
questions above must have affirmative answers.

Now suppose that one could prove all four questions above and we
will show that it would follow that every $A$ as above must
necessarily be tracially AF.  By problem (1) we would have that
$A$ is necessarily quasidiagonal since simplicity implies that any
trace is faithful and we already saw in the proof of Proposition
\ref{thm:rosenberg} that the existence of a faithful quasidiagonal
trace implies that the C$^*$-algebra is quasidiagonal. If we knew
that (2) also held then we would have that $$\TA = \TAuQD =
\TAlfd$$ since nuclearity (trivially) implies $\TA = \TAim$ and
(non-trivially) $\TAQD = \TAuQD$.  In particular, problem
(3) would then tell us that $$\TA = \TAuQD \cap \TAlfd = \TAulfd$$
and hence one would only have to remove the norm separability
hypothesis from Proposition \ref{thm:characterization} (i.e.\
prove problem (4)) in order to deduce that $A$ was tracially AF.
Of course, Lin's classification theorem would then imply this case
of Elliott's conjecture.
\end{proof}

\section{Counterexamples to questions of Lin and Popa}

Now that Huaxin Lin has essentially completed the tracially AF
case of Elliott's conjecture (cf.\ \cite{lin:TAFclassification})
it is clear that having abstract hypotheses which imply that a
given C$^*$-algebra is tracially AF is of fundamental importance.
One result along these lines is Proposition
\ref{thm:characterization}. On page 694 of \cite{lin:TAF} Lin
wrote ``It is certainly tempting to conjecture that every
quasidiagonal simple C$^*$-algebra of real rank zero, stable rank
one and with weakly unperforated $K_0$ is TAF (tracially AF).''
Sorin Popa has asked (private communication) whether every Popa
algebra with real rank zero and unique trace is tracially AF.  The
purpose of this section is to resolve these questions
negatively. Essentially we show that no amount of
K-theoretic assumptions alone will allow one to deduce that a Popa
algebra is tracially AF.  C$^*$-algebraic structure (e.g.\
nuclearity or, perhaps, local reflexivity) must also be assumed.
As we have tried to argue in the previous section, we believe
that approximation properties of traces provide a reasonable
strategy for proving that certain C$^*$-algebras are tracially AF.
In this section we will see that our study of finite
representation theory also shows why certain algebras can't be
tracially AF.

Before constructing our counterexamples we need a few remarks
concerning the Universal Coefficient Theorem (UCT) of
Rosenberg-Schochet (cf.\ \cite{rosenberg-schochet}).  A very nice
discussion of the UCT can be found in
\cite{phillips:nonclassification}. Indeed, in this manuscript
Phillips not only describes what it means for a C$^*$-algebra to
satisfy the UCT but he also proves the following two results which
we will need.

\begin{prop}[Two out of Three Principle]
\label{thm:twooutofthreeprinciple} If $0 \to I \to E \to B \to 0$
is a semi-split extension (i.e.\ there exists a contractive c.p.\
splitting $B \to E$) and any two of $I$, $E$ and $B$ satisfy the
UCT then so does the third.
\end{prop}

\begin{prop}
\label{thm:UCTinductivelimits} Let $\phi_n : A_n \to A_{n+1}$ be
injective $*$-homomorphisms and let $A$ denote the inductive limit
of this sequence.  Assume further that for each $n$, there exists
a contractive c.p.\ map $\psi_n : A_{n+1} \to A_n$ such that
$\psi_n \circ \phi_n = id_{A_n}$ for all $n$.  If each $A_n$
satisfies the UCT then so does $A$.
\end{prop}

\begin{rem}
\label{thm:UCTremark} Note that the inductive limit construction
used in Theorem \ref{thm:basicconstruction} does have the
one-sided c.p.\ inverses required in Proposition
\ref{thm:UCTinductivelimits} and hence the Popa algebra
constructed in that theorem will satisfy the UCT whenever the
original residually finite dimensional algebra does.  To see that
such c.p.\ inverses exist, we let $\pi: E \to M_k ({\mathbb C})$
be a $*$-homomorphism and $\rho : E \to E\otimes M_n ({\mathbb
C})$ be as in the basic construction.  The desired map $E\otimes
M_n ({\mathbb C}) \to E \cong E \otimes e_{1,1}$ is just given by
compressing to the (1,1) corner.
\end{rem}

We now answer Lin's question negatively.

\begin{thm}
\label{thm:exoticpopaalgebra} There exists an exact, Popa algebra,
$A$, with real rank zero, stable rank one, UCT, Blackadar's
fundamental comparison property (i.e.\ if $p,q \in A$ are
projections such that $\tau(q) < \tau(p)$ for all $\tau \in \TA$
then $q$ is equivalent to a subprojection of $p$), unperforated
K-theory, Riesz interpolation property and which is approximately
divisible and an increasing union of residually finite dimensional
subalgebras such that $\TA \neq \TAim$.
\end{thm}

\begin{proof} We first claim that it suffices to find a
C$^*$-algebra $C$ which is residually finite dimensional, exact,
real rank zero, satisfies the UCT and such that $\mathrm{T}(C)
\neq \mathrm{AT}(C)$.  Indeed, if we can find such a $C$ then by
applying Theorem \ref{thm:basicconstruction} to the class of all
exact, real rank zero C$^*$-algebras which satisfy the UCT we can
find an exact Popa algebra with real rank zero, UCT and such that
$\TA \neq \TAim$.  (See Remark \ref{thm:UCTremark} for the UCT
assertion and \cite{brown-pedersen} for a proof that matrices over
a real rank zero algebra also have real rank zero.) Then replacing
$A$ with $A\otimes \mathcal{U}$, where $\mathcal{U}$ is some UHF
algebra, will be the desired example since this operation
preserves Popa's property, exactness, real rank zero, UCT and
picks up stable rank one, Riesz interpolation (cf.\
\cite[Corollary 3.15]{blackadar-kumjian-rordam}), Blackadar's
fundamental comparison property and hence unperforated K-theory
(cf.\ \cite{rordam:tensorUHFI}, \cite{rordam:tensorUHFII}).
Moreover, it is clear that this example will be an inductive limit
of residually finite dimensional subalgebras and be approximately
divisible in the sense of \cite{blackadar-kumjian-rordam}.

The construction of the desired residually finite dimensional
C$^*$-algebra requires another one of Kirchberg's
characterizations of exactness \cite[Theorem
1.3]{kirchberg:commutantsofunitaries}: A separable C$^*$-algebra
$A$ is exact if and only if there exists a subalgebra $B$ of the
CAR algebra, $M_{2^{\infty}}$, and an AF ideal $J \subset B$ such
that $A \cong B/J$.  We remark if $A$ is exact and $0 \to J \to B
\to A \to 0$ is the short exact sequence given by Kirchberg's
theorem, then this sequence is automatically semi-split (i.e.\
there exists a c.p. splitting $A \to B$) since $B$ is exact and
$J$ is nuclear (cf.\ the bottom of page 41 in
\cite{kirchberg:commutantsofunitaries}).

Since $C^*_r ({\mathbb F}_2)$ has a unique trace, it follows from
\cite[Theorem 7.2]{rordam:tensorUHFII} that $C^*_r ({\mathbb
F}_2)\otimes M_{2^{\infty}}$ is exact and has real rank zero.  By
Kirchberg's characterization, we can find an exact, quasidiagonal
C$^*$-algebra $B$ with an AF ideal $J \subset B$ such that $B/J
\cong C^*_r ({\mathbb F}_2)\otimes M_{2^{\infty}}$.  Moreover, the
short exact sequence $0 \to J \to B \to C^*_r ({\mathbb
F}_2)\otimes M_{2^{\infty}} \to 0$ is semi-split.

Since $J$ is AF, it follows from \cite[Theorem 3.14 and Corollary
3.16]{brown-pedersen} that $B$ also has real rank zero.  Note that
$B$ has a trace which gives $L(\mathbb{F}_2)\overline{\otimes} R$
in its GNS representation and, by Corollary
\ref{cor:locallyreflexive}, this trace is not amenable.  Hence
$\mathrm{T}(B) \neq \mathrm{AT}(B)$. Finally, since $C^*_r
({\mathbb F}_2)$ satisfies the UCT, it follows from the `two out
of three principle' (cf.\ Proposition
\ref{thm:twooutofthreeprinciple}) that $B$ also satisfies the UCT.
In other words, $B$ is almost the desired algebra; we only have to
replace $B$ with something residually finite dimensional ($B$ is
only quasidiagonal).

From the proof of (2) $\Longrightarrow$ (1) in Lemma
\ref{thm:lemmasection4} we can use $B$ to construct a residually
finite dimensional C$^*$-algebra $C$ such that $C$ is exact, real
rank zero, satisfies the UCT and such that $\mathrm{T}(C) \neq
\mathrm{AT}(C)$.  Indeed, that proof produces a short exact
sequence $$0 \to J \to C \to B \to 0,$$ where $C$ is residually
finite dimensional, $J$ is an AF algebra (hence $C$ has real rank
zero by \cite[Theorem 3.14 and Corollary 3.16]{brown-pedersen}),
the sequence is semi-split (hence $C$ is exact and satisfies the
UCT) and, finally, $C$ will have a non-amenable trace which is
gotten by composing the quotient map $C \to B$ with the
non-amenable trace on $B$ (since $C$ is locally reflexive and this
trace has a non-hyperfinite GNS representation).
\end{proof}

\begin{cor}
\label{thm:notTAF} There exists an exact, Popa algebra, $A$, with
real rank zero, stable rank one, UCT, unperforated K-theory, Riesz
interpolation and Blackadar's fundamental comparison property
which is approximately divisible and an increasing union of
residually finite dimensional subalgebras but such that $A$ is not
tracially AF.
\end{cor}

\begin{proof}
Since $\TA = \TAulfd$ for every tracially AF algebra (see
Proposition \ref{thm:characterization}), the example in the
previous theorem can't be tracially AF.
\end{proof}

It may seem unusual to mention in the previous two results that
$A$ is an inductive limit of residually finite dimensional
subalgebras.  Our main reason for pointing out this fact is that
some of Lin's recent structural work on the class of tracially AF
C$^*$-algebras relies heavily on a theorem of Blackadar and
Kirchberg stating that every simple, nuclear, quasidiagonal
C$^*$-algebra is an inductive limit of such subalgebras.  In fact,
for some of Lin's structural work, this is the only place that
nuclearity is used (i.e.\ his results hold more generally if one
replaces the assumption of nuclearity by the assumption of an
inductive limit decomposition by residually finite dimensional
subalgebras).

\begin{rem}
In \cite{bedos:hypertraces} B\'{e}dos asked whether or not every
separable, unital hypertracial C$^*$-algebra is nuclear (see
\cite[Section 3]{bedos:hypertraces} - in the language of the
present paper, a C$^*$-algebra is hypertracial if every quotient
has at least one amenable trace).  It is easy to see that every
simple, unital, quasidiagonal C$^*$-algebra is hypertracial and
hence Dadarlat's examples of non-nuclear Popa algebras provide
counterexamples to this question
\cite{dadarlat:nonnuclearsubalgebras}.  Theorem
\ref{thm:exoticpopaalgebra} above provides further examples.
Indeed, for every non-hyperfinite II$_1$-factor $M$ which contains
a weakly dense exact C$^*$-subalgebra the proof of Theorem
\ref{thm:exoticpopaalgebra} shows that we can construct an exact
Popa algebra with stable rank one, Blackadar's comparison property
(hence unperforated K-theory), Riesz property, approximate
divisibility and which is an increasing union of residually finite
dimensional subalgebras but which is not nuclear since it will
have $M\bar{\otimes} R$ as the weak closure of some GNS
representation.
\end{rem}

We now turn to Popa's question.  Note that we can't use the
examples from Corollary \ref{thm:notTAF} to produce
counterexamples to Popa's question since the GNS representation of
any quasidiagonal, locally reflexive (e.g.\ exact) C$^*$-algebra
with unique trace will yield the hyperfinite II$_1$-factor (which
was the obstruction used in the proof of Theorem
\ref{thm:exoticpopaalgebra}).  Hence to get a Popa algebra with
unique trace which is not tracially AF we are forced to step
outside the world of exact C$^*$-algebras.

\begin{thm}
\label{thm:uniquetracenotTAF} There exists a quasidiagonal
C$^*$-algebra $A$ with real rank zero, stable rank one, strong
Dixmier property (i.e.\ for every $a \in A$ the norm closed,
convex hull of $\{ uau^*: u \in {\mathcal U}(A)\}$ intersects
${\mathbb C}1_{A}$) -- hence $A$ is simple, has unique trace and
is a Popa algebra -- but such that $A$ is not tracially AF.
\end{thm}

\begin{proof} Though the construction is somewhat technical (being
quite similar to the inductive limit techniques from Section \ref{sec:Popareps}),
the main idea is not too hard to see when presented abstractly.
Our goal is to construct an inductive system $A_0 \to A_1 \to A_2
\to \cdots$, with connecting maps $\sigma_{n,m} : A_m \to A_n$ and
$\{ a^{(m)}_k \}_{k \in {\mathbb N}}$ denoting dense sequences in
the unit balls of the $A_m$'s, where all of the following
conditions are satisfied:

\begin{enumerate}
\item For each $n$ there is a finite set of self-adjoint elements
${\mathcal S}_n \subset (A_n)_{sa}$ such that:

  \begin{enumerate}
  \item Each element of ${\mathcal S}_n$ has finite spectrum.

  \item For each $x \in \sigma_{n,0}(\{ a^{(0)}_k \}_{k = 1}^n) \cup \ldots \cup
  \sigma_{n,n-1}(\{ a^{(n-1)}_k \}_{k = 1}^n)$ there exists an element
  $y \in {\mathcal S}_n$ such that $\| (x + x*) - y \| < 1/n$.
  \end{enumerate}

\item For each $n$ there is a finite set of invertible elements
${\mathcal I}_n \subset A_n$ such that for each $x \in
\sigma_{n,0}(\{ a^{(0)}_k \}_{k = 1}^n) \cup \sigma_{n,1}(\{
a^{(1)}_k \}_{k = 1}^n) \cup \ldots \cup \sigma_{n,n-1}(\{
a^{(n-1)}_k \}_{k = 1}^n)$ there exists an element $y \in
{\mathcal I}_n$ with $\| x - y \| < 1/n$.

\item For each $n$ there is a finite set of unitary elements
${\mathcal U}_n \subset A_n$ such that for each $x \in
\sigma_{n,0}(\{ a^{(0)}_k \}_{k = 1}^n) \cup \sigma_{n,1}(\{
a^{(1)}_k \}_{k = 1}^n) \cup \ldots \cup \sigma_{n,n-1}(\{
a^{(n-1)}_k \}_{k = 1}^n)$ there exists a complex number
$\lambda(x)$ and positive real numbers $\{ \theta^{(x)}_1, \ldots,
\theta^{(x)}_{k(x)} \}$ and unitaries $\{ u_1, \ldots, u_{k(x)} \}
\subset {\mathcal U}_n$ such that $\sum \theta^{(x)}_i = 1$ and
$$\|\lambda(x) 1_{A_n} - \sum_{i=1}^{k(x)} \theta^{(x)}_i u_i x
u_i^* \| < 1/n.$$

\item For each $n$ there is a nonzero, finite dimensional
subalgebra $B_n \subset A_n$ with unit $e_n$ such that for each $x
\in \sigma_{n,0}(\{ a^{(0)}_k \}_{k = 1}^n) \cup \sigma_{n,1}(\{
a^{(1)}_k \}_{k = 1}^n) \cup \ldots \cup \sigma_{n,n-1}(\{
a^{(n-1)}_k \}_{k = 1}^n)$, $[x, e_n] = 0$ and $e_n x \in B_n$.

\item For each $n$, the $*$-homomorphism $\sigma_{n+1,n}$ will be
faithful when restricted to the finite dimensional subalgebras
$\sigma_{n,0}(B_0), \ \sigma_{n,1}(B_1), \ \ldots,
\sigma_{n,n-1}(B_{n-1})$, (though the connecting maps will not be
injective on all of $A_n$).
\end{enumerate}

If we are successful in constructing such an inductive system and
we let $A$ denote the inductive limit then it is clear that
condition (1) will imply that $A$ has real rank zero, (2) will
imply stable rank one, (3) will force the strong Dixmier property
(which will imply simplicity of $A$) and thus (by simplicity) (4)
and (5) will imply that $A$ is a Popa algebra (hence quasidiagonal
-- which implies the existence of a tracial state -- which is
necessarily unique by the Dixmier property).  In other words, $A$
satisfies all of the hypotheses stated in the theorem.

In order to deduce that $A$ is not tracially AF, it suffices to
show that its unique trace $\gamma$ is {\em not}  uniform locally
finite dimensional.  To deduce this, we will also arrange that
there exists a C$^*$-algebra $C$ with a trace $\tau$, where $\tau$
is known to {\em not} be uniform amenable, and there exists a
sequence of $*$-homomorphisms $\rho_n : C \to A$ such that $\|
\gamma\circ\rho_n - \tau \|_{C^*} \to 0$. Since $\tau$ is not
uniform amenable it will follow that $\gamma$ can't satisfy the
approximation property of a uniform amenable trace either (for if
it did then this approximation property would apparently pass to
$\tau$).  In particular, $\gamma$ can't have the stronger property
of being uniform locally finite dimensional and since every trace
on a tracially AF algebra has this stronger property it follows
that $A$ can't be tracially AF.

We will use a recursive procedure to construct our inductive
system which we now describe abstractly (as we again feel that
this is easier to digest).  The starting point is to consider a
sequence of matrix algebras $\{ M_{k(n)}({\mathbb C})\}$ and a
(separable, unital) subalgebra $A_0 \subset \Pi M_{k(n)}({\mathbb
C})$.  We claim that if ${\mathfrak F} \subset A_0$ is any finite
set then we can find a (separable, unital) C$^*$-algebra $A_1
\subset \Pi M_{l(n)}({\mathbb C})$ (note that the dimensions of
the matrix algebras have changed) and a $*$-homomorphism $\sigma :
A_0 \to A_1$ such that the five properties listed above all hold.

When $A_0 \subset \Pi M_{k(n)}({\mathbb C})$ is given we can
define tracial states $\tau_n$ on $A_0$ by first cutting to the
$n^{th}$ summand of $\Pi M_{k(n)}({\mathbb C})$ and then applying
the unique tracial state on $M_{k(n)}({\mathbb C})$.  By
compactness, $\{\tau_n \}$ has a weak-$*$ convergent subsequence,
and hence for any $\epsilon
> 0$ we can find a further subsequence $\{ \tau_{n(j)} \}$ such that
for every $x \in {\mathfrak F}$, $| \tau_{n(i)} (x) - \lim_{j \to
\infty} \tau_{n(j)}(x) | < \epsilon$ for all $i \in {\mathbb N}$.
Cutting $\Pi_{n \in {\mathbb N}} M_{k(n)}({\mathbb C})$ down to
the subsequence $\{ n(j) \}$ we get a $*$-homomorphism $\eta : A_0
\to \Pi_{j \in {\mathbb N}} M_{k(n(j))}({\mathbb C})$.  Since von
Neumann algebras have the usual Dixmier property (cf.\
\cite[Theorem 8.3.5]{kadison-ringrose}) we can find a finite set
of unitaries ${\mathcal U} \subset \Pi_{j \in {\mathbb N}}
M_{k(n(j))}({\mathbb C})$ such that for each $x \in {\mathfrak F}$
some convex combination of unitary conjugates of $x$ by unitaries
from ${\mathcal U}$ will be within $\epsilon$ of the center of
$\Pi_{j \in {\mathbb N}} M_{k(n(j))}({\mathbb C})$.  Since we have
arranged that $| \tau_{n(i)} (x) - \lim_{j \to \infty}
\tau_{n(j)}(x) | < \epsilon$ for all $i \in {\mathbb N}$ and $x
\in {\mathfrak F}$ this means that we can find positive real
numbers $\{ \theta^{(x)}_1, \ldots, \theta^{(x)}_{k(x)} \}$ and
unitaries $\{ u_1, \ldots, u_{k(x)} \} \subset {\mathcal U}$ such
that $\sum \theta^{(x)}_i = 1$ and $$\|(\lim_{j \to \infty}
\tau_{n(j)}(x)) 1 - \sum_{i=1}^{k(x)} \theta^{(x)}_i u_i \eta(x)
u_i^* \| < \epsilon.$$ Since $\Pi_{j \in {\mathbb N}}
M_{k(n(j))}({\mathbb C})$ has both real rank zero and stable rank
one, we can also find a finite set of self-adjoints ${\mathcal S}
\subset \Pi_{j \in {\mathbb N}} M_{k(n(j))}({\mathbb C})$ and a
finite set of invertibles ${\mathcal I} \subset \Pi_{j \in
{\mathbb N}} M_{k(n(j))}({\mathbb C})$ which approximate
$\eta({\mathfrak F}$ appropriately.  Letting $D_0 = C^*(
\eta(A_0), {\mathcal U}, {\mathcal S}, {\mathcal I})$ we get a
(separable, unital) algebra with the first three of the desired
approximation properties.  At this point it would be most natural
to just repeat this procedure, but unfortunately this will not
necessarily produce a limit algebra which is quasidiagonal (cf.\
\cite[Remark 9.4]{brown:QDsurvey}).  Hence we are forced to add
one more complication.  Namely, consider the natural finite
dimensional representation $\pi_{n(1)} : D_0 \to
M_{k(n(1))}({\mathbb C})$. Let $A_1 \subset
\bigg(M_{k(n(1))}({\mathbb C}) \otimes M_t ({\mathbb C})\bigg)
\otimes \Pi_{j \in {\mathcal N}} M_{k(n(j))}({\mathbb C})$ be the
C$^*$-algebra generated by $M_{k(n(1))}({\mathbb C}) \otimes M_t
({\mathbb C})$ and $D_0$ (here $t > 1/\epsilon$ is just some
really big integer) and now define a (unital) $*$-homomorphism
$D_0 \to A_1$ exactly as we did in the basic construction employed
in the proof of Theorem \ref{thm:basicconstruction} (i.e.\ replace
the natural inclusion by one twisted by the finite dimensional
representation $\pi_{n(1)}$). Composing $D_0 \to A_1$ with the map
$\eta$ we get a $*$-homomorphism $\sigma : A_0 \to A_1$ where we
have now arranged properties (1) - (4). Note that (at the next
step) condition (5) will be automatic since the finite dimensional
$B_1 \cong M_{k(n(1))}({\mathbb C})$ lives in the left tensor
factor of $\bigg(M_{k(n(1))}({\mathbb C}) \otimes M_t ({\mathbb
C})\bigg) \otimes \Pi_{j \in {\mathcal N}} M_{k(n(j))}({\mathbb
C}) \cong \Pi_{j \in {\mathcal N}} M_{k(n(1))}({\mathbb C})
\otimes M_t ({\mathbb C})\otimes M_{k(n(j))}({\mathbb C})$.  Also
note that if we choose the integer $t$ large enough and if $\tau$
denotes the weak-$*$ limit of $\{\tau_n \}$ (from the beginning of
this paragraph) it is clear from the construction that $A_1$ has a
trace $\gamma$ such that $\| \tau - \gamma \circ \sigma
\|_{A_0^*}$ is as small as you like.  Finally, we remark that
there is also the natural embedding $D_0 \hookrightarrow A_1 \cong
M_{k(n(1))}({\mathbb C}) \otimes M_t ({\mathbb C})\otimes D_0$
(which should not be forgotten as we will need it later).

The (actual!) construction of our mystery algebra $A$ begins with
any residually finite, discrete, non-amenable group $\Gamma$.  Let
$\Gamma \triangleright \Gamma_1 \triangleright \Gamma_2
\triangleright \cdots$ be a sequence of normal subgroups, each of
finite index, such that their intersection is the neutral element.
Let $\pi_i : C^*(\Gamma) \to B(l^2(\Gamma / \Gamma_i))$ be the
representations induced from the left regular representations of
the $\Gamma / \Gamma_i$'s and let $\pi = \oplus \pi_i$.  A key
remark is that if we choose any subsequence $\{ i_k \}$ then
$\mathrm{tr}_{i_k} \circ \pi_{i_k} (x) \to \tau(x)$ for all $x \in
C^*(\Gamma)$, where $\mathrm{tr}_{i_k}$ is the unique tracial
state on $B(l^2(\Gamma / \Gamma_{i_k}))$ and $\tau$ is the trace
on $C^*(\Gamma)$ whose GNS representation gives the left regular
representation.

Letting $A_0 = \pi(C^*(\Gamma))$ and recursively applying the
scheme outlined above, we get a sequence $A_0 \to A_1 \to \cdots$
with all five properties listed at the beginning of the proof and
hence we get a limit algebra $A$ satisfying all of the hypotheses
of the theorem. Finally, letting $C = C^*(\Gamma)$ it is easy to
see that there exists a sequence of $*$-homomorphisms $\psi_n : C
\to A_n$ and traces $\gamma_n$ on $A_n$ such that $\tau = \gamma_n
\circ \psi_n$ (use the natural embedding $D_0 \hookrightarrow A_1$
from above) and hence (composing with the natural maps $A_n \to
A$) we get a sequence of $*$-homomorphisms $\rho_n : C \to A$ such
that $\| \tau - \gamma\circ\rho_n \|_{C^*} \to 0$ (this is
basically the same as statement (1) in Theorem
\ref{thm:basicconstruction}).  Since $\Gamma$ was assumed
non-amenable (hence $\tau$ is not a uniform amenable trace -- see
Proposition \ref{thm:what?}) the proof is now complete.
\end{proof}

\begin{rem} Note that the examples of the previous theorem also
answer Popa's weaker question of whether or not a Popa algebra
with unique trace must give the hyperfinite II$_1$-factor in its
GNS representation. Indeed, we deduced that the previous examples
were not tracially AF from the fact that their unique traces were
not uniform amenable. This is equivalent, by Theorem
\ref{thm:mainthmUTAwafd}, to asserting that their GNS
constructions do not yield hyperfinite factors and hence there
exists a Popa algebra $A$ with unique trace $\tau$ such that
$\pi_{\tau}(A)^{\prime\prime} \ncong R$.
\end{rem}

\section{Connes' embedding problem}
\label{thm:connesproblem}

In this section we show that the techniques of this paper yield
new characterizations of those II$_1$-factors which are embeddable
into $R^{\omega}$.  Connes' embedding problem asks whether {\em
every} separable II$_1$-factor embeds into $R^{\omega}$.  Though
it appears to be a technical question only of interest to
specialists, it is a remarkable fact that this problem is
equivalent to numerous other questions ranging from uniqueness of
certain tensor product norms, existence of faithful traces,
(local) universality of the trace class operators among all
noncommutative $L^1$-spaces, Kirchberg's `QWEP conjecture',
density of finite dimensional factors (in the Effros-Mar\'echal
topology) in the space of all finite factors and a certain
statement about the space of moments of noncommutative random
variables in an arbitrary finite W$^*$-probability space.  (We
highly recommend Ozawa's survey article \cite{ozawa:QWEP} for
precise statements of these equivalences -- with a number of
simplified proofs -- as well as a comprehensive bibliography on
the subject.)  Moreover, thanks to the deep work of Haagerup
\cite{haagerup:invariantsubspaces}, a positive answer to Connes'
problem would nearly complete the relative invariant subspace
problem for II$_1$-factors while a negative answer would imply a
negative answer to Voiculescu's `unification problem' (i.e.\ the
question of whether or not the microstates and non-microstates
approaches to free entropy theory agree).  Finally we mention that
it seems quite likely that geometric group theory -- in
particular, Gromov's hyperbolic groups -- will play an important
role in future developments around this problem.  In short, this
seemingly technical question is quite deep, intimately related to
a surprising variety of other problems and, as such, is a very
important open question.

II$_1$-factors which are embeddable into $R^{\omega}$ already admit a
number of characterizations (see \cite{kirchberg:invent},
\cite{haagerup-winslow}, \cite{ozawa:QWEP}) but the results of this
section show that the difference between embeddability and
hyperfiniteness is rather delicate.

\begin{thm}
\label{thm:embeddableMcDuff}
If $M \subset B(L^2(M))$ is a II$_1$-factor with trace $\tau_M$ then
the following are equivalent:

\begin{enumerate}
\item $M$ is embeddable into $R^{\omega}$.

\item There exists a weakly dense C$^*$-subalgebra $A \subset M$ and a
sequence of normal, u.c.p.\ maps $\varphi_n : M \to M_{k(n)} ({\mathbb
C})$ such that

\begin{enumerate}
\item $\| \varphi_n (ab) - \varphi_n (a) \varphi_n (b) \|_2 \to 0$ and

\item $| \tau_{k(n)} \circ \varphi_n (a) - \tau_M (a) | \to 0$ for all
$a,b \in A$.
\end{enumerate}

\item There exists a weakly dense C$^*$-subalgebra $A \subset M$ and
finite rank projections $\{ P_n \}$ such that

\begin{enumerate}
\item $\frac{ \| [P_n, a] \|_{HS}}{\| P_n \|_{HS}} \to 0$ and

\item $\frac{\langle aP_n, P_n\rangle_{HS}}{\langle
P_n,P_n\rangle_{HS}} \to \tau_M (a)$ for all $a \in A$.
\end{enumerate}

\item There exists a weakly dense C$^*$-subalgebra $A \subset M$ and a
state $\phi$ on $B(L^2(M))$ such that

\begin{enumerate}
\item $A \subset B(L^2(M))_{\phi}$\footnote{Contrary to previous
notation, $B(L^2(M))_{\phi}$ does {\em not} denote the
multiplicative domain. We are following classical notation and now
using $B(L^2(M))_{\phi}$ to denote the centralizer of the state
$\phi$ just as in Corollary \ref{thm:BFDV}.} and

\item $\phi|_A = \tau_M.$
\end{enumerate}

\item There exists a weakly dense C$^*$-subalgebra $A \subset M$ and
an embedding $\Phi : M \hookrightarrow R^{\omega}$ such that $\Phi|_A$
is completely positively liftable.

\item There exists a weakly dense C$^*$-subalgebra $A \subset M$ such
that $C^*(A, JAJ) \cong A \otimes A^{op}$.

\item $M$ has a weak expectation relative to a weakly dense
C$^*$-subalgebra. (i.e. There exists a weakly dense C$^*$-subalgebra
$A \subset M$ and a u.c.p.\ map $\Phi : B(L^2(M)) \to M$ such that
$\Phi|_A = id_A$.)

\item There exists a weakly dense operator system $X \subset M$ such
that $X$ is injective.

\item For every finite set ${\mathfrak F} \subset M$ and $\varepsilon
> 0$ there exists a complete order embedding $\Phi : R \hookrightarrow
M$ (i.e.\ $\Phi$ is an operator system isomorphism from $R$ to
$\Phi(R)$) such that for each $x \in {\mathfrak F}$ there exists $r
\in R$ such that $\| x - \Phi(r) \|_2 < \varepsilon$.
\end{enumerate}
\end{thm}

\begin{proof}  Our study of representation theory will provide the
key to proving this result.  Indeed, the first thing we will do is
identify $C^*({\mathbb F}_{\infty})$ with a weakly dense
subalgebra of $M$ (which is possible thanks to Proposition
\ref{thm:universalfaithful}).

(1) $\Longrightarrow$ (2). If $M$ embeds into $R^{\omega}$ it
follows that $\tau_M|_{C^* ({\mathbb F}_{\infty})}$ is an amenable
trace on $C^* ({\mathbb F}_{\infty})$ (this fact was already used
in the proof of Proposition \ref{thm:freegroupexample}).  Hence we
can find u.c.p.\ maps (even $*$-homomorphisms if you like)
$\varphi_n : C^* ({\mathbb F}_{\infty}) \to M_{k(n)} ({\mathbb
C})$ which are asymptotically multiplicative (in 2-norm) and
recover $\tau_M|_{C^* ({\mathbb F}_{\infty})}$ as a weak-$*$ limit
after composing with $\mathrm{tr}_{k(n)}$.  As we have seen
before, one completes the proof by first extending the
$\varphi_n$'s to all of $M$ (via Arveson's theorem) and then
replacing the (not necessarily normal) extensions with normal
u.c.p.\ maps (which is possible since there is a one-to-one
correspondence between c.p.\ maps $M \to M_k(\mathbb{C})$ and
positive linear functionals on $M_k(M)$ and any positive linear
functional can be approximated by normal, positive linear
functionals).

Note that statement (2) above is exactly the same as saying that
$M$ contains a weakly dense subalgebra $A \subset M$ such that
$\tau_M|_A$ is an amenable trace on $A$.  With this reformulation
it follows easily from Theorems \ref{thm:amenabletraces} and
\ref{thm:amenabletracesII} that (2), (3), (4) and (5) are
equivalent.  Since (5) obviously implies (1) we have shown the
equivalence of (1) - (5).

We now wish to add (6) and (7) to the list of equivalent
statements so first assume that there exists a weakly dense
algebra $A \subset M$ such that $\tau_M|_A$ is an amenable trace
on $A$.  By the fourth statement in Theorem
\ref{thm:amenabletraces} it follows that the natural map $A \odot
A^{op} \subset M\odot M^{op} \to B(L^2(M))$ is continuous with
respect to the minimal tensor product norm on $A \odot A^{op}$.
However, a classical result of Murray and von Neumann states that
$M\odot M^{op} \to B(L^2(M))$ is injective and hence $C^*(A, JAJ)$
is necessarily isomorphic to $A \otimes A^{op}$.   In other words,
(1) $\Longrightarrow$ (6).  That (6) $\Longrightarrow$ (7) is just
another application of Lance's trick (cf.\ the proof of (4)
$\Longrightarrow$ (5) from Theorem \ref{thm:amenabletraces}) while
that fact that (7) implies $\tau_M|_A$ is an amenable trace on $A$
is already contained in the proof of the last assertion of Theorem
\ref{thm:amenabletraces}.  Hence we have shown that (1) - (7) are
all equivalent.

(8) $\Longrightarrow$ (1).  In \cite{kirchberg:invent} Kirchberg
proved the following result:  $M \subset R^{\omega}$ if and only
if there exists a C$^*$-algebra $B$ with Lance's WEP and a u.c.p.\
map $\Phi:B \to M$ such that the unit ball of $B$ gets mapped onto
a weakly dense subset of the unit ball of $M$ (cf.\ \cite[Theorem
1.4]{kirchberg:invent} and the equivalence of ($vi$) and ($iii$)
in \cite[Proposition 1.3]{kirchberg:invent}).  Since $B(H)$ is
injective it has the WEP and hence this implication follows from
Kirchberg's result.

(7) $\Longrightarrow$ (8).  If $A \subset M$ is weakly dense and
$\Phi : B(L^2(M)) \to M$ is a u.c.p.\ map such that $\Phi(a) = a$
for all $a \in A$ then \cite[Theorem 2.1]{blackadar:Proc.WEP}
ensures that we can find an {\em idempotent} u.c.p.\ map $\Psi :
B(L^2(M)) \to M$ such that $\Psi (a) = a$ for all $a \in A$.  The
desired injective operator system is then $X = \Psi(B(L^2(M)))$
(since $\Psi \circ \Psi = \Psi$).

The proof will be complete once we observe the implications (7)
$\Longrightarrow$ (9) and (9) $\Longrightarrow$ (1).

(7) $\Longrightarrow$ (9).  Given ${\mathfrak F} \subset M$ and
$\varepsilon > 0$, let $p \in M$ be a projection whose trace is
very close to 1.  By the equivalence of (7) and (1), together with
Proposition \ref{thm:universalfaithful}, we can find a weakly
dense copy of $C^* ({\mathbb F}_{\infty})$ inside (the
II$_1$-factor) $pMp$ such that there exists a u.c.p.\ map $\Phi :
B(L^2(pMp)) \to pMp$ with the property that $\Phi(x) = x$ for all
$x$ in the dense copy of $C^* ({\mathbb F}_{\infty}) \subset pMp$.
Since $C^* ({\mathbb F}_{\infty})$ also embeds into $R$ we can use
Arveson's extension theorem followed by $\Phi$ to construct a
u.c.p.\ map $\Psi : R \to pMp$ with weakly dense range.  Looking
to the orthogonal side, we can also find a unital, normal,
embedding $\rho : R \hookrightarrow (1-p)M(1-p)$.  It follows that
$\Psi \oplus \rho : R \to pMp \oplus (1-p)M(1-p) \subset M$ is a
unital, complete order embedding (since $\rho$ is an injective
$*$-homomorphism) whose image nearly contains ${\mathfrak F}$ in
2-norm.

(9) $\Longrightarrow$ (1). Assuming (9), we will show that one can
construct a u.c.p.\ map $$\Psi : \Pi_{\mathbb N} R \to M$$ such
that the image of the unit ball of $\Pi_{\mathbb N} R$ is weakly
dense in the unit ball of $M$.  From Kirchberg's result quoted
above it will follow that $M$ is embeddable into $R^{\omega}$
since $\Pi_{\mathbb N} R$ is injective and hence has Lance's WEP.
The construction of $\Psi$ is a fairly standard
density/ultrafilter argument which goes as follows.  Fix a
sequence $a_1, a_2, \ldots$ which is weakly dense in the unit ball
of $M$ and for each $n \in \mathbb{N}$ let $\phi_n:R \to M$ be a
complete order embedding whose range almost contains the first $n$
elements in the sequence $\{a_i\}$ (to within, say,
$\frac{1}{n}$). Now fix a free ultrafilter $\omega$ and define a
u.c.p.\ map $$\Psi : \Pi_{\mathbb N} R \to M$$ by declaring
$$\Psi( (r_n)_{n\in \mathbb{N}}) = \lim_{n \to \omega} \phi_n
(r_n)$$ and it is easily seen that this takes the unit ball onto a
weakly dense subset of the unit ball of $M$.
\end{proof}

\begin{rem} In \cite{lance:nuclear} Lance introduced the WEP as a
property of C$^*$-algebras.  In \cite{blackadar:Proc.WEP}
Blackadar changed the point of view to make this a W$^*$-notion as
follows: A von Neumann algebra $M \subset B(H)$ has the WEP
relative to a weakly dense C$^*$-subalgebra $A \subset M$ if there
exists a u.c.p.\ map $\Phi : B(H) \to M$ such that $\Phi(a) = a$
for all $a \in A$.  Lance had earlier speculated that this
`W$^*$-WEP' would actually imply injectivity (indeed, a von
Neumann algebra has the `C$^*$-WEP' if and only if it is
injective) however a (non-factor) counterexample was given in
\cite{blackadar:Proc.WEP}.  In \cite{blackadar:WEP} and
\cite[Corollary 3.5]{kirchberg:invent} it was shown that there
exist non-hyperfinite {\em factors} with the W$^*$-WEP, but the
proofs were non-constructive and concrete examples remained
elusive. In a preliminary version of this paper we gave the first
concrete examples of non-injective factors (actually McDuff
II$_1$-factors) with the W$^*$-WEP.  Shortly after that, in joint
work with Dykema \cite{BD}, it was shown that all the
(interpolated) free group factors on a finite number of generators
have the W$^*$-WEP.  However, the previous theorem shows that most
of the standard examples of II$_1$-factors (e.g.\ things coming
from residually finite groups like free groups, $SL(n, {\mathbb
Z})$ or ${\mathbb Z}^2 \rtimes SL(2, {\mathbb Z})$) have the
W$^*$-WEP and hence we get examples with property T and even
factors with trivial fundamental group \cite{popa:fundgroup}.

In other words, the W$^*$-WEP is not an exotic property and it is
enjoyed by many II$_1$-factors which are quite far from being
hyperfinite.  On the other hand, we strongly believe that this notion
is still `too close' to injectivity to expect that {\em every}
II$_1$-factor has this property.  For example, from part (9) of
Theorem \ref{thm:embeddableMcDuff} we get that a II$_1$-factor has the
W$^*$-WEP if and only if it is quite literally built out the
hyperfinite II$_1$-factor in a way which naturally mixes von Neumann
algebraic and operator space notions.  Our feeling is that this is not
that far ({\em very} loosely speaking) from being injective -- far
enough to exhibit vastly different properties like no Cartan
subalgebras, property T, trivial fundamental group, etc.\ -- but,
still, not that far.
\end{rem}

In order to illustrate just how delicate the results above are, we
remind the reader of the various characterizations of the hyperfinite
II$_1$-factor (most of which are due to Alain Connes).

\begin{thm}
\label{thm:R}
If $M \subset B(L^2(M))$ is a II$_1$-factor with trace $\tau_M$ then
the following are equivalent:

\begin{enumerate}
\item $M \cong R$.

\item There exists a weakly dense C$^*$-subalgebra $A \subset M$ and a
sequence of normal, u.c.p.\ maps $\varphi_n : M \to M_{k(n)} ({\mathbb
C})$ such that

\begin{enumerate}
\item $\| \varphi_n (ab) - \varphi_n (a) \varphi_n (b) \|_2 \to 0$ for
all $a,b \in A$ and

\item $| \tau_{k(n)} \circ \varphi_n (x) - \tau_M (x) | \to 0$ for all
$x \in M$.
\end{enumerate}

\item There exists a weakly dense C$^*$-subalgebra $A \subset M$ and
finite rank projections $\{ P_n \}$ such that

\begin{enumerate}
\item $\frac{ \| [P_n, a] \|_{HS}}{\| P_n \|_{HS}} \to 0$ for every $a
\in A$ and

\item $\frac{\langle xP_n, P_n\rangle_{HS}}{\langle
P_n,P_n\rangle_{HS}} \to \tau_M (x)$ for all $x \in M$.
\end{enumerate}

\item There exists a weakly dense C$^*$-subalgebra $A \subset M$ and a
state $\phi$ on $B(L^2(M))$ such that

\begin{enumerate}
\item $A \subset B(L^2(M))_{\phi}$ and

\item $\phi|_M = \tau_M.$
\end{enumerate}

\item There exists a completely positively liftable embedding $\Phi :
M \hookrightarrow R^{\omega}$.

\item $C^*(M, JMJ) \cong M \otimes M^{op}$.

\item There exists a u.c.p.\ map $\Phi : B(L^2(M)) \to M$ such that
$\Phi|_M = id_M$ (i.e.\ $M$ is injective).

\item For every finite set ${\mathfrak F} \subset M$ and $\varepsilon
> 0$ there exists a $*$-monomorphism $\Phi : R \hookrightarrow M$ such
that for each $x \in {\mathfrak F}$ there exists $r \in R$ such that
$\| x - \Phi(r) \|_2 < \varepsilon$.

\end{enumerate}
\end{thm}

\begin{proof} (1) $\Longrightarrow$ (2) is obvious.
(2) $\Longrightarrow$ (5) is contained in the argument in the proof of
(1) $\Longrightarrow$ (2) from Theorem \ref{thm:mainthmUTAwafd}.
Since (5) is equivalent to hyperfiniteness for general finite von
Neumann algebras, we see that (1), (2) and (5) are equivalent.

The equivalence of (1), (6) and (7) is due to Connes (cf.\
\cite[Theorem 5.1]{connes:classification}).  This paper also contains
his adaptation of Day's trick to deduce (3) from injectivity (7).
Note also that the equivalence of (8) and (1) goes all the way back to
Murray and von Neumann.  Hence we are left to prove (3)
$\Longrightarrow$ (4) and (4) $\Longrightarrow$ (1).

(3) $\Longrightarrow$ (4) is well known: simply take a cluster
point of the states $$T \mapsto \frac{\langle TP_n, P_n
\rangle_{HS}}{\langle P_n, P_n \rangle_{HS}}$$ on $B(L^2(M))$.
Finally, (4) $\Longrightarrow$ (1) also follows from Connes' work
since the density of $A$ in $M$, together with the fact that the
hypertrace takes the correct value on all of $M$, implies that
actually $M$ is contained in the centralizer of $\varphi$ and
hence is hyperfinite (cf.\ proof of \cite[Theorem
5.1]{connes:compactmetricspaces}).
\end{proof}

We close this section with a simple result relating amenable
traces on C$^*$-algebras with the {\em local lifting property}
(LLP) and Connes' embedding problem.  We have already seen that in
many instances every trace on a particular C$^*$-algebra will
enjoy some type of amenability (e.g.\ every trace on a
C$^*$-algebra with the WEP is amenable).  We now observe that
Connes' embedding problem predicts that every trace on a
C$^*$-algebra with the LLP is amenable and hence it is possible
that a counterexample to Connes' problem could be constructed by
finding a non-amenable trace on some C$^*$-algebra with the LLP.

By definition, $A$ has the local lifting property if for every
C$^*$-algebra $B$, ideal $J \triangleleft B$, u.c.p.\ map $\phi :
A \to B/J$ and finite dimensional operator system $X \subset A$
there exists a u.c.p.\ lifting of $\phi|_X$.

\begin{prop}
\label{thm:LLPexample}  Assume that $A$ has the LLP and $\tau \in
\TA$. Then $\tau$ is amenable if and only if there exists a
$\tau^{\prime\prime}$-preserving embedding
$$\pi_{\tau}(A)^{\prime\prime} \hookrightarrow R^{\omega},$$ where
$\tau^{\prime\prime}$ denotes the normal extension of $\tau$ to
$\pi_{\tau}(A)^{\prime\prime}$.
\end{prop}

\begin{proof} This is very similar to the equivalence of statements
(1) and (2) from Theorem \ref{thm:amenabletracesII} and hence the
details are left to the reader.
\end{proof}

For all intents and purposes the following result is due to
Kirchberg (cf.\ \cite{kirchberg:invent}) however the proof is very
simple so we include it.

\begin{prop}
\label{thm:BFDII}
 The following statements are equivalent.
\begin{enumerate}
\item Every (separable) II$_1$-factor embeds into $R^{\omega}$.

\item $\TA = \TAim$ for every $A$ with the LLP.

\item $\mathrm{T}(C^*(\mathbb{F}_{\infty})) =
\mathrm{AT}(C^*(\mathbb{F}_{\infty}))$.
\end{enumerate}
\end{prop}

\begin{proof}
(1) $\Longrightarrow$ (2).  Since $\TAim$ is always a face in
$\TA$ it follows that every trace on $A$ is amenable if and only
if every factorial trace on $A$ is amenable.  With this
observation and Proposition \ref{thm:LLPexample} in hand, it is
easily seen that an affirmative answer to Connes' problem would
imply (2).

(2) $\Longrightarrow$ (3). In \cite[Lemma
3.3]{kirchberg:commutantsofunitaries} Kirchberg showed that
$C^*(\mathbb{F}_{\infty})$ has the LLP.

(3) $\Longrightarrow$ (1).  We have already observed that every
II$_1$-factor arises as the GNS representation of
$C^*(\mathbb{F}_{\infty})$ with respect to some trace and hence
this part follows from Theorem \ref{thm:amenabletracesII}.
\end{proof}

\section{Amenable traces and numerical analysis}
\label{thm:numericalanalysis}

In this section we observe that the results and techniques of
these notes are relevant to some natural and important questions
which lie in the intersection of operator theory and numerical
analysis. In order to stay reasonably self-contained we will
recall the definitions and (statements of) results that we need,
however we highly recommend looking at the recent book of Hagen,
Roch and Silbermann \cite{HRS} (see also B\"{o}ttcher's survey
article \cite{bottcher}) for more on this subject.

Let $T \in B(H)$ be given.  There are  two important
questions that one may hope to (approximately) solve via a combination
of operator theoretic and numerical analytic techniques.  Namely, (1)
given a vector $v \in H$ find a vector $w \in H$ such that $Tw = v$
and (2) compute the spectrum, $\sigma(T)$, of $T$.  (In fact, there
are several other questions one may ask and we will see another
concerning spectral distributions of self-adjoint operators later in
this section. However, for now, we will stick to these two basic
problems.) If $H$ happens to be a finite dimensional Hilbert space
then, in many instances, there are efficient algorithms for attacking
these problems numerically.  For infinite dimensional $H$ the `finite
section method' provides a natural {\em strategy} for reducing these
questions to the finite dimensional case where one can then apply the
established numerical algorithms.  Of course, the real question -- at
least from a theoretical point of view -- is whether or not this
strategy will work.

To make this more precise, we need to introduce some terminology.  A
{\em filtration} of $H$ is an increasing sequence of finite rank
projections $P_1 \leq P_2 \leq \cdots$ such that $\| P_n (x) - x \|
\to 0$, as $n \to \infty$, for every vector $x \in H$.  The {\em
finite section method} is to take an operator $T \in B(H)$, cut it
down to a filtration (i.e.\ consider the sequence $P_n TP_n$ -- now
regarded as acting on finite dimensional Hilbert spaces) and hope that
by solving questions (1) and (2) -- presumably numerically -- for each
$P_n T P_n$ we can pass to a limit to recover the solutions for $T$.
While this idea could not be more natural, great care must be taken in
choosing the `right' filtration for a particular operator $T$ as it
can easily happen that $P_n T P_n$ do not provide good approximations
to $T$. (Think of cutting down the bilateral shift on $l^2({\mathbb
Z})$ to a filtration defined by the canonical basis.  The spectrum of
the bilateral shift is ${\mathbb T}$ while the cut-downs will all be
nilpotent and hence have spectrum 0.)

In relation to solving general vector equations, the following
definition naturally arises: we say that the finite section method
(w.r.t.\ a given filtration) is {\em stable} (for $T$) if there exists
a number $n_0$ such that for all $n \geq n_0$, $P_n T P_n$ are
invertible (as operators from $P_n H \to P_n H$) and $$\sup_{n\geq
n_0} \| (P_n T P_n)^{-1} \| < \infty.$$ The following result
summarizes the usefulness of this notion (cf.\ \cite[Theorem
1.4]{HRS}, \cite[Proposition 1.2]{bottcher}).

\begin{thm} Let $T \in B(H)$ and a filtration $P_1 \leq P_2 \leq \cdots$
be given.  If the finite section method is stable then:

\begin{enumerate}
\item $T$ is invertible.

\item Given a vector $v \in H$, let $v_n = P_n v$.  Then, for
sufficiently large $n$, there exists a (necessarily unique) vector
$w_n \in P_n H$ such that $(P_n T P_n)w_n = v_n$ and, moreover, $\{
w_n \}$ is a norm convergent sequence (in $H$) with limit $w$, where
$w$ is the (necessarily unique) solution to the equation $Tw = v$.
\end{enumerate}
\end{thm}

In other words, stability is the notion which ensures that there is
hope of solving (infinite dimensional) equations numerically.

To handle the spectral approximation problem we will need more terminology.

\begin{defn}

\begin{enumerate}
\item If ${X_n}$ is a sequence of subsets of the complex plane then we
let $\limsup\limits_{n \to \infty} X_n$ (resp.\ $\liminf\limits_{n \to
\infty} X_n$) denote the set of all cluster points (resp.\ limits) of
sequences $\{x_n \}$ with $x_n \in X_n$.  Another way of saying this
is that $\liminf\limits_{n \to \infty} X_n$ is the set of limits along
all possible sequences of points in $\{ X_n \}$ while
$\limsup\limits_{n \to \infty} X_n$ is the set of limits along all
possible subsequences of points.  Evidently we always have
$\liminf\limits_{n \to \infty} X_n \subset \limsup\limits_{n \to
\infty} X_n$.  (An example: Let $X_n = \{ (-1)^n + 1/n, 0 \}$.  Then
$\liminf\limits_{n \to \infty} X_n = \{ 0 \}$ while $\limsup\limits_{n
\to \infty} X_n = \{ -1, 0, 1\}$.)

\item Given $S \in B(H)$ and $\varepsilon > 0$, the {\em
$\varepsilon$-pseudospectrum of $S$}, denoted
$\sigma^{(\varepsilon)}(S)$, is the union of the usual spectrum,
$\sigma(S)$, together with the set of points $\lambda \in {\mathbb
C}\backslash \sigma(S)$ such that $\| (\lambda - S)^{-1} \| \geq
1/\varepsilon$.  Note that if $\varepsilon_1 \leq \varepsilon_2$ then
$\sigma^{(\varepsilon_1)}(S) \subset \sigma^{(\varepsilon_2)}(S)$ and,
also, $$\sigma(S) = \bigcap\limits_{\varepsilon > 0}
\sigma^{(\varepsilon)}(S).$$
\end{enumerate}
\end{defn}

The first definition above is simply a precise formulation of the
notion of convergence of spectra which will appear in our results.
The notion of pseudospectrum is relevant to our discussion because it
turns out that these sets behave better than actual spectra when
passing to limits.  More precisely, we will need the following result.

\begin{thm}\cite[Theorem 3.31]{HRS}
\label{thm:convergenceofspectra}
Let $(T_n) \in \Pi M_{k(n)}
({\mathbb C})$ and $\varepsilon > 0$ be given.  Then
$$\sigma^{(\varepsilon)}((T_n) + \oplus M_{k(n)} ({\mathbb C})) =
\limsup_{n \to \infty} \sigma^{(\varepsilon)}(T_n),$$ where $(T_n) +
\oplus M_{k(n)} ({\mathbb C})$ is the image of $(T_n)$ in the quotient
algebra $\Pi M_{k(n)} ({\mathbb C})/\oplus M_{k(n)} ({\mathbb C})$.
\end{thm}

For actual spectra one always has the inclusion $$\sigma((T_n) +
\oplus M_{k(n)} ({\mathbb C})) \supset \limsup_{n \to \infty}
\sigma(T_n),$$ but the inclusion can be proper.  However, it is an
interesting fact (cf.\ \cite[Corollary 3.8]{HRS}) that if each matrix
$T_n$ is normal then we get the equation $\sigma((T_n) + \oplus
M_{k(n)} ({\mathbb C})) = \limsup_{n \to \infty} \sigma(T_n).$

For an arbitrary element $T$ in a C$^*$-algebra we will let
$C^*(T)$ denote the unital C$^*$-algebra generated by $T$.  Also,
if $\{ P_n \}$ is a filtration of $H$ and $T \in B(H)$ then
$\varphi_n$ will denote the state on $C^*(T)$ given by $\varphi_n
(X) = \mathrm{tr}_{rank(P_n)} (P_n X P_n)$ for all $X \in C^*(T)$.

We are finally in a position to state the main theorem of this
section.

\begin{thm}
\label{thm:finitesection}
Let $T \in B(H)$ be a quasidiagonal operator and assume
that the C$^*$-algebra generated by $T$ is exact.  Then there exists a
filtration $P_1 \leq P_2 \leq \cdots$ such that:

\begin{enumerate}
\item The finite section method is stable if and only if $T$ is
invertible.

\item $\sigma(T) = \bigcap\limits_{\varepsilon > 0} \bigg(
\limsup\limits_{n \to \infty} \sigma^{(\varepsilon)}(P_n T
P_n)\bigg)$.

\item For every $\tau \in {\rm T(C^*(\pi(T)))_{QD}}$ there exists a
subsequence $\{n(k)\}$ such that $\varphi_{n(k)} \to \tau \circ
\pi|_{C^*(T)}$ in the weak-$*$ topology, where $\pi : B(H) \to Q(H)$
is the quotient mapping onto the Calkin algebra.
\end{enumerate}
\end{thm}

\begin{proof} We first claim that it will suffice to prove the
following local statement: Given $\tau \in {\rm T(C^*(\pi(T)))_{QD}}$,
a finite set ${\mathfrak F} \subset C^*(T)$, a finite set $F \subset
H$ and $\delta > 0$ there exists a finite rank projection $P \in B(H)$
such that:

\begin{enumerate}
\item[(a)] $\| [P, X] \| < \delta$ for all $X \in {\mathfrak F}$.

\item[(b)] $Px = x$ for all $x \in F$.

\item[(c)] $| \mathrm{tr}_{rank(P)} (PXP) - \tau(\pi(X)) | <
\delta$ for all $X \in {\mathfrak F}$.
\end{enumerate}

Assuming that this local statement holds, an argument similar to
the one given in the proof of Proposition \ref{thm:QDcase} would
show that one can construct a filtration $P_1 \leq P_2 \leq
\cdots$ such that $\| [P_n, T] \| \to 0$ and for every $\tau \in
{\rm T(C^*(\pi(T)))_{QD}}$ there exists a subsequence $\{n(k)\}$
such that $\varphi_{n(k)} \to \tau \circ \pi|_{C^*(T)}$ in the
weak-$*$ topology.  This evidently yields part (3) of the theorem.
However it also gives statements (1) and (2). Indeed, letting
$\Phi : C^*(T) \to \Pi P_n B(H) P_n / \oplus P_n B(H) P_n$ be the
u.c.p.\ map which takes $X$ to $P_n X P_n + \oplus P_n B(H) P_n$
we have that in fact $\Phi$ is a faithful $*$-homomorphism.  The
faithfulness follows from the fact that the $P_n$'s are a
filtration and hence $\|X \| = \lim \|P_n X P_n\|$ for every $X
\in B(H)$.  That $\Phi$ is multiplicative follows from the
quasidiagonality assumption since $$\| P_n XY P_n - P_nX P_n Y P_n
\| = \| P_n X\big(YP_n - P_nY\big)P_n \| \leq \|X\|\|[P_n,Y]\| \to
0.$$

From the theorem quoted above it follows that for every $\epsilon
> 0$, $$\sigma^{(\varepsilon)}(\Phi(T)) = \limsup_{n \to \infty}
\sigma^{(\varepsilon)}(P_n T P_n).$$ Since $\sigma(S) = \cap \
\sigma^{(\varepsilon)}(S)$ holds in general and since $\Phi$ is a
faithful homomorphism we conclude that $$\sigma(T) =
\sigma(\Phi(T)) = \cap \ \sigma^{(\varepsilon)}(\Phi(T)) = \cap
\limsup_{n \to \infty} \sigma^{(\varepsilon)}(P_n T P_n).$$ To see
that part (1) also follows we observe that the finite section
method is stable if and only if $(P_n T P_n) + \oplus P_n B(H)
P_n$ is an invertible element of $\Pi P_n B(H) P_n / \oplus P_n
B(H) P_n$ (this is a straightforward exercise -- or see
\cite[Theorem 1.15]{HRS}). Hence stability is equivalent to the
invertibility of $\Phi(T)$ which, by preservation of spectra, is
equivalent to the invertibility of $T$.

The proof of the required local statement is also similar to the
proof of Proposition \ref{thm:QDcase} however there are a few
additional wrinkles to iron out.  First of all, since $C^*(T)$ is
exact (hence locally reflexive), we can apply the Effros-Haagerup
lifting theorem (cf.\ \cite{wassermann:exactbook}) to find a
u.c.p.\ splitting $\Phi : C^* (\pi(T)) \to B(H)$.  It follows that
$C^* (\pi(T))$ is also an exact, quasidiagonal C$^*$-algebra.
Indeed, exactness always passes to quotients but we also get
quasidiagonality from Voiculescu's abstract characterization since
we can compose the splitting $\Phi$ with compression by an
appropriately chosen sequence of finite rank projections (coming
from the quasidiagonality of the operator $T$) to produce some
u.c.p.\ maps to matrix algebras which are asymptotically
multiplicative in norm and asymptotically isometric.

Since $T$ is a quasidiagonal operator, if we are given finite sets
${\mathfrak F} \subset C^* (T)$ and $F \subset H$ we can find a
finite rank projection $Q$ such that $Qx = x$ for all $x \in F$,
$\| [Q, X] \|$ is as small as we like and (hence) $\| X - (QXQ +
Q^{\perp} X Q^{\perp}) \|$ is also as small as desired for every
operator $X \in {\mathfrak F}$. Note that the u.c.p.\ map
$\Phi_{Q^{\perp}} : C^* (\pi(T)) \to B(Q^{\perp}H)$,
$\Phi_{Q^{\perp}} (Y) = Q^{\perp} \Phi(Y) Q^{\perp}$ is still a
u.c.p.\ splitting (i.e.\ a faithful $*$-homomorphism modulo the
compacts) and, moreover, is nearly multiplicative on
$\pi({\mathfrak F})$.  Letting $\rho : C^*(\pi(T)) \to B(K)$ be
any faithful representation which contains no non-trivial compact
operators, we can, thanks to Proposition \ref{thm:QDcase}, find a
filtration $Q_1 \leq Q_2 \leq \cdots$ which asymptotically
commutes (in norm) with $\rho(C^*(\pi(T)))$ and, furthermore,
which recaptures any fixed quasidiagonal trace $\tau \in {\rm
T(C^*(\pi(T)))_{QD}}$.  From Voiculescu's Theorem (version
\ref{thm:technicalVoiculescuThm}) we can find a unitary $U : K \to
Q^{\perp}H$ such that $\| Q^{\perp} X Q^{\perp} - U\rho(\pi(X))U^*
\|$ is very small for all $X \in {\mathfrak F}$.  Defining $P_n =
Q \oplus UQ_n U^*$ we get a filtration of $H$ where each
projection almost commutes with ${\mathfrak F}$ (they are {\em
not} asymptotically commuting, but there is a small, uniform upper
bound on the norms of commutators) and cutting by these
projections will almost recover $\tau$ (since $Q$ is fixed and
finite dimensional, its contribution becomes negligible as $n \to
\infty$).
\end{proof}

\begin{rem} In most cases, part (3) above is a much stronger
statement than part (2).  Indeed, in \cite{arveson:numerical} Arveson
observes how statement (3), which he regards as an analogue of
Szeg\"{o}'s Limit Theorem, can be used to recover the essential
spectrum of a self-adjoint operator.

Note that we only used exactness of C$^*(T)$ to deduce the existence
of a u.c.p.\ splitting C$^*(\pi(T)) \to B(H)$ (which, in turn, was
only needed in part (3) above).  Hence the result above holds under
this weaker assumption.  However, many standard examples in operator
theory are easily seen to generate exact C$^*$-algebras.  For example,
essentially normal operators (hence Toeplitz operators with continuous
symbol) and all weighted shift operators have this property (since
essentially normal operators generate nuclear C$^*$-algebras while
every weighted shift is an element of the nuclear C$^*$-algebra
$l^{\infty}({\mathbb N}) \rtimes {\mathbb N}$ and hence generates
something exact). Many of the concrete operators considered in
\cite{HRS} can be realized as weighted shifts (or, more generally, as
elements in one of the algebras $l^{\infty}({\mathbb N}) \rtimes
{\mathbb N}$ or $l^{\infty}({\mathbb Z}) \rtimes {\mathbb Z}$) and
hence generate exact C$^*$-algebras.  Of course, quasidiagonality of
concrete examples is a much harder question.  However for the classes
of essentially normal and (finite direct sums of) weighted shift
operators there are nice characterizations of quasidiagonality as well
(cf.\ \cite[Corollary IX.7.4, Section IX.8]{davidson}, \cite{smucker},
\cite{narayan}).  Moreover, we should also point out that if a
weighted shift satisfies the appropriate hypotheses (cf.\
\cite{smucker}) then it is possible to actually construct (as opposed
to proving the existence of) an asymptotically commuting filtration
via Berg's technique (this would be enough for statements (1) and (2)
above, but not necessarily enough for conclusion (3)).
\end{rem}

Specializing to the case of a self-adjoint operator Theorem
\ref{thm:finitesection} takes an especially nice form.  Note that if
$T \in B(H)$ is self-adjoint and $P \in B(H)$ is a projection then
$PTP$ is still a self-adjoint operator (in particular, it is normal --
see the remark after Theorem \ref{thm:convergenceofspectra}).  Hence
if $\{ P_n \}$ is a filtration and the rank of $P_n$ is some integer
$d_n$ then $P_n T P_n$ will have $d_n$ (not necessarily distinct)
eigenvalues which we will denote by $\lambda^{(n)}_1, \ldots,
\lambda^{(n)}_{d_n}$.  To complete the spectral approximation picture
we will need one more of Arveson's notions: A real number $\lambda$
will be called {\em essential} (for $T$ with respect to $\{P_n\}$) if for
every open interval $U$ containing $\lambda$ we have $$\lim_{n \to
\infty} N_n (U) = \infty,$$ where $N_n(U)$ is the number of
eigenvalues of $P_n T P_n$ (counted with multiplicities) in the
interval $U$.  We then let $\liminf_{ess} \sigma(P_nTP_n)$ denote the
set of essential points.  Note that $\liminf_{ess}
\sigma(P_nTP_n) \subset \liminf \sigma(P_nTP_n)$.

\begin{cor} Let $T \in B(H)$ be a self-adjoint operator.  Then there
exists a filtration $P_1 \leq P_2 \leq P_3 \cdots$ such that:

\begin{enumerate}
\item The finite section method is stable if and only if $T$ is
invertible.

\item $\sigma(T) = \limsup\limits_{n \to \infty} \sigma(P_n T P_n) =
\liminf\limits_{n \to \infty} \sigma(P_n T P_n)$.

\item $\sigma_{ess}(T) = \liminf_{ess} \sigma(P_n T P_n)$.

\item If $\mu$ is a regular, Borel, probability measure on $\sigma(T)$
then there exists a subsequence $\{n(k)\}$ such that for every $f \in
C(\sigma(T))$ we have $$\lim_{k \to \infty}
\frac{f(\lambda^{(n(k))}_1) + f(\lambda^{(n(k))}_2) + \cdots +
f(\lambda^{(n(k))}_{d_{n(k)}})}{d_{n(k)}} = \int_{\sigma(T)} f(x)
d\mu(x)$$ if and only if $\mu$ is supported on the essential spectrum
of $T$.
\end{enumerate}
\end{cor}

\begin{proof} Let $P_1 \leq P_2 \leq \cdots$ be the filtration given
by Theorem \ref{thm:finitesection}.  According to \cite[Theorem
7.2]{HRS} one always has the inclusion $\sigma(T) \subset \liminf
\sigma(P_n T P_n)$ -- so long as $T$ is self-adjoint and $\{ P_n \}$
is a filtration -- and hence we get $$\sigma(T) \subset \liminf
\sigma(P_n T P_n) \subset \limsup \sigma(P_n T P_n) = \sigma(T),$$
where we have used the fact that $P_n TP_n$ are normal and hence the
$\varepsilon$-psuedospectra in Theorem \ref{thm:finitesection} are
unnecessary (see the remark after Theorem
\ref{thm:convergenceofspectra}).

Conclusion (3) is essentially due to Arveson.  The set-up here is not
exactly the same as that in \cite[Theorem 7.10]{HRS} but the
difference is minor and the details will be left to the reader.

For the final statement above first note that there is a
one-to-one correspondence between regular, Borel, probability
measures on the essential spectrum of T and quasidiagonal traces
on $C^*(\pi(T))$ (since every state on an abelian C$^*$-algebra
is, in fact, a uniform locally finite dimensional trace). The
conclusion now follows from the corresponding statement in Theorem
\ref{thm:finitesection}.
\end{proof}

\begin{rem} In \cite{arveson:numerical} Arveson proved the essential
parts of the corollary above under the following additional
hypotheses: $T$ should be a self-adjoint element of a simple, unital
C$^*$-algebra $A$ which has a unique tracial state and which admits an
`$A$-filtration' (see \cite[Theorems 3.8 and 4.5]{arveson:numerical}).
Of course, at the time Arveson was introducing these notions and the
point of his results was to show that there was an abstract set of
hypotheses which would ensure convergence results (and hence give hope
of attacking infinite dimensional problems numerically).  Moreover, he
showed that a number of interesting examples where covered by his
results.  In the remarks after \cite[Proposition
4.4]{arveson:numerical} Arveson notes that it is natural to ask which
C$^*$-algebras admit an `$A$-filtration' (hence fall under his
results) and points out that things like Cuntz algebras do not (the
existence of an `$A$-filtration' implies the existence of a tracial
state -- which no purely infinite C$^*$-algebra will have).  However,
from the theoretical point of view, it does not follow that
self-adjoint operators in Cuntz algebra are beyond the reach of these
techniques.  Indeed, the corollary above holds for a completely
arbitrary self-adjoint operator and hence there is {\em always} hope
of attacking infinite dimensional self-adjoint operators numerically.
\end{rem}

\section{Amenable traces and obstructions in K-homology}
\label{thm:k-homology}

In this section we will show how the theory of amenable traces can
be used to solve some K-homological questions and a natural
question regarding the passage of quasidiagonality to the Calkin
algebra. The main results of this section are taken from
\cite{anderson}, \cite{connes:compactmetricspaces},
\cite{wassermann:annals} and \cite{wassermann:nonQDstablyfinite}.
However, the invariant mean point of view will turn these
``technical'' results into very natural and intuitively clear
theorems.

The three questions which are easily solved using amenable traces
are the following:
\begin{enumerate}
\item (cf.\ \cite{connes:compactmetricspaces}) If $\Gamma$ is a discrete
group can one always construct a finitely summable, unbounded Fredholm
module over $C^*_r (\Gamma)$?

\item (cf.\ \cite{anderson}, \cite{wassermann:annals},
\cite{wassermann:nonQDstablyfinite}) Is the BDF semigroup
Ext($\cdot$) always a group?\footnote{See \cite{davidson} for a
nice introduction to Ext and BDF theory.  In particular we refer
to the same text for the definitions and basic facts regarding the
extension semigroup.  All we will need in this paper is the
following result of Arveson \cite{arveson:extensions}: For a given
C$^*$-algebra $A$, Ext($A$) is {\em not} a group if and only if
there exists a $*$-homomorphism $\pi:A \to Q(H)$ to the Calkin
algebra with the property that there is {\em no} u.c.p.\ lifting
into $B(H)$.}

\item (cf.\ \cite{wassermann:annals},
\cite{wassermann:nonQDstablyfinite}) If $A \subset B(H)$ is a
quasidiagonal set of operators is it true that the image of $A$ in the
Calkin algebra is quasidiagonal?
\end{enumerate}

The answer to all the questions above is `No'.  It is intriguing,
however, that there is a single reason why all three questions are
false: {\em Not all traces are amenable.}

This obstruction is explicitly pointed out in
\cite{connes:compactmetricspaces} and implicitly used in
\cite{anderson}, \cite{wassermann:annals} and
\cite{wassermann:nonQDstablyfinite}.  Indeed we can reformulate
\cite[Lemma 9]{connes:compactmetricspaces} as follows.

\begin{lem}(See \cite[Lemma 9, Remark 10(b)]{connes:compactmetricspaces})
Let $A$ be a C$^*$-algebra and assume that $A$ admits a finitely
summable, unbounded Fredholm module.  Then $\TAim \neq \emptyset$.
\end{lem}

Since the {\em reduced} group C$^*$-algebra of any discrete,
non-amenable group has no amenable traces (cf.\ Proposition
\ref{thm:reducedgroupalgebras}) we immediately deduce the
following result.

\begin{cor}(cf.\ \cite[Theorem 19]{connes:compactmetricspaces})
Let $\Gamma$ be any discrete, non-amenable group.  Then $C^*_r
(\Gamma)$ has no unbounded, finitely summable Fredholm modules.
\end{cor}

\begin{rem}
In \cite{dvv:quasicentralapproximateunits} Voiculescu defined a
notion of `subexponential', unbounded Fredholm module and
generalized Connes' result to this setting (i.e.\ the reduced
group C$^*$-algebra of a non-amenable, discrete group never admits
one of these either -- see \cite[Proposition
4.10]{dvv:quasicentralapproximateunits}).  Amenable traces explain
this result in exactly the same way: If a C$^*$-algebra $A$ admits
a subexponential, unbounded Fredholm module then $\TAim \neq
\emptyset$ (cf.\ \cite[Proposition
4.6]{dvv:quasicentralapproximateunits}) and hence $C^*_r (\Gamma)$
has no such modules whenever $\Gamma$ is not amenable.
\end{rem}

The Ext and quasidiagonality questions above are easily handled
simultaneously.  We believe that the ideas involved become most
transparent when treated in generality and so this is what we will
do. The main idea is to construct a C$^*$-algebra which has a
tracial state which is not amenable.  The easiest example of such
an algebra is the {\em reduced} group C$^*$-algebra of any
discrete, non-amenable group.  Unfortunately these algebras don't
quite do the trick, but certain extensions of them
will.\footnote{Actually, in \cite{haagerup-thorbjornsen:II}
Haagerup and Thorbj{\o}rnsen have recently shown that
$C^*_r({\mathbb F}_2)$ will do the trick! They have shown the
remarkable fact that $C^*_r({\mathbb F}_2) \subset \Pi M_{k(n)}/
\oplus M_{k(n)}$ and hence Lemma \ref{thm:HT} applies.}

We begin with some general remarks which underlie the main ideas.
Let $M_{k(n)}({\mathbb C})$ be a sequence of matrix algebras and
$$\Pi M_{k(n)} = \{ (x_n) : \sup_n \| x_n \| < \infty \}$$ denote
the $l^{\infty}$-product.  The $c_0$-product, $$\oplus M_{k(n)} =
\{ (x_n) : \| x_n \| \to 0, n \to \infty \},$$ sits as an ideal in
$\Pi M_{k(n)}$ and we will identify the quotient C$^*$-algebra,
$$\Pi M_{k(n)}/ \oplus M_{k(n)},$$ with a subalgebra of the Calkin
algebra (associated to the Hilbert space $H = \oplus {\mathbb
C}^{k(n)}$). Note that if $\omega \in \beta {\mathbb N} \backslash
{\mathbb N}$ is a free ultrafilter then we can define a trace
$\tau_{\omega}$ on $\Pi M_{k(n)}/ \oplus M_{k(n)}$ by
$$\tau_{\omega} ((x_n) + \oplus M_{k(n)}) = \lim_{n \to \omega}
\mathrm{tr}_{k(n)} (x_n).$$ The following lemma explains why Ext
need not be a group and why quasidiagonality need not pass to the
Calkin algebra.

\begin{lem}
\label{thm:HT} Let $A \subset \Pi M_{k(n)}/ \oplus M_{k(n)}$ be
given and assume that $\tau_{\omega}|_A$ is not an amenable trace
on $A$. Then Ext$(A)$ is not a group.  If $A$ has no amenable
traces whatsoever then $A$ is not quasidiagonal.
\end{lem}

\begin{proof} Assume that Ext$(A)$ is a group.  Then there exists a
u.c.p.\ map $\Phi:A \to B(H)$ which is a splitting for the
inclusion $A \hookrightarrow \Pi M_{k(n)}/ \oplus M_{k(n)} \subset
Q(H)$ (cf.\ \cite{arveson:extensions}).  Recall that $H = \oplus
{\mathbb C}^{k(n)}$ and we have a natural inclusion $\Pi M_{k(n)}
\subset B(H)$ and a natural conditional expectation $E:B(H) \to
\Pi M_{k(n)}$ given by $$E(T) = (P_n T P_n)_{n \in \mathbb{N}}$$
where $P_n$ denotes the unit of $M_{k(n)}$.  This conditional
expectation evidently has the property that it maps compact
operators into compact operators (i.e.\ into $\oplus M_{k(n)}$)
and this implies that $E\circ\Phi: A \to \Pi M_{k(n)}$ is also a
splitting.  Indeed, since each $a \in A$ has a lift, say
$\tilde{a}$, in $\Pi M_{k(n)}$ we have that $\Phi(a) - \tilde{a}$
is compact and hence $$E(\Phi(a) - \tilde{a}) = E\circ\Phi(a) -
\tilde{a} \in \oplus M_{k(n)}.$$  This implies that
$E\circ\Phi(a)$ is also a lift of $a$ for each $a \in A$.

To ease notation, we will now forget about $E$ and just let $\Phi
:A \to \Pi M_{k(n)}$ be a u.c.p.\ map which lifts the natural
inclusion $A \subset Q(H)$.  For each $n\in \mathbb{N}$ we now
define a u.c.p.\ map $\phi_n:A \to M_{k(n)}$ by $$\phi_n(a) = P_n
\Phi(a) P_n.$$  Since $\Phi$ is a splitting we have that $\Phi(ab)
- \Phi(a)\Phi(b) \in \oplus M_{k(n)}$ and this implies that the
maps $\phi_n$ are asymptotically multiplicative in norm.  Also, it
is not hard to see that $\tau_{\omega}|_A$ is in the weak-$*$
closure of the states $\{ \mathrm{tr}_{k(n)} \circ \phi_n \}
\subset S(A)$ and this implies that $\tau_{\omega}|_A$ is a
quasidiagonal trace on $A$ (hence amenable) which contradicts our
assumption.

The second statement above is trivial since every (unital)
quasidiagonal C$^*$-algebra has at least one quasidiagonal tracial
state.
\end{proof}

Thus our strategy is clear -- construct some algebra $A \subset
\Pi M_{k(n)}/ \oplus M_{k(n)}$ such that $\tau_{\omega}|_A$ is not
an amenable trace on $A$ or, better yet, so that $A$ has no
amenable traces at all. There are two ways of using discrete
groups to realize this strategy.  One starts from a discrete,
non-amenable, residually finite group and the other starts from a
discrete, property T group with infinitely many non-equivalent,
finite dimensional, irreducible, unitary representations.  We
treat the former case first.

So let $\Gamma$ be a non-amenable, residually finite group (e.g.\ a
free group) and $\Gamma \vartriangleright \Gamma_1 \vartriangleright
\Gamma_2 \vartriangleright \cdots$ be a descending sequence of normal
 subgroups of finite index  such that the
intersection of the $\Gamma_n$'s is the neutral element.  Let
$\pi_n : C^*(\Gamma) \to B(l^2(\Gamma/ \Gamma_n)) \cong M_{k(n)}$
be the representations induced by the left regular representations
of the $\Gamma / \Gamma_n$'s.  Define $$\pi = \oplus \pi_n :
C^*(\Gamma) \to \Pi M_{k(n)}.$$ Now consider the algebra $$B =
\sigma \circ \pi (C^*(\Gamma)) \subset \Pi M_{k(n)}/ \oplus
M_{k(n)},$$ where $\sigma : \Pi M_{k(n)} \to \Pi M_{k(n)}/ \oplus
M_{k(n)}$ is the quotient map. Since $B$ sits where we want it to
(by construction) the only question is whether or not $B$ has any
amenable traces.  Unfortunately, it seems possible that $B$ could
have such traces.  However, the following observation will keep
our hopes alive.

\begin{lem}
\label{thm:whatever} The image of $B$ in the GNS representation of
$\Pi M_{k(n)}/ \oplus M_{k(n)}$ with respect to the trace
$\tau_{\omega}$ is isomorphic to $C^*_r (\Gamma)$ (which has no
amenable traces).
\end{lem}

\begin{proof}
Evidently we may identify the GNS representation of $\Pi M_{k(n)}/
\oplus M_{k(n)}$, with respect to $\tau_{\omega}$, with a subalgebra
(in fact, II$_1$-subfactor) of $R^{\omega}$. Using this
identification, it is clear that the resulting trace on $C^*(\Gamma)$
(which is just $\tau_{\omega} \circ \sigma \circ \pi$) is the
canonical trace which vanishes on all non-trivial group elements.  By
uniqueness of GNS representations, the lemma follows.
\end{proof}

Thus we only need a result which states that amenable traces pass
to quotients (for then $B$ can't have any amenable traces since
one of its quotients has no such traces).  Unfortunately this is
not true (e.g.\ the canonical trace on $C^*({\mathbb F}_2)$ is
amenable while it is not when regarded as a trace on $C^*_r
({\mathbb F}_2)$).  Thus we will have to use a bit of trickery.
The next observation gives us a way of forcing amenable traces to
pass to quotients.

\begin{lem}
Let $X$ be a C$^*$-algebra, $J \subset X$ be an ideal, $J^{\perp}
\subset X$ be the ideal perpendicular to $J$ (i.e.\ $J^{\perp} =
\{ x \in X: xJ=Jx=0 \}$) and $\tau$ be an amenable trace on $X$.
If $\tau|_{J^{\perp}} \neq 0$ then $X/J$ has at least one amenable
trace.
\end{lem}

\begin{proof} Let $\gamma = \frac{1}{\| \tau|_{J^{\perp}} \|}
\tau|_{J^{\perp}}$ and Proposition \ref{thm:restrictiontoideal}
ensures that $\gamma$ is an amenable trace on $J^{\perp}$. Since
$J \cap J^{\perp} = \{ 0 \}$ we may identify $J^{\perp}$ with an
ideal in $X/J$. However, thanks to  Proposition
\ref{thm:extensionfromideal} we may extend the amenable trace
 $\gamma$ to an amenable trace on $X/J$.
\end{proof}

We almost have all the pieces to the puzzle.  We just have to find an
appropriate ideal to add on to the algebra $B$ above.  To do this
requires a bit more notation.  Let $I_{\omega} \subset \Pi M_{k(n)}$
be the ideal of sequences $(x_n)$ such that $\| x_n \|_2 \to 0$ as $n
\to \omega$.  Let $J = B \cap \sigma(I_{\omega})$.  Note that $B/J
\cong C^*_r (\Gamma)$ by Lemma \ref{thm:whatever}.

\begin{lem}(cf.\ \cite{anderson})
There exists a projection $p \in \Pi M_{k(n)}/ \oplus M_{k(n)}$ and a
unitary $u \in \Pi M_{k(n)}/ \oplus M_{k(n)}$ such that:

\begin{enumerate}
\item $pJ = Jp = 0$.

\item $\tau_{\omega} (p) \geq 1/2$.

\item $p + upu^* \geq 1$.
\end{enumerate}
\end{lem}

Please see \cite[Proposition, pg.\ 456]{anderson} for the proof of
this elementary (but crucial!) fact. (Actually, Anderson only
constructs the projection and in \cite{wassermann:nonQDstablyfinite}
Wassermann makes the observation that the unitary $u$ exists.)  We now
define the right algebras.  Put $$A = C^*(B, p) \subset W = C^*(B, p,
u).$$

\begin{cor}(cf.\ \cite{anderson}, \cite{wassermann:nonQDstablyfinite})
$Ext(A)$ and $Ext(W)$ are not groups.  $W$ is not quasidiagonal,
though it is the image (in the Calkin algebra) of a quasidiagonal
set of operators on $B(H)$.
\end{cor}

\begin{proof}
First note that $J$ is also an ideal in $A$ (since $p$ is
orthogonal to $J$) and $A/J$ contains a unital copy of
$C^*_r(\Gamma) \cong B/J$. By the sequence of lemmas above we have
the following chain of implications: $Ext(A)$ is a group
$\Longrightarrow \tau_{\omega}|_A$ is amenable $\Longrightarrow
A/J$ has an amenable trace (since $\tau_{\omega}$ restricted to
the ideal generated by $p$ -- which is orthogonal to $J$ -- does
not vanish).  However the last statement provides a contradiction
since $C^*_r(\Gamma) \subset A/J$ has no amenable traces.

More generally, note that this argument shows that any amenable
trace on $A$ must vanish on $p$.  This implies that $W$ has no
amenable traces whatsoever since the equation $p + upu^* \geq 1$
implies that no trace on $W$ can vanish on $p$.  Hence Ext$(W)$ is
not a group and $W$ is not quasidiagonal.
\end{proof}

\begin{rem}
In the case that $\Gamma = {\mathbb F}_2$, the algebra $A$ is
(essentially) the example given in \cite{anderson} and the algebra $W$
is the one given in \cite{wassermann:nonQDstablyfinite}.
\end{rem}

We now turn to the property T case.  We will need the following
structure theorem for full group C$^*$-algebras of property T groups.

\begin{thm}(cf.\ \cite[3.7.6]{higson-roe})
Let $\Gamma$ be a discrete group with Kazhdan's property T.  Then for
each finite dimensional, irreducible, unitary representation $\pi :
\Gamma \to B(H_{\pi})$ there is a central projection $p_{\pi} \in
C^*(\Gamma)$ with the following properties:
\begin{enumerate}
\item $p_{\pi_1} = p_{\pi_2}$ if and only if $\pi_1$ and $\pi_2$ are
unitarily equivalent and $p_{\pi_1} \perp p_{\pi_2}$ otherwise.

\item $p_{\pi} C^*(\Gamma) = p_{\pi} C^*(\Gamma) p_{\pi} \cong
B(H_{\pi}) = M_{n}$, where $n = dim(\pi)$.

\item If $\rho\colon C^*(\Gamma) \to B(K)$ is a representation which
contains the finite dimensional, irreducible representation $\pi$ as a
subrepresentation then $\rho(p_{\pi}) \neq 0$ (actually,
$\rho(p_{\pi})$ acts as the orthogonal projection onto the
$\pi$-isotypical subspace).
\end{enumerate}
\end{thm}

With this result in hand, it is an easy matter to construct a
C$^*$-algebra with traces which are not amenable.  Indeed, let $J
\subset C^*(\Gamma)$ be the ideal generated by all the (central)
projections coming from finite dimensional, irreducible
representations.  Note that $J \cong \oplus M_{k(n)}$ (possibly a
finite direct sum) for some integers $k(n)$.  Since the multiplier
algebra of $\oplus M_{k(n)}$ is just $\Pi M_{k(n)}$ we get a
$*$-homomorphism $\pi : C^* (\Gamma) \to \Pi M_{k(n)}$ such that
$\pi(J) = \oplus M_{k(n)}$.  Identifying $\Pi M_{k(n)}$ with the
block diagonal operators on $H = \oplus {\mathbb C}^{k(n)}$ and
letting $\sigma : B(H) \to Q(H)$ be the quotient map onto the
Calkin algebra we define the right C$^*$-algebra as: $$W_{\Gamma}
= \sigma \circ \pi (C^*(\Gamma)) \subset \Pi M_{k(n)}/ \oplus
M_{k(n)}.$$

\begin{cor}(cf.\ \cite{wassermann:annals})
If $\Gamma$ is a discrete group with Kazhdan's property T and $\Gamma$
has infinitely many non-equivalent, finite dimensional, irreducible
representations (e.g.\ $SL(3,{\mathbb Z})$) then $Ext(W_{\Gamma})$ is
not a group and $W_{\Gamma}$ is not quasidiagonal.
\end{cor}

\begin{proof} It suffices to show that $W_{\Gamma}$ has no amenable
traces.  So assume that $\tau$ is amenable on $W_{\Gamma}$.
Representing $W_{\Gamma}$ faithfully on some Hilbert space $K$
(i.e.\ $W_{\Gamma} \subset B(K)$) we can, by part (3) of Theorem
\ref{thm:amenabletracesII}, find unit vectors $v_n \in HS(K)$
(where $HS(K)$ is the Hilbert space of Hilbert-Schmidt operators
on K) which are asymptotically invariant under the (conjugation)
action of the unitary group of $W_{\Gamma}$ on
$HS(K)$.\footnote{Since the Hilbert-Schmidt norm is invariant
under left and right multiplication by unitaries, Connes'
F{\o}lner condition is equivalent to $\| u v_n u^* - v_n \|_{HS}
\to 0$, where $v_n = \frac{1}{\|P_n\|_{HS}} P_n$, whenever $u$ is
unitary.} Since $\Gamma$ has property $T$ it follows that this
conjugation action on $HS(K)$ must have a fixed point.  But this
is precisely the same thing as saying that there is a
Hilbert-Schmidt operator in the commutant $W_{\Gamma}^{\prime}$
and once you get a compact in the commutant you must also have a
finite rank (spectral) projection in the commutant of $W_{\Gamma}$
as well.  But commuting finite rank projections give finite
dimensional representations and this will give our contradiction.
Indeed, we constructed $W_{\Gamma}$ by first factoring out all of
the irreducible finite dimensional representations from
$C^*(\Gamma)$ (which implies that we factored out all finite
dimensional representations as well).  Hence, $W_{\Gamma}$ can't
have any finite dimensional representations since these would
induce finite dimensional representations on $C^*(\Gamma)$ which
factorize through $W_{\Gamma}$.
\end{proof}

\begin{rem}
In \cite{kirchberg:invent} Kirchberg proved that a C$^*$-algebra
$A$ has the local lifting property if and only if Ext$(S(A))$ is a
group, where $S(A)$ denotes the unitization of $C_0 ({\mathbb R})
\otimes A$. Ozawa has pointed out to us that the class of
C$^*$-algebras with the local lifting property is, in some sense,
small and hence Kirchberg's result shows that Ext is not a group
for lots of C$^*$-algebras. (For example, `most' finite
dimensional operator spaces will generate C$^*$-algebras without
the local lifting property since the space of n-dimensional
operator spaces is not separable -- cf.\ \cite{junge-pisier} --
while every finite dimensional subspace of a C$^*$-algebra with
the local lifting property can be identified with a subspace of
$C^* ({\mathbb F}_{\infty})$ and the space of n-dimensional
operator subspaces of any separable C$^*$-algebra is a separable
set.)

However, it is interesting to note that all of the {\em concrete}
examples of C$^*$-algebras without the local lifting property fail
to have this property because of invariant mean considerations (we
are referring to concrete {\em separable} examples -- $B(H)$ does
not have the LLP, by \cite{junge-pisier}).  Indeed, the only
examples of C$^*$-algebras (that we are aware of) without the
local lifting property are the reduced group C$^*$-algebras of
non-amenable, residually finite discrete groups or groups which
contain a subgroup isomorphic to such a group. Actually the most
general statement we know is: If $\Gamma$ contains a non-amenable
subgroup $H$ such that $C^*_r (H)$ embeds into $R^{\omega}$ then
$C^*_r(\Gamma)$ does not have the LLP and hence
Ext$(S(C^*_r(\Gamma)))$ is not a group. (Proof: Since $C^*_r (H)$
has no amenable traces, yet it does embed into $R^{\omega}$ it
follows that $C^*_r (H)$ does not have the LLP. Since there is a
conditional expectation $C^*_r(\Gamma) \to C^*_r (H)$ it follows
that $C^*_r(\Gamma)$ does not have the LLP either -- cf.\
\cite[Corollary 2.6.(v) and Proposition 3.1]{kirchberg:invent}.)
\end{rem}

\section{Stable finiteness versus quasidiagonality}
\label{thm:stablyfinitevsQD}

We now observe that approximation properties of traces are related
to some important questions about quasidiagonal C$^*$-algebras. At
the moment there are only two known obstructions to
quasidiagonality; (1) the existence of an infinite projection (in
some matrix algebra over the given algebra) and (2) the absence of
a quasidiagonal trace. In other words, every quasidiagonal
C$^*$-algebra is stably finite and has a quasidiagonal tracial
state (recall that the latter statement is not true in the
non-unital case).  In \cite{dvv:QDsurvey} Voiculescu posed the
challenging problem of finding a complete set of obstructions --
preferably of a topological nature -- to quasidiagonality.  This
ambitious goal is still out of our reach as the class of
quasidiagonal C$^*$-algebras is deceptively difficult to get a
handle on. In this section we wish to describe two of the major
open questions around quasidiagonality and point out the role that
approximation properties of traces have already played as well as
the role they may play in the future.

The two questions we will address are:

\begin{enumerate}
\item Is every nuclear, stably finite C$^*$-algebra necessarily
quasidiagonal?

\item Is the hyperfinite II$_1$-factor quasidiagonal?
\end{enumerate}

The first question was asked by Blackadar and Kirchberg in
\cite{blackadar-kirchberg}.  This question is of fundamental
importance in Elliott's classification program.  Indeed, suppose
that one could even show the weaker statement that every unital,
simple, nuclear, stably finite C$^*$-algebra with real rank zero
was quasidiagonal.  Then by Popa's result it would follow that
every such algebra was a Popa algebra.  Hence every such algebra
would satisfy an internal finite dimensional approximation
property which is similar to that which defines Lin's tracially AF
algebras.  In other words, question (1) above is `half way' to
completing the simple, stably finite, real rank zero case of
Elliott's conjecture (the other half being a proof that nuclear,
Popa algebras with real rank zero are tracially AF and hence fall
under Lin's classification theorem). Clearly weaker variations of
question (1) would already be of great interest.  For example,
what if one assumes a faithful trace?  How about just the simple,
real rank zero case?

Note also that question (1), or even the weaker version where one
assumes a faithful trace, would immediately imply Rosenberg's
conjecture that the reduced group C$^*$-algebra of any discrete,
amenable group is quasidiagonal.  The second question above would
also imply Rosenberg's conjecture as well as yield the faithful
trace case of question (1) since every finite, hyperfinite von
Neumann algebra can be embed into $R$. (Note that it seems
unlikely that question (2) would imply question (1) in general
since there are nuclear, stably finite -- even quasidiagonal --
C$^*$-algebras which can't be embed into $R$ since they have no
{\em faithful} trace.)

In our opinion, it is far from clear which way either of the
questions above will go.  Suppose one wanted to give negative
answers.  One natural strategy for doing this would be to embed a
C$^*$-algebra which was not quasidiagonal into $R$ or into some
nuclear, stably finite C$^*$-algebra (since quasidiagonality
passes to subalgebras, this would give a contradiction).  However,
amenable trace considerations show that we cannot realize this
strategy with our current set of examples.

\begin{prop}
\label{thm:BFDIII}
Let $A$ be any C$^*$-algebra which is known to
not be quasidiagonal.  Then it is impossible to embed $A$ into $R$
or into any nuclear, stably finite C$^*$-algebra.
\end{prop}

\begin{proof} We caution the reader that there are plenty of stably
finite C$^*$-algebras for which it is not known whether  they are
quasidiagonal (cf.\ Rosenberg's Conjecture), but we now describe
all the examples which are known to not be quasidiagonal.

The basic examples of non-quasidiagonal C$^*$-algebras are: (1)
Anything with an infinite projection (e.g.\ Cuntz algebras), (2)
$C_r^* (\Gamma)$, where $\Gamma$ is a discrete, non-amenable group
and (3) the algebras constructed by Wassermann and denoted by $W$
and $W_{\Gamma}$ in the previous section of these notes.  Since
quasidiagonality passes to subalgebras, anything containing one of
the examples above will also not be quasidiagonal (e.g.\ most
reduced free products are not quasidiagonal since these tend to
either be purely infinite or contain a copy of $C^*_r ({\mathbb
F}_2)$). However, since the proposition at hand regards embeddings
it will be sufficient to observe that none of the three basic
examples above can be embed into $R$ or into a nuclear, stably
finite C$^*$-algebra.

Excluding case (1) is trivial since stably finite C$^*$-algebras
never contain infinite projections.  Excluding the reduced group
C$^*$-algebras of non-amenable, discrete groups has already been
observed in Corollary \ref{thm:nohomomorphisms2} (recall the
proof: $C_r^* (\Gamma)$ has no amenable traces, while $R$ and any
stably finite, nuclear C$^*$-algebra does have amenable traces).
Finally, Wassermann's examples also can't be embed into anything
with an amenable trace because it was precisely the absence of
such traces which allowed us to deduce that they were not
quasidiagonal C$^*$-algebras.
\end{proof}

\begin{rem} We wish to emphasize that the result above does not imply
that one could not answer the questions above negatively by
embedding a non-quasidiagonal C$^*$-algebra into one of the
algebras in question.  It only says that we can't use any of the
examples which are currently known to not be quasidiagonal.  Of
course, this takes us back to Voiculescu's original problem: what
is a complete set of obstructions to quasidiagonality?
Constructing new examples of non-quasidiagonal algebras would
presumably require gaining insight to Voiculescu's question.
\end{rem}

We should mention that several experts we have talked to feel that
it is very unlikely that $R$ is quasidiagonal.  Their intuition
seems to be based on the following fact which is known to many,
though we were unable to find it written down anywhere.  We thank
George Elliott for showing us a very nice proof.  We have only
modified the last few lines of Elliott's argument so that we can
deduce a slightly stronger statement.

\begin{lem}
Let $R$ act on $L^2(R)$ via the GNS construction.  There is no {\em
 sequence} of nonzero, finite rank projections $P_1, P_2, \ldots$ such
 that $\| [x,P_n] \| \to 0$ for all $x \in R$.
\end{lem}

\begin{proof}  The proof goes by contradiction.  So let $P_1, P_2,
\ldots$ be finite rank projections such that $\| [x,P_n] \| \to 0$
for all $x \in R$.  Put $K = \oplus_{n \in \mathbb N} L^2(R) =
L^2(R) \otimes_2 l^2({\mathbb N})$ and $P = \oplus_{n \in \mathbb
N} P_n$. Then $(x\otimes 1)P - P(x\otimes 1)$ is a compact
operator for every $x \in R$.  Hence, down in the Calkin algebra
$P$ will land in the commutant of $R\otimes 1$.  But then by a
theorem of Johnson and Parrott (see the remarks after \cite[Lemma
3.3]{johnson-parret}) it follows that $P$ is a compact
perturbation of an element in the commutant of $R\otimes 1$. That
is, there exists an infinite matrix $T = (T_{i,j})_{i,j \in
{\mathbb N}}$ such that each $T_{i,j} \in R^{\prime} \subset
B(L^2(R))$ and $P - T$ is compact on $K$.  In particular, this
implies that $\| P_n - T_{n,n} \| \to 0$.  Thus $\| T_{n,n} \| \to
1$ and down in the Calkin algebra the norm of $T_{n,n}$ is tending
to zero.  However this is a contradiction since the commutant of
$R$ is a II$_1$-factor (isomorphic to $R$) and hence a simple
C$^*$-algebra.  Thus the mapping to the Calkin algebra is
isometric.
\end{proof}

Note that the proof above never used the fact that the $P_n$'s are
projections and hence also holds for {\em sequences} of finite rank
operators whose norms are tending to one.  However, since we can
always construct a quasicentral {\em net} of finite rank operators for
$R$ we are left to conclude that the lemma above has more to do with
sequences versus nets (i.e.\ separable versus non-separable Hilbert
spaces) than it does with quasidiagonality.

The point of the above discussion is that providing negative
answers to either of the questions posed seems difficult at this
point. The theory of amenable traces tell us that the embedding
strategy is (currently!) hopeless, while the W$^*$-obstruction
appears to be less related to quasidiagonality issues than had
been expected. As mentioned at the beginning, the only other
obvious strategy is to prove that there exists a stably finite,
nuclear C$^*$-algebra with no quasidiagonal traces as this would
certainly imply that such an example was not a quasidiagonal
C$^*$-algebra (recall, however, that all traces on nuclear
C$^*$-algebras are uniform amenable) or to show that the unique
trace on $R$ is not quasidiagonal. However we are not aware of a
single example of an amenable trace which is not quasidiagonal. We
believe such things exist, but constructing one seems quite
difficult (note the parallel with Voiculescu's question). Our next
result shows that uniform amenable traces may help unlock the
quasidiagonality question for $R$.

\begin{prop}
$R$ is quasidiagonal if and only if for every (separable)
C$^*$-algebra $A$ we have $\TAQD \supset \TAuim$.
\end{prop}

\begin{proof}
We begin with the necessity. Let $A$ be arbitrary.  It suffices to
show that the extreme points of $\TAuim$ belong to $\TAQD$.
However, every extreme point of $\TAuim$ is also an extreme point
of $\TA$, since $\TAuim$ is a face, and hence gives $R$ in the GNS
representation.  Thus, if we assume that $R$ is quasidiagonal then
its unique trace must belong to $\mathrm{T}(R)_{\rm QD}$ and this
completes the proof. (One may worry about non-separability issues
here, but everything works fine with nets.)

For the sufficiency, we first point out that $R$ is quasidiagonal
if and only if all of its separable C$^*$-subalgebras are
quasidiagonal.  So let $A \subset R$ be an arbitrary separable,
unital subalgebra.  Let $\tau \in \TA$ be the restriction of the
unique trace on $R$ to $A$.  Clearly $\tau$ is faithful and
belongs to $\TAuim$.  Hence it also belongs to $\TAQD$.  But we
already saw in the proof of Proposition \ref{thm:rosenberg} that
the existence of a faithful quasidiagonal trace implies
C$^*$-quasidiagonality of the algebra and hence the proof is
complete.
\end{proof}

Applying the proposition above and Theorem
\ref{thm:locallyreflexive} we immediately get the following
corollary.  Note, in particular, that if $R$ is quasidiagonal then
if would follow that every trace on a nuclear C$^*$-algebra is
uniform quasidiagonal (since we already know every trace on such
an algebra is amenable).

\begin{cor}
If $R$ is quasidiagonal then for every locally reflexive
C$^*$-algebra $A$ we have $$\TAim = \TAuim = \TAQD = \TAuQD.$$
\end{cor}

\section{Questions}
\label{thm:questions}

We wish to end these notes with a series of questions.  Preliminary 
versions of this paper appeared in 2001 and hence some of the 
questions have been answered.  However, we list all of our original 
questions just for completeness.

\begin{enumerate}
\item Is every II$_1$-factor representation of a Popa algebra
McDuff?\footnote{This question has been resolved in joint work with
Ken Dykema: $L({\mathbb F}_2)$, which is not McDuff, contains a weakly
dense Popa algebra \cite{BD}.}

\item Is there an example such that $\TAim \neq \TAQD$? How about $\TAuim
\neq \TAuQD$?  The only obvious obstruction is related to
quasidiagonality since the existence of a {\em faithful} trace in
$\TAQD$ implies quasidiagonality.  However we just saw in the
previous section that $\TAim = \emptyset$ for any of the
standard examples of stably finite, non-quasidiagonal
C$^*$-algebras.  For example, can one construct a C$^*$-algebra
which has the WEP and a {\em faithful} trace but which is not
quasidiagonal? The most natural candidate is the hyperfinite
II$_1$-factor.

\item Can a free group or property T II$_1$-factor contain a weakly
dense Popa algebra?\footnote{Free group factors -- at least on
finitely many generators -- do contain weakly dense Popa algebras
\cite{BD}, but the $L({\mathbb F}_{\infty})$ and property T cases are
still open.}

\item Can one give estimates of the free entropy dimension of a finite
set of elements in a Popa algebra which is independent of the
particular trace?  This is related to the semicontinuity/invariance
problem for free entropy dimension.  Constructing a counterexample may
be more reasonable and this would also be of interest.

\item Can one prove a classification theorem for simple, nuclear, real
rank zero C$^*$-algebras which satisfy the UCT and such that $\TA =
\TAQD$ ($= \TAuQD$, by Theorem \ref{thm:locallyreflexive})? In
\cite[Theorem 3.3]{popa:simpleQD} Popa proves that such algebras have
an internal finite dimensional approximation property which should be
of use.  Presumably the role of nuclearity needs to be clarified as
Popa never assumes nuclearity in \cite{popa:simpleQD}.

\item  If $A$ is locally reflexive (or even nuclear) do we always
have that $$\TAuQD \cap \TAlfd = \TAulfd?$$  Feel free to add on
as many additional assumptions, such as real rank zero or stable
rank one, as are necessary! For example, if $A$ is a nuclear Popa
algebra with real rank zero, stable rank one, unperforated
K-theory and unique trace $\tau$, is it necessarily true that
$\tau$ is  uniform locally finite dimensional? (In this setting we
do know that $\tau \in \TAuQD \cap \TAlfd$.) See Proposition
\ref{thm:BFD} for related questions.

\item Can an infinite, simple, discrete group with Kazhdan's
property T be embed into the unitary group of an
$R^{\omega}$-embeddable factor? (Compare with \cite{robertson}
where it is shown that no such embedding exists into the unitary
group of $L({\mathbb F}_n)$ or, more generally, $L(\Gamma)$ for
any a-T-menable discrete group $\Gamma$.)

\item Let $\Gamma$ be a discrete group such that the group von
Neumann algebra of $\Gamma$ has the weak expectation property
relative to a weakly dense C$^*$-subalgebra (i.e.\ $L(\Gamma)
\subset R^{\omega}$). Does it follow that $\Gamma$ is an exact
group? Perhaps just uniformly embeddable into Hilbert space?  An
affirmative answer would have three important consequences: (a)
the Novikov conjecture is true for all residually finite groups
(and every other group which embeds into $R^{\omega}$) (b)
counterexamples to Connes' embedding problem (hence a negative
answer to Voiculescu's `unification problem') since Gromov has
constructed discrete groups which can't be uniformly embed into
Hilbert space and (3) there exist hyperbolic groups which are not
residually finite (since Gromov's examples are apparently
inductive limits of hyperbolic groups and it is easy to see that
inductive limits of groups which embed into $R^{\omega}$ are also
embeddable into $R^{\omega}$). For some time we felt confident
that an affirmative answer was not far off.  However, we have been
unable to fill some gaps in our original strategy and the
statement ``$R^{\omega}$-embeddable implies uniformly embeddable''
has been downgraded to a conjecture.  Recently Guentner, Higson
and Weinberger showed that every linear group is exact and, hence,
uniformly embeddable into Hilbert space (cf.\ \cite{GHW}). Is
there an extremely clever choice of microstates together with an
ultraproduct argument which will work here?

\item Is the canonical trace on $C^*(\Gamma)$ (giving the left
regular representation) always an amenable trace whenever $\Gamma$
is a hyperbolic group? Is every hyperbolic group embeddable into
the unitary group of $R^{\omega}$?  An affirmative answer would
imply that our previous question has a negative answer (since
Gromov's non-exact groups are inductive limits of hyperbolic
groups and hence would have to embed into $R^{\omega}$) while a
negative answer would imply that not every hyperbolic group is
residually finite (since the canonical trace on $C^*(\Gamma)$ is
amenable for every residually finite group).
\end{enumerate}

\bibliographystyle{amsplain}

\providecommand{\bysame}{\leavevmode\hbox to3em{\hrulefill}\thinspace}

\end{document}